\documentclass[12pt,leqno]{article}
\usepackage{amssymb}
\usepackage[mathscr]{eucal}
\usepackage{amsmath,amssymb,latexsym,theorem,bbm}
\usepackage{amsmath,amssymb,latexsym,theorem}
\usepackage{color,url}
\usepackage{appendix}

\newcommand{\CC}{\mathbb{C}}

\newcommand{\NN}{\mathbb{N}}

\newcommand{\RR}{\mathbb{R}}

\newcommand{\ZZ}{\mathbb{Z}}

\newcommand{\bA}{{\boldsymbol{A}}}

\newcommand{\bB}{{\boldsymbol{B}}}

\newcommand{\tb}{\widetilde{b}}

\newcommand{\tbB}{{\widetilde{\bB}}}

\newcommand{\bc}{{\boldsymbol{c}}}

\newcommand{\tC}{\widetilde{C}}

\newcommand{\be}{{\boldsymbol{e}}}

\newcommand{\Bf}{{\boldsymbol{f}}}

\newcommand{\bI}{{\boldsymbol{I}}}

\newcommand{\bJ}{{\boldsymbol{J}}}

\newcommand{\bP}{{\boldsymbol{P}}}

\newcommand{\br}{{\boldsymbol{r}}}

\newcommand{\bu}{{\boldsymbol{u}}}
\newcommand{\tbu}{\widetilde{\bu}}
\newcommand{\bv}{{\boldsymbol{v}}}

\newcommand{\bx}{{\boldsymbol{x}}}
\newcommand{\bX}{{\boldsymbol{X}}}
\newcommand{\by}{{\boldsymbol{y}}}
\newcommand{\bY}{{\boldsymbol{Y}}}

\newcommand{\bz}{{\boldsymbol{z}}}
\newcommand{\bZ}{{\boldsymbol{Z}}}

\newcommand{\Bbeta}{{\boldsymbol{\beta}}}

\newcommand{\tBbeta}{\widetilde{\Bbeta}}

\newcommand{\blambda}{{\boldsymbol{\lambda}}}

\newcommand{\bmu}{{\boldsymbol{\mu}}}

\newcommand{\bbeta}{{\boldsymbol{\beta}}}
\newcommand{\bDelta}{{\boldsymbol{\Delta}}}

\newcommand{\bzero}{{\boldsymbol{0}}}

\newcommand{\cA}{{\mathcal A}}

\newcommand{\cD}{{\mathcal D}}

\newcommand{\cF}{{\mathcal F}}

\newcommand{\cU}{{\mathcal U}}

\newcommand{\cc}{\mathrm{c}}
\newcommand{\dd}{\mathrm{d}}
\newcommand{\ee}{\mathrm{e}}

\newcommand{\EE}{\operatorname{\mathbb{E}}}

\newcommand{\PP}{\operatorname{\mathbb{P}}}

\renewcommand{\Re}{{\operatorname{Re}}}
\renewcommand{\Im}{{\operatorname{Im}}}

\newcommand{\tM}{\widetilde{M}}
\newcommand{\tN}{\widetilde{N}}

\newcommand{\vare}{\varepsilon}

\renewcommand{\mid}{\,|\,}

\renewcommand{\leq}{\leqslant}
\renewcommand{\geq}{\geqslant}

\newcommand{\distr}{\stackrel{\cD}{\longrightarrow}}
\newcommand{\distre}{\stackrel{\cD}{=}}
\newcommand{\mean}{\stackrel{L_1}{\longrightarrow}}
\newcommand{\qmean}{\stackrel{L_2}{\longrightarrow}}

\newcommand{\as}{\stackrel{{\mathrm{a.s.}}}{\longrightarrow}}
\newcommand{\ase}{\stackrel{{\mathrm{a.s.}}}{=}}

\newcommand{\bbone}{\mathbbm{1}}

\newcommand{\proofend}{\hfill\mbox{$\Box$}}

\numberwithin{equation}{section}

\theoremstyle{change} \theorembodyfont{\em}
\newtheorem{Lem}{Lemma.}[section]
\newtheorem{Thm}[Lem]{Theorem.}
\newtheorem{Pro}[Lem]{Proposition.}

\newtheorem{Def}[Lem]{Definition.}

\theorembodyfont{\rm}
\newtheorem{Rem}[Lem]{Remark.}

\begin{document}

\begin{center}
 {\bfseries\Large
   Almost sure, $L_1$- and $L_2$-growth behavior of supercritical\\[1mm]
    multi-type continuous state and continuous time\\[3mm]
    branching processes with immigration}

\vspace*{3mm}

 {\sc\large
  M\'aty\'as $\text{Barczy}^{*,\diamond}$,
  \ Sandra $\text{Palau}^{**}$,
  \ Gyula $\text{Pap}^{***}$}

\end{center}

\vskip0.2cm

\noindent
 * MTA-SZTE Analysis and Stochastics Research Group,
   Bolyai Institute, University of Szeged,
   Aradi v\'ertan\'uk tere 1, H--6720 Szeged, Hungary.

\noindent
 ** Department of Statistics and Probability, Instituto de Investigaciones en Matem\'aticas Aplicadas y en Sistemas,
 Universidad Nacional Aut\'onoma de M\'exico, M\'exico.

\noindent
 *** Bolyai Institute, University of Szeged,
    Aradi v\'ertan\'uk tere 1, H-6720 Szeged, Hungary.

\noindent e-mail: barczy@math.u-szeged.hu (M. Barczy),
                  sandra@sigma.iimas.unam.mx (S. Palau),
                  papgy@math.u-szeged.hu (G. Pap).

\noindent $\diamond$ Corresponding author.



\renewcommand{\thefootnote}{}
\vspace*{0.2cm}
\footnote{M\'aty\'as Barczy is supported by the J\'anos Bolyai Research Scholarship
 of the Hungarian Academy of Sciences.
 Sandra Palau was supported by the Royal Society Newton International Fellowship
 and by the EU-funded Hungarian grant EFOP-3.6.1-16-2016-00008.}

\vspace*{-10mm}

\begin{abstract}
Under a first order moment condition on the immigration mechanism, we show that an
 appropriately scaled supercritical and irreducible multi-type continuous state and
 continuous time branching process with immigration (CBI process) converges almost
 surely.
If an $x \log(x)$ moment condition on the branching mechanism does not hold, then
 the limit is zero.
If this $x \log(x)$ moment condition holds, then we prove $L_1$ convergence as
 well.
The projection of the limit on any left non-Perron eigenvector of the branching
 mean matrix is vanishing.
If, in addition, a suitable extra power moment condition on the branching
 mechanism holds, then we provide the correct scaling for the projection of a CBI
 process on certain left non-Perron eigenvectors of the branching mean matrix in
 order to have almost sure and $L_1$ limit.
Moreover, under a second order moment condition on the branching and immigration
 mechanisms, we prove $L_2$ convergence of an appropriately scaled process and the
 above mentioned projections as well.
A representation of the limits is also provided under the same moment conditions.
\end{abstract}

\section{Introduction}
\label{section_intro}

The description of the asymptotic behavior of branching processes without or with
 immigration has a long history.
For multi-type Galton--Watson processes without immigration see, e.g., Athreya and
 Ney \cite[Sections 4--8 in Chapter V]{AthNey}.
For supercritical multi-type Galton--Watson processes with immigration see, e.g.,
 Kaplan \cite{Kap}.

Let us consider a multi-type continuous state and continuous time branching process
 with immigration (CBI process) which can be represented as a pathwise unique
 strong solution of the stochastic differential equation (SDE)
 \begin{align}\label{SDE_atirasa_dimd}
  \begin{split}
   \bX_t
   &=\bX_0
     + \int_0^t (\Bbeta + \tbB \bX_u) \, \dd u
     + \sum_{\ell=1}^d
        \int_0^t \sqrt{2 c_\ell \max \{0, X_{u,\ell}\}} \, \dd W_{u,\ell}
        \, \be_\ell \\
   &\quad
      + \sum_{\ell=1}^d
         \int_0^t \int_{\cU_d} \int_{\cU_1}
          \bz \bbone_{\{w\leq X_{u-,\ell}\}} \, \tN_\ell(\dd u, \dd\bz, \dd w)
      + \int_0^t \int_{\cU_d} \br \, M(\dd u, \dd\br)
  \end{split}
 \end{align}
 for \ $t \in [0, \infty)$, \ see, Theorem 4.6 and Section 5 in
 Barczy et al.~\cite{BarLiPap2}, where \eqref{SDE_atirasa_dimd} was proved only for
 \ $d \in \{1, 2\}$, \ but their method clearly works for all
 \ $d \in \{1, 2, \ldots\}$.
\ Here \ $d \in \{1, 2, \ldots\}$ \ is the number of types, \ $X_{t,\ell}$,
 \ $\ell \in \{1, \ldots, d\}$, \ denotes the $\ell^{\mathrm th}$ coordinate of
 \ $\bX_t$, \ $\PP(\bX_0\in[0,\infty)^d)=1$,
  \ $\Bbeta \in [0, \infty)^d$,
 \ $c_1, \ldots, c_d \in [0, \infty)$,
 \ $\be_1, \ldots, \be_d$ \ denotes the natural basis in \ $\RR^d$,
 \ $\cU_d := [0, \infty)^d \setminus \{(0, \ldots, 0)\}$, \ $(W_{t,1})_{t\geq0}$,
 \ \ldots, \ $(W_{t,d})_{t\geq0}$ \ are independent standard Wiener processes,
 \ $N_\ell$, \ $\ell \in \{1, \ldots, d\}$, \ and \ $M$ \ are independent Poisson
 random measures on \ $(0, \infty) \times \cU_d \times (0, \infty)$ \ and on
 \ $(0, \infty) \times \cU_d$ \ with intensity measures
 \ $\dd u \, \mu_\ell(\dd\bz) \, \dd w$, \ $\ell \in \{1, \ldots, d\}$, \ and
 \ $\dd u \, \nu(\dd\br)$, \ respectively, and
 \ $\tN_\ell(\dd u, \dd\bz, \dd w)
    := N_\ell(\dd u, \dd\bz, \dd w) - \dd u \, \mu_\ell(\dd\bz) \, \dd w$,
 \ $\ell \in \{1, \ldots, d\}$.
\ We suppose that \ $\EE(\Vert\bX_0\Vert)<\infty$, \
 the Borel measures \ $\mu_\ell$, \ $\ell \in \{1, \ldots, d\}$,
 \ and \ $\nu$ \ on \ $\cU_d$ \ satisfy the moment conditions given in parts
 (v), (vi) of Definition \ref{Def_admissible} and \eqref{moment_condition_m_new},
  and \ $\bX_0$, \ $(W_{t,1})_{t\geq0}$, \ldots, \ $(W_{t,d})_{t\geq0}$, \ $N_1$,
 \ldots, \ $N_d$ \ and \ $M$ \ are independent.
Moreover, \ $\tbB = (\tb_{i,j})_{i,j\in\{1,\ldots,d\}} \in \RR^{d\times d}$ \ is a
 matrix satisfying
 \ $\tb_{i,j} \geq \int_{\cU_d} z_i \, \mu_j(\dd\bz)$ \ for all
 \ $i, j \in \{1, \ldots, d\}$ \ with \ $i \ne j$.

A multi-type CBI process \ $(\bX_t)_{t\in\RR_+}$ \ is called irreducible if
 \ $\tbB$ \ is irreducible, see Definition \ref{Def_irreducible}.
An irreducible multi-type CBI process is called subcritical, critical or
 supercritical if the logarithm \ $s(\tbB)$ \ of the Perron eigenvalue of the
 branching mean matrix \ $\ee^\tbB$ \ is negative, zero or positive, respectively,
 see Definition \ref{Def_indecomposable_crit}.
A multi-type CBI process \ $(\bX_t)_{t\in\RR_+}$ \ is called a multi-type CB
 process if there is no immigration, i.e., \ $\Bbeta = \bzero$ \ and \ $\nu = 0$.

In case of a subcritical or critical single-type CBI process (when it is necessarily irreducible)
 with a non-vanishing branching mechanism, \ $\bX_t \distr \pi$ \ as \ $t \to \infty$ \ with a probability measure \ $\pi$
 \ on \ $[0, \infty)$ \ if and only if certain integrability condition holds for
 the branching and immigration mechanisms, see, e.g., Li \cite[Theorem 3.20]{Li}.

In case of a supercritical single-type CB process, under the $x \log(x)$ moment
 condition \eqref{moment_condition_xlogx} with \ $\lambda = s(\tbB)$ \ on the
 branching mechanism, Li \cite[Corollary 3.16 and Theorem 3.8]{Li} proved that
 \ $\ee^{-s(\tbB)t} \bX_t$ \ converges almost surely as \ $t \to \infty$ \ towards
 a non-negative random variable, and the probability that this limit is zero equals
 to the probability of the event that the extinction time is finite.

In case of a critical and irreducible multi-type CBI process, under fourth order
 moment conditions on the branching and immigration mechanisms, Barczy and Pap
 \cite[Theorem 4.1]{BarPap} proved that the sequence
 \ $(n^{-1} \bX_{\lfloor nt \rfloor})_{t\in[0,\infty)}$,
 \ $n \in \{1, 2, \ldots\}$, \ of scaled random step functions converges weakly
 towards a squared Bessel process (in other words, a Feller diffusion) supported by
 a ray determined by the right Perron vector \ $\tbu$ \ of the branching mean
 matrix \ $\ee^{\tbB}$.

Recently, there is a renewed interest for studying asymptotic behavior of
 supercritical branching processes.
In case of a supercritical and irreducible multi-type CB process,
 Kyprianou et al.\ \cite[Theorem 1.3]{KypPalRen} described the asymptotic behavior
 of the projection \ $\langle\bu, \bX_t\rangle$ \ as \ $t \to \infty$, \ where
 \ $\bu$ \ denotes the left Perron eigenvector of the branching mean matrix
 \ $\ee^{\tbB}$.
\ Namely, they proved that if an $x \log(x)$ moment condition on the branching
 mechanism holds, then
 \ $\ee^{-s(\tbB)t} \langle\bu, \bX_t\rangle \to w_{\bu,\bX_0}$ \ almost surely and
 in $L_1$ as \ $t \to \infty$, \ where \ $w_{\bu,\bX_0}$ \ is a non-negative random
 variable, otherwise \ $\ee^{-s(\tbB)t} \langle\bu, \bX_t\rangle \to 0$ \ almost
 surely as \ $t \to \infty$.
\ Note that their $x \log(x)$ moment condition is equivalent to our moment
 condition \eqref{moment_condition_xlogx} with \ $\lambda = s(\tbB)$, \ since for
 \ $\RR^d$, \ all norms are equivalent.
Moreover, in case of a supercritical and irreducible multi-type CB process,
 Kyprianou et al.\ \cite[Theorem 1.4]{KypPalRen} proved that
 \ $\ee^{-s(\tbB)t} \bX_t \to w_{\bu,\bX_0} \tbu$ \ almost surely as
 \ $t \to \infty$.

Ren et al.\ \cite{RenSonZha1} investigated central limit theorems for supercritical
 branching Markov processes, and Ren et al.\ \cite{RenSonZha3, RenSonZha2} studied some
 properties of strong limits for supercritical superprocesses.
Moreover, Chen et al.\ \cite{CheRenSon} and Ren et al.\ \cite{RenSonYan}
 studied spine decomposition and an $x\log x$ criterion for supercritical
 superprocesses with non-local branching mechanisms.

Recently, Marks and Mi{\l}o\'{s} \cite[Theorem 3.2]{MarMil} considered a branching
 particle system with particles moving according to a multi-dimensional Ornstein-Uhlenbeck process
 with a positive drift and branching according to a law in the domain of attraction of a stable law
 having stability index in \ $(1,2)$, \ and in the so-called ''large branching case'' (see \cite[page 3]{MarMil})
 they proved almost sure and $L_1$ convergence of appropriately normalized projections
 of the particle system in question onto certain twice differentiable real-valued functions defined on the real line
 of polynomial growth together with a description of the limit in which the whole genealogical structure is somewhat preserved.
These projections include projections onto certain eigenfunctions of the semigroup associated to the infinitesimal
 generator of the underlying Ornstein-Uhlenbeck process.

Very recently, Ren et al.\ \cite{RenSonSunZhao} derived stable central limit theorems for some kind of projections of
 (measure-valued) super Ornstein-Uhlenbeck processes having a branching mechanism which is close to a function
 of the form \ $-a_1 z + a_2 z^2 + a_3 z^{1+\alpha}$, \ $z\geq 0$, \ with \ $a_1 > 0$, \ $a_2 \geq  0$, \ $a_3 > 0$ \ and \ $\alpha\in(0,1)$
 \ in some sense (see the Assumption 2 in Ren et al.\ \cite{RenSonSunZhao}).

As a new result, in case of a supercritical and irreducible multi-type CBI process,
 under the first order moment condition \eqref{moment_condition_m_new} on the
 immigration mechanism, we show \ $\ee^{-s(\tbB)t} \bX_t \to w_{\bu,\bX_0} \tbu$
 \ almost surely as \ $t \to \infty$, \ where \ $w_{\bu,\bX_0}$ \ is a non-negative
 random variable, see Theorem \ref{convCBIwasL1}.
If the $x \log(x)$ moment condition \eqref{moment_condition_xlogx} with
 \ $\lambda = s(\tbB)$ \ does not hold, then \ $\PP(w_{\bu,\bX_0} = 0) = 1$, \ see
 Theorem \ref{convCBIasL1}.
If this $x \log(x)$ moment condition holds, then we prove $L_1$ convergence, see
 Theorem \ref{convCBIwasL1}, and we give a representation of \ $w_{\bu,\bX_0}$ \ as
 well, see \eqref{reprCBIw}.
Note that \ $\PP(w_{\bu,\bX_0} = 0) = 1$ \ if and only if
 \ $\PP(\bX_t = \bzero) = 1$ \ for all \ $t \in \RR_+$, \ see Theorem \ref{convCBIasL1}.
Hence here the scaling factor \ $\ee^{-s(\tbB)t}$ \ is correct.
If \ $\bv$ \ is a left non-Perron eigenvector of the branching mean matrix
 \ $\ee^{\tbB}$, \ then this result implies that
 \ $\ee^{-s(\tbB)t} \langle\bv, \bX_t\rangle
    \to w_{\bu,\bX_0} \langle\bv, \tbu\rangle = 0$
 \ almost surely as \ $t \to \infty$, \ since \ $\langle\bv, \tbu\rangle = 0$ \ due
 to the so-called principle of biorthogonality (see, e.g., Horn and Johnson
 \cite[Theorem 1.4.7(a)]{HorJoh}), consequently, the scaling factor
 \ $\ee^{-s(\tbB)t}$ \ is not appropriate for describing the asymptotic behavior of
 the projection \ $\langle\bv, \bX_t\rangle$ \ as \ $t \to \infty$.
\ It turns out that, under the extra power moment condition
 \eqref{moment_condition_xlogx} with
 \ $\Re(\lambda) \in \bigl(\frac{1}{2} s(\tbB), s(\tbB)\bigr)$ \ on the branching
 mechanism and the first order moment condition \eqref{moment_condition_m_new} on
 the immigration mechanism, we can show
 \ $\ee^{-\lambda t} \langle\bv, \bX_t\rangle \to w_{\bv,\bX_0}$ \ almost surely
 and in \ $L_1$ \ as \ $t \to \infty$, \ where \ $\lambda$ \ is a non-Perron
 eigenvalue of the branching mean matrix \ $\ee^{\tbB}$ \ with
 \ $\Re(\lambda) \in \bigl(\frac{1}{2} s(\tbB), s(\tbB)\bigr)$, \ $\bv$ \ is a left
 eigenvector corresponding to \ $\lambda$, \ and \ $w_{\bv,\bX_0}$ \ is a complex
 random variable, see Theorem \ref{convCBIasL1}, where we give a representation of
 \ $w_{\bv,\bX_0}$ \ as well, see \eqref{reprCBIw}.
Here the scaling factor \ $\ee^{-\lambda t}$ \ is correct if
 \ $\langle\bv, \EE(\bX_0) + \lambda^{-1} \tBbeta\rangle \ne 0$, \ since then
 \ $\PP(w_{\bv,\bX_0} = 0) < 1$, \ see Theorem \ref{convCBIasL1}.
In Remark \ref{Rem_lambda_not_sB} we explain why we do not have any result in the case when the moment condition \eqref{moment_condition_xlogx}
 does not hold for \ $\lambda \in \big(\frac{1}{2}s(\tbB), s(\tbB)\big)$ \ formulating some open problems as well.
Note that the asymptotic behavior of the second moment
 \ $\EE(|\langle\bv, \bX_t\rangle|^2)$ \ as \ $t \to \infty$ \ explains the role of
 the assumption \ $\Re(\lambda) \in \bigl(\frac{1}{2} s(\tbB), s(\tbB)\bigr]$,
 \ see Proposition \ref{second_moment_asymptotics_CBI}.

Further, in case of a supercritical and irreducible multi-type CBI process, under
 the second order moment condition \eqref{moment_condition_CBI2} on the branching
 and immigration mechanisms, we show
 \ $\ee^{-s(\tbB)t} \bX_t \to w_{\bu,\bX_0} \tbu$ \ and
 \ $\ee^{-\lambda t} \langle\bv, \bX_t\rangle \to w_{\bv,\bX_0}$ \ in \ $L_2$ \ as
 \ $t \to \infty$ \ as well, where \ $\lambda$ \ is a eigenvalue of the branching
 mean matrix \ $\ee^{\tbB}$ \ with
 \ $\Re(\lambda) \in \bigl(\frac{1}{2} s(\tbB), s(\tbB)\bigr]$ \ and \ $\bv$ \ is a
 left eigenvector corresponding to \ $\lambda$, \ see Theorem \ref{convCBIL2}.

The paper is structured as follows.
In Section \ref{section_CBI}, for completeness and better readability, from
 Barczy et al.~\cite{BarLiPap2}, we recall some notions and statements for
 multi-type CBI processes such as a formula for their first moment, an appropriate
 transformation which results in a $d$-dimensional martingale in Lemma \ref{mart},
 a useful representation of \ $(\bX_t)_{t\in\RR_+}$ \ in Lemma
 \ref{SDE_transform_sol}, the definition of subcritical, critical and supercritical
 irreducible CBI processes, see Definitions \ref{Def_irreducible} and
 \ref{Def_indecomposable_crit}.
Section \ref{section_CBI_as_growth} contains our main results detailed above, see
 Theorems \ref{convCBIasL1}, \ref{convCBIwasL1} and \ref{convCBIL2}.
For the proofs, we use heavily the representation of \ $(\bX_t)_{t\in\RR_+}$ \ in
 Lemma \ref{SDE_transform_sol} based on the SDE \eqref{SDE_atirasa_dimd}.
In the course of the proof of Theorem \ref{convCBIwasL1}, we follow the steps and
 methods of the proof of Theorem 1.4 in Kyprianou et al.\ \cite{KypPalRen}.
We close the paper with two Appendices.
We present a useful decomposition of a CBI process as an independent sum of a CBI
 process starting from $\bzero$ and a CB process, see Appendix \ref{deco_CBI}.
In Appendix \ref{second_moment_CBI}, we describe the asymptotic behavior of the
 second moment of \ $|\langle\bv, \bX_t\rangle|$ \ as \ $t \to \infty$ \ for each
 left eigenvector \ $\bv \in \CC^d$ \ of \ $\tbB$ \ corresponding to an arbitrary
 eigenvalue \ $\lambda \in \sigma(\tbB)$ \ in case of a supercritical and
 irreducible CBI process.

Now, we summarize the novelties of the paper.
We point out that we investigate the asymptotic behavior of the projections of a
 multi-type CBI process on certain left non-Perron eigenvectors of its branching
 mean matrix.
According to our knowledge, this type of question has not been studied so far for
 multi-type CBI processes.
A new phenomenon appears compared to the left Perron eigenvector case, namely, a
 moment type condition on the branching mechanism of the CBI process in question.
Furthermore, if the $x \log(x)$ moment condition \eqref{moment_condition_xlogx} with
 \ $\lambda = s(\tbB)$ \ on the branching mechanism does not hold, then one
 usually uses a so-called spine decomposition technique in order to show that \ $w_{\bu,\bX_0} \ase 0$ \
 (see, e.g., the proof of Theorem 1.3 in Kyprianou et al.\ \cite{KypPalRen} or that of Theorem 6.2
  in Ren et al.\ \cite{RenSonYan}).
In this paper, we use that the law of a multi-type CBI process \ $(\bX_t)_{t\in\RR_+}$ \
 at time \ $t+T$, \ $t,T\in\RR_+$, \ coincides with the law of an independent sum of a multi-type CB process at time \ $t$ \ starting from
 an initial value having distribution as that of \ $\bX_T$ \ and a multi-type CBI process  at time \ $t$ \ starting from
 \ $\bzero$, \ presented in Lemma \ref{decomposition_CBI}, and that the corresponding result
 \ $w_{\bu,\bX_0} \ase 0$ \ is already known for CB processes due to Kyprianou et al.\ \cite[Theorem 1.3]{KypPalRen}.

Finally, we mention a possible extension of the present results which can be a topic of future work.
Since the \ $d$-dimensional matrix \ $\tbB$ \ is not symmetric in general, its left eigenvectors may not generate \ $\CC^d$,
 \ so it is natural to study the asymptotic behaviour of \ $\langle \bv, \bX_t\rangle$ \ as \ $t\to\infty$, \
 where \ $\bv$ \ is an arbitrary vector in \ $\CC^d$.
This type of question was investigated by Kesten and Stigun \cite{KesSti} and Badalbaev and Mukhitdinov \cite{BadMuk}
 for supercritical and irreducible multi-type discrete time Galton--Watson processes without immigration under second order moment assumptions,
 and, by Athreya \cite{Ath1, Ath2}, for supercritical and positively regular
 multi-type continuous time Markov branching processes without immigration under second order moment assumptions.
The above mentioned four references are for some branching processes without immigration, we do not know any corresponding result
 for branching processes with immigration.
Motivated by these references, we think that the Jordan normal form of \ $\tbB$ \ may be well-used in our case as well,
 where we consider multi-type CBI processes with immigration.

\section{Multi-type CBI processes}
\label{section_CBI}

Let \ $\ZZ_+$, \ $\NN$, \ $\RR$, \ $\RR_+$, \ $\RR_{++}$ \ and \ $\CC$ \ denote the
 set of non-negative integers, positive integers, real numbers, non-negative real
 numbers, positive real numbers and complex numbers, respectively.
For \ $x , y \in \RR$, \ we will use the notations
 \ $x \land y := \min \{x, y\} $ \ and \ $x^+:= \max \{0, x\}$.
\ By \ $\langle\bx, \by\rangle := \sum_{j=1}^d x_j \overline{y_j}$, \ we denote the
 Euclidean inner product of \ $\bx = (x_1, \ldots, x_d)^\top \in \CC^d$ \ and
 \ $\by = (y_1, \ldots, y_d)^\top \in \CC^d$, \ and by \ $\|\bx\|$ \ and
 \ $\|\bA\|$, \ we denote the induced norm of \ $\bx \in \CC^d$ \ and
 \ $\bA \in \CC^{d\times d}$, \ respectively.
The null vector and the null matrix will be denoted by \ $\bzero$.
\ Moreover, \ $\bI_d \in \RR^{d\times d}$ \ denotes the identity matrix.
By \ $C^2_\cc(\RR_+^d,\RR)$, \ we denote the set of twice continuously
 differentiable real-valued functions on \ $\RR_+^d$ \ with compact support.
Convergence almost surely, in $L_1$ and in $L_2$ will be denoted by
 \ $\as$, \ $\mean$ \ and \ $\qmean$, \ respectively.
Almost sure equality will be denoted by \ $\ase$.
\ Throughout this paper, we make the conventions \ $\int_a^b := \int_{(a,b]}$
 \ and \ $\int_a^\infty := \int_{(a,\infty)}$ \ for any \ $a, b \in \RR$ \ with
 \ $a < b$.

\begin{Def}\label{Def_essentially_non-negative}
A matrix \ $\bA = (a_{i,j})_{i,j\in\{1,\ldots,d\}} \in \RR^{d\times d}$ \ is
 called essentially non-negative if \ $a_{i,j} \in \RR_+$ \ whenever
 \ $i, j \in \{1,\ldots,d\}$ \ with \ $i \ne j$, \ that is, if \ $\bA$ \ has
 non-negative off-diagonal entries.
The set of essentially non-negative \ $d \times d$ \ matrices will be denoted by
 \ $\RR^{d\times d}_{(+)}$.
\end{Def}

\begin{Def}\label{Def_admissible}
A tuple \ $(d, \bc, \Bbeta, \bB, \nu, \bmu)$ \ is called a set of admissible
 parameters if
 \renewcommand{\labelenumi}{{\rm(\roman{enumi})}}
 \begin{enumerate}
  \item
   $d \in \NN$,
  \item
   $\bc = (c_i)_{i\in\{1,\ldots,d\}} \in \RR_+^d$,
  \item
   $\Bbeta = (\beta_i)_{i\in\{1,\ldots,d\}} \in \RR_+^d$,
  \item
   $\bB = (b_{i,j})_{i,j\in\{1,\ldots,d\}} \in \RR^{d \times d}_{(+)}$,
  \item
   $\nu$ \ is a Borel measure on \ ${\cU_d := \RR_+^d \setminus \{\bzero\}}$
    \ satisfying \ $\int_{\cU_d} (1 \land \|\br\|) \, \nu(\dd\br) < \infty$,
  \item
   $\bmu = (\mu_1, \ldots, \mu_d)$, \ where, for each
    \ $i \in \{1, \ldots, d\}$, \ $\mu_i$ \ is a Borel measure on
    \ $\cU_d$ \ satisfying
    \begin{align}\label{help_moment_mu}
      \int_{\cU_d}
       \biggl[\|\bz\| \land \|\bz\|^2
              + \sum_{j \in \{1, \ldots, d\} \setminus \{i\}} (1 \land z_j)\biggr]
       \mu_i(\dd\bz)
      < \infty .
    \end{align}
  \end{enumerate}
\end{Def}

\begin{Rem}
Our Definition \ref{Def_admissible} of the set of admissible parameters is a
 special case of Definition 2.6 in Duffie et al.~\cite{DufFilSch}, which is
 suitable for all affine processes, see Barczy et al.~\cite[Remark 2.3]{BarLiPap2}.
Further, due to Remark 2.3 and (2.12) in Barczy et al. \cite{BarLiPap2}, the
 condition \eqref{help_moment_mu} is equivalent to
 \[
   \int_{\cU_d}
    \biggl[\|\bz\| \land \|\bz\|^2
           + \sum_{j\in\{1,\ldots,d\} \setminus \{i\}} z_j\biggr]
    \mu_i(\dd\bz)
   < \infty ,
 \]
 and also to
 \[
   \int_{\cU_d}
    \biggl[(1 \land z_i)^2
           + \sum_{j\in\{1,\ldots,d\} \setminus \{i\}} (1 \land z_j)\biggr]
    \mu_i(\dd\bz)
   < \infty \qquad \text{and} \qquad
   \int_{\cU_d}
    \|\bz\| \bbone_{\{\|\bz\|\geq1\}} \, \mu_i(\dd\bz) < \infty .
    \vspace*{-4mm}
 \]
  \proofend
\end{Rem}

\begin{Thm}\label{CBI_exists}
Let \ $(d, \bc, \Bbeta, \bB, \nu, \bmu)$ \ be a set of admissible parameters.
Then, there exists a unique conservative transition semigroup \ $(P_t)_{t\in\RR_+}$
 \ acting on the Banach space (endowed with the supremum norm) of real-valued
 bounded Borel-measurable functions on the state space \ $\RR_+^d$ \ such that its
 infinitesimal generator is
  \begin{align*}
   (\cA f)(\bx)
   &= \sum_{i=1}^d c_i x_i f_{i,i}''(\bx)
      + \langle \Bbeta + \bB \bx, \Bf'(\bx) \rangle
      + \int_{\cU_d} \bigl( f(\bx + \br) - f(\bx) \bigr) \, \nu(\dd\br) \\
   &\phantom{\quad}
      + \sum_{i=1}^d
         x_i
         \int_{\cU_d}
          \bigl( f(\bx + \bz) - f(\bx) - f'_i(\bx) (1 \land z_i) \bigr)
          \, \mu_i(\dd \bz)
  \end{align*}
 for \ $f \in C^2_\cc(\RR_+^d,\RR)$ \ and \ $\bx \in \RR_+^d$, \ where \ $f_i'$
 \ and \ $f_{i,i}''$, \ $i \in \{1, \ldots, d\}$, \ denote the first and second
 order partial derivatives of \ $f$ \ with respect to its \ $i$-th variable,
 respectively, and \ $\Bf'(\bx) := (f_1'(\bx), \ldots, f_d'(\bx))^\top$.
\ Moreover, the Laplace transform of the transition semigroup \ $(P_t)_{t\in\RR_+}$
 \ has a representation
 \begin{equation}\label{Laplace_transform}
  \int_{\RR_+^d} \ee^{- \langle \blambda, \by \rangle} P_t(\bx, \dd \by)
  = \ee^{- \langle \bx, \bv(t, \blambda) \rangle
         - \int_0^t \psi(\bv(s, \blambda)) \, \dd s} , \qquad
  \bx \in \RR_+^d, \quad \blambda \in \RR_+^d , \quad t \in \RR_+ ,
 \end{equation}
 where, for any \ $\blambda \in \RR_+^d$, \ the continuously differentiable
 function
 \ $\RR_+ \ni t \mapsto \bv(t, \blambda)
    = (v_1(t, \blambda), \ldots, v_d(t, \blambda))^\top \in \RR_+^d$
 \ is the unique locally bounded solution to the system of differential equations
 \[
   \partial_t v_i(t, \blambda) = - \varphi_i(\bv(t, \blambda)) , \qquad
   v_i(0, \blambda) = \lambda_i , \qquad i \in \{1, \ldots, d\} ,
 \]
 with
 \[
   \varphi_i(\blambda)
   := c_i \lambda_i^2 -  \langle \bB \be_i, \blambda \rangle
      + \int_{\cU_d}
         \bigl( \ee^{- \langle \blambda, \bz \rangle} - 1
                + \lambda_i (1 \land z_i) \bigr)
         \, \mu_i(\dd \bz)
 \]
 for \ $\blambda \in \RR_+^d$, \ $i \in \{1, \ldots, d\}$, \ and
 \[
   \psi(\blambda)
   := \langle \bbeta, \blambda \rangle
      + \int_{\cU_d}
         \bigl( 1 - \ee^{- \langle\blambda, \br\rangle} \bigr)
         \, \nu(\dd\br) , \qquad
   \blambda \in \RR_+^d .
 \]
\end{Thm}

Theorem \ref{CBI_exists} is a special case of Theorem 2.7 of Duffie et al.\
 \cite{DufFilSch} with \ $m = d$, \ $n = 0$ \ and zero killing rate.
For more details, see Remark 2.5 in Barczy et al.\ \cite{BarLiPap2}.

\begin{Def}\label{Def_CBI}
A conservative Markov process with state space \ $\RR_+^d$ \ and with transition
 semigroup \ $(P_t)_{t\in\RR_+}$ \ given in Theorem \ref{CBI_exists} is called a
 multi-type CBI process with parameters \ $(d, \bc, \Bbeta, \bB, \nu, \bmu)$.
\ The function
 \ $\RR_+^d \ni \blambda
    \mapsto (\varphi_1(\blambda), \ldots, \varphi_d(\blambda))^\top \in \RR^d$
 \ is called its branching mechanism, and the function
 \ $\RR_+^d \ni \blambda \mapsto \psi(\blambda) \in \RR_+$ \ is called its
 immigration mechanism.
A multi-type CBI process with parameters \ $(d, \bc, \Bbeta, \bB, \nu, \bmu)$ \ is
 called a CB process (a continuous state and continuous time branching process
 without immigration) if \ $\Bbeta = \bzero$ \ and \ $\nu = 0$.
\end{Def}

Let \ $(\bX_t)_{t\in\RR_+}$ \ be a multi-type CBI process with parameters
 \ $(d, \bc, \Bbeta, \bB, \nu, \bmu)$ \ such that \ $\EE(\|\bX_0\|) < \infty$ \ and
 the moment condition
 \begin{equation}\label{moment_condition_m_new}
  \int_{\cU_d} \|\br\| \bbone_{\{\|\br\|\geq1\}} \, \nu(\dd\br) < \infty
 \end{equation}
 holds.
Then, by formula (3.4) in Barczy et al. \cite{BarLiPap2},
 \begin{equation}\label{EXcond}
  \EE(\bX_t \mid \bX_0 = \bx)
  = \ee^{t\tbB} \bx + \int_0^t \ee^{u\tbB} \tBbeta \, \dd u ,
  \qquad \bx \in \RR_+^d , \quad t \in \RR_+ ,
 \end{equation}
 where
 \begin{gather*}
  \tbB := (\tb_{i,j})_{i,j\in\{1,\ldots,d\}} , \qquad
  \tb_{i,j}
  := b_{i,j} + \int_{\cU_d} (z_i - \delta_{i,j})^+ \, \mu_j(\dd\bz) , \qquad
  \tBbeta := \Bbeta + \int_{\cU_d} \br \, \nu(\dd\br) ,
 \end{gather*}
 with \ $\delta_{i,j}:=1$ \ if \ $i = j$, \ and \ $\delta_{i,j} := 0$ \ if
 \ $i \ne j$.
\ Note that \ $\tbB \in \RR^{d \times d}_{(+)}$ \ and \ $\tBbeta \in \RR_+^d$,
 \ since
 \[
   \int_{\cU_d} \|\br\| \, \nu(\dd\br) < \infty , \qquad
   \int_{\cU_d} (z_i - \delta_{i,j})^+ \, \mu_j(\dd \bz) < \infty , \quad
   i, j \in \{1, \ldots, d\} ,
 \]
 see Barczy et al. \cite[Section 2]{BarLiPap2}.
Further, \ $\EE(\bX_t \mid \bX_0 = \bx)$, \ $\bx \in \RR_+^d$, \ does not depend on
 the parameter \ $\bc$.
\ One can give probabilistic interpretations of the modified parameters \ $\tbB$
 \ and \ $\tBbeta$, \ namely, \ $\ee^\tbB \be_j = \EE(\bY_1 \mid \bY_0 = \be_j)$,
 \ $j \in \{1, \ldots, d\}$, \ and \ $\tBbeta = \EE(\bZ_1 \mid \bZ_0 = \bzero)$,
 \ where \ $(\bY_t)_{t\in\RR_+}$ \ and \ $(\bZ_t)_{t\in\RR_+}$ \ are multi-type CBI
 processes with parameters \ $(d, \bc, \bzero, \bB, 0, \bmu)$ \ and
 \ $(d, \bzero, \Bbeta, \bzero, \nu, \bzero)$, \ respectively, see formula
 \eqref{EXcond}.
The processes \ $(\bY_t)_{t\in\RR_+}$ \ and \ $(\bZ_t)_{t\in\RR_+}$ \ can be
 considered as pure branching (without immigration) and pure immigration (without
 branching) processes, respectively.
Consequently, \ $\ee^\tbB$ \ and \ $\tBbeta$ \ may be called the branching mean
 matrix and the immigration mean vector, respectively.
Note that the branching mechanism depends only on the parameters \ $\bc$, \ $\bB$
 \ and \ $\bmu$, \ while the immigration mechanism depends only on the parameters
 \ $\Bbeta$ \ and \ $\nu$.

\begin{Lem}\label{mart}
Let \ $(\bX_t)_{t\in\RR_+}$ \ be a multi-type CBI process with parameters
 \ $(d, \bc, \Bbeta, \bB, \nu, \bmu)$ \ such that \ $\EE(\|\bX_0\|) < \infty$ \ and
 the moment condition \eqref{moment_condition_m_new} holds.
Then the process
 \ $\bigl(\ee^{-t\tbB} \bX_t
          - \int_0^t \ee^{-u\tbB} \tBbeta \, \dd u\bigr)_{t\in\RR_+}$
 \ is a \ $d$-dimensional martingale with respect to the filtration
 \ $\cF_t^\bX := \sigma(\bX_u : u \in [0, t])$, \ $t \in \RR_+$.
\end{Lem}

\noindent
\textbf{Proof.}
First, note that for all \ $t \in \RR_+$, \ $\bX_t$ \ is measurable with respect to
 \ $\cF_t^{\bX}$, \ and due to \ $\EE(\|\bX_0\|) < \infty$ \ and
 \eqref{moment_condition_m_new}, by Lemma 3.4 in Barczy et al.\ \cite{BarLiPap2},
 we have \ $\EE(\|\bX_t\|) < \infty$.
\ For each \ $v, t \in \RR_+$ \ with \ $v \leq t$, \ we have
 \[
   \EE(\bX_t \mid \cF_v^\bX)
   = \EE(\bX_t \mid \bX_v)
   = \ee^{(t-v)\tbB} \bX_v + \int_0^{t-v} \ee^{w\tbB} \tBbeta \, \dd w ,
 \]
 since \ $(\bX_t)_{t\in\RR_+}$ \ is a time-homogeneous Markov process, and we can
 apply \eqref{EXcond}.
Thus for each \ $v,t\in\RR_+$ \ with \ $v\leq t$, \ we obtain
 \begin{align*}
  &\EE\biggl(\ee^{-t\tbB} \bX_t - \int_0^t \ee^{-u\tbB} \tBbeta \, \dd u
            \,\bigg|\, \cF_v^\bX\biggr)
   = \ee^{-t\tbB} \ee^{(t-v)\tbB} \bX_v
     + \ee^{-t\tbB} \int_0^{t-v} \ee^{w\tbB} \tBbeta \, \dd w
     - \int_0^t \ee^{-u\tbB} \tBbeta \, \dd u \\
  &= \ee^{-v\tbB} \bX_v + \int_0^{t-v} \ee^{(w-t)\tbB} \tBbeta \, \dd w
     - \int_0^t \ee^{-u\tbB} \tBbeta \, \dd u
   = \ee^{-v\tbB} \bX_v - \int_0^v \ee^{-u\tbB} \tBbeta \, \dd u ,
 \end{align*}
 and consequently, the process
 \ $\bigl(\ee^{-t\tbB} \bX_t
          - \int_0^t \ee^{-u\tbB} \tBbeta \, \dd u\bigr)_{t\in\RR_+}$ \ is a
 martingale with respect to the filtration \ $(\cF_t^\bX)_{t\in\RR_+}$.
\proofend

By an application of the multidimensional It\^o's formula one can derive the
 following useful representation of \ $(\bX_t)_{t\in\RR_+}$, \ where the drift part
 is deterministic.
The proof can be found in Barczy et al.\ \cite[Lemma 4.1]{BarLiPap3}.

\begin{Lem}\label{SDE_transform_sol}
Let \ $(\bX_t)_{t\in\RR_+}$ \ be a multi-type CBI process with parameters
 \ $(d, \bc, \Bbeta, \bB, \nu, \bmu)$ \ such that \ $\EE(\|\bX_0\|) < \infty$ \ and
 the moment condition \eqref{moment_condition_m_new} holds.
Then, for each \ $s, t \in \RR_+$ \ with \ $s \leq t$, \ we have
 \begin{align*}
  \bX_t
  &= \ee^{(t-s)\tbB} \bX_s + \int_s^t \ee^{(t-u)\tbB} \tBbeta \, \dd u
     + \sum_{\ell=1}^d
        \int_s^t \ee^{(t-u)\tbB}
         \be_\ell \sqrt{2 c_\ell X_{u,\ell}} \, \dd W_{u,\ell} \\
  &\quad
     + \sum_{\ell=1}^d
        \int_s^t \int_{\cU_d} \int_{\cU_1}
         \ee^{(t-u)\tbB} \bz \bbone_{\{w\leq X_{u-,\ell}\}}
         \, \tN_\ell(\dd u, \dd\bz, \dd w)
      + \int_s^t \int_{\cU_d} \ee^{(t-u)\tbB} \br \, \tM(\dd u, \dd\br) ,
 \end{align*}
 where \ $\tM(\dd u, \dd\br) := M(\dd u, \dd\br) - \dd u \, \nu(\dd\br)$.
\end{Lem}

Note that the formula for \ $(\bX_t)_{t\in\RR_+}$ \ in Lemma
 \ref{SDE_transform_sol} is a generalization of the formula (3.1) in Xu \cite{Xu},
 and the formula (1.5) in Li and Ma \cite{LiMa}.

Next we recall a classification of multi-type CBI processes.
For a matrix \ $\bA \in \RR^{d \times d}$, \ $\sigma(\bA)$ \ will denote the
 spectrum of \ $\bA$, \ that is, the set of all \ $\lambda\in\CC$ \ that are eigenvalues of \ $\bA$.
\ Then \ $r(\bA) := \max_{\lambda \in \sigma(\bA)} |\lambda|$ \ is the spectral
 radius of \ $\bA$.
\ Moreover, we will use the notation
 \[
   s(\bA) := \max_{\lambda \in \sigma(\bA)} \Re(\lambda) .
 \]
A matrix \ $\bA \in \RR^{d\times d}$ \ is called reducible if there exist a
 permutation matrix \ $\bP \in \RR^{ d \times d}$ \ and an integer \ $r$ \ with
 \ $1 \leq r \leq d-1$ \ such that
 \[
  \bP^\top \bA \bP
   = \begin{bmatrix} \bA_1 & \bA_2 \\ \bzero & \bA_3 \end{bmatrix},
 \]
 where \ $\bA_1 \in \RR^{r\times r}$, \ $\bA_3 \in \RR^{ (d-r) \times (d-r) }$,
 \ $\bA_2 \in \RR^{r \times (d-r) }$, \ and \ $\bzero \in \RR^{(d-r)\times r}$ \ is
 a null matrix.
A matrix \ $\bA \in \RR^{d\times d}$ \ is called irreducible if it is not
 reducible, see, e.g., Horn and Johnson
 \cite[Definitions 6.2.21 and 6.2.22]{HorJoh}.
We do emphasize that no 1-by-1 matrix is reducible.

\begin{Def}\label{Def_irreducible}
Let \ $(\bX_t)_{t\in\RR_+}$ \ be a multi-type CBI process with parameters
 \ $(d, \bc, \Bbeta, \bB, \nu, \bmu)$ \ such that the moment condition
 \eqref{moment_condition_m_new} holds.
Then \ $(\bX_t)_{t\in\RR_+}$ \ is called irreducible if \ $\tbB$ \ is irreducible.
\end{Def}

Recall that if \ $\tbB \in \RR^{d\times d}_{(+)}$ \ is irreducible, then
 \ $\ee^{t\tbB} \in \RR^{d \times d}_{++}$ \ for all \ $t \in \RR_{++}$, \ and
 \ $s(\tbB)$ \ is an eigenvalue of \ $\tbB$, \ the algebraic and geometric
 multiplicities of \ $s(\tbB)$ \ is 1, and the real parts of the other eigenvalues
 of \ $\tbB$ \ are less than \ $s(\tbB)$.
\ Moreover, corresponding to the eigenvalue \ $s(\tbB)$ \ there exists a unique
 (right) eigenvector \ $\tbu \in \RR^d_{++}$ \ of \ $\tbB$ \  such that the sum of
 its coordinates is 1 which is also the unique (right) eigenvector of \ $\ee^\tbB$,
 \ called the right Perron vector of \ $\ee^\tbB$, \ corresponding to the
 eigenvalue \ $r(\ee^\tbB) = \ee^{s(\tbB)}$ \ of \ $\ee^\tbB$ \ such that the sum
 of its coordinates is 1.
\ Further, there exists a unique left eigenvector \ $\bu \in \RR^d_{++}$ \ of
 \ $\tbB$ \ corresponding to the eigenvalue \ $s(\tbB)$ \ with
 \ $\tbu^\top \bu = 1$, \ which is also the unique (left) eigenvector of
 \ $\ee^\tbB$, \ called the left Perron vector of \ $\ee^\tbB$, \ corresponding to
 the eigenvalue \ $r(\ee^\tbB) = \ee^{s(\tbB)}$ \ of \ $\ee^\tbB$ \ such that
 \ $\tbu^\top \bu = 1$.
\ Moreover, we have
 \begin{align*}
  \ee^{-s(\tbB)t} \ee^{t\tbB} \to \tbu \bu^\top \in \RR^{d \times d}_{++} \qquad
  \text{as \ $t \to \infty$,}
 \end{align*}
 and there exist \ $C_1, C_2, C_3 \in \RR_{++}$ \ such that
 \begin{equation}\label{C}
  \|\ee^{-s(\tbB)t} \ee^{t\tbB} - \tbu \bu^\top\| \leq C_1 \ee^{-C_2 t} , \qquad
  \|\ee^{t\tbB}\| \leq C_3 \ee^{s(\tbB)t} , \qquad t \in \RR_+ .
 \end{equation}
These Frobenius and Perron type results can be found, e.g., in Barczy and Pap
 \cite[Appendix A]{BarPap}.

\begin{Def}\label{Def_indecomposable_crit}
Let \ $(\bX_t)_{t\in\RR_+}$ \ be an irreducible multi-type CBI process with
 parameters \ $(d, \bc, \Bbeta, \bB, \nu, \bmu)$ \ such that
 \ $\EE(\|\bX_0\|) < \infty$ \ and the moment condition
 \eqref{moment_condition_m_new} holds.
Then \ $(\bX_t)_{t\in\RR_+}$ \ is called
 \[
   \begin{cases}
    subcritical & \text{if \ $s(\tbB) < 0$,} \\
    critical & \text{if \ $s(\tbB) = 0$,} \\
    supercritical & \text{if \ $s(\tbB) > 0$.}
   \end{cases}
 \]
\end{Def}
For motivations of Definitions \ref{Def_irreducible} and
 \ref{Def_indecomposable_crit}, see Barczy and Pap \cite[Section 3]{BarPap}.
Here we only point out that our classification of multi-type CBI processes is based
 on the asymptotic behaviour of \ $\EE(\bX_t)$ \ as \ $t \to \infty$, \ and this
 asymptotics is available at the moment only under the assumption of irreducibility
 of \ $(\bX_t)_{t\in\RR_+}$.

\section{Main results}
\label{section_CBI_as_growth}

First we present almost sure and \ $L_1$-convergence results for supercritical and
 irreducible multi-type CBI processes.

\begin{Thm}\label{convCBIasL1}
Let \ $(\bX_t)_{t\in\RR_+}$ \ be a supercritical and irreducible multi-type CBI
 process with parameters \ $(d, \bc, \Bbeta, \bB, \nu, \bmu)$ \ such that
 \ $\EE(\|\bX_0\|) < \infty$ \ and the moment condition
 \eqref{moment_condition_m_new} holds.
Then, there exists a non-negative random variable \ $w_{\bu,\bX_0}$ \ with
 \ $\EE(w_{\bu,\bX_0}) < \infty$ \ such that
 \begin{equation}\label{convwu}
  \ee^{-s(\tbB)t} \langle \bu, \bX_t \rangle \as w_{\bu,\bX_0}
  \qquad \text{as \ $t \to \infty$.}
 \end{equation}

Moreover, for each \ $\lambda \in \sigma(\tbB)$ \ such that
 \ $\Re(\lambda) \in \bigl(\frac{1}{2} s(\tbB), s(\tbB)\bigr]$ \ and the
 moment condition
 \begin{equation}\label{moment_condition_xlogx}
  \sum_{\ell=1}^d
   \int_{\cU_d} g(\|\bz\|) \bbone_{\{\|\bz\|\geq1\}} \, \mu_\ell(\dd \bz)
  < \infty
 \end{equation}
 with
 \[
   g(x)
   := \begin{cases}
       x^{\frac{s(\tbB)}{\Re(\lambda)}}
        & \text{if \ $\Re(\lambda)
                      \in \bigl(\frac{1}{2} s(\tbB), s(\tbB)\bigr)$,} \\
       x \log(x)
        & \text{if \ $\Re(\lambda) = s(\tbB)$
                \ ($\Longleftrightarrow$ $\lambda = s(\tbB)$),}
       \end{cases}
   \qquad x \in \RR_{++}
 \]
 holds, and for each left eigenvector \ $\bv \in \CC^d$ \ of \ $\tbB$
 \ corresponding to the eigenvalue \ $\lambda$, \ there exists a complex random
 variable \ $w_{\bv,\bX_0}$ \ with \ $\EE(|w_{\bv,\bX_0}|) < \infty$ \ such that
 \begin{equation}\label{convwv}
  \ee^{-\lambda t} \langle \bv, \bX_t \rangle \to w_{\bv,\bX_0} \qquad
  \text{as \ $t \to \infty$ \ in \ $L_1$ \ and almost surely,}
 \end{equation}
 and
 \begin{equation}\label{reprCBIw}
  \begin{aligned}
   w_{\bv,\bX_0}
   &\ase
    \langle \bv, \bX_0 \rangle
    + \frac{\langle \bv, \tBbeta \rangle}{\lambda}
    + \sum_{\ell=1}^d
       \langle\bv, \be_\ell\rangle
       \int_0^\infty
        \ee^{-\lambda u} \sqrt{2 c_\ell X_{u,\ell}} \, \dd W_{u,\ell} \\
   &\phantom{\ase}
    + \sum_{\ell=1}^d
       \int_0^\infty \int_{\cU_d} \int_{\cU_1}
        \ee^{-\lambda u} \langle\bv, \bz\rangle \bbone_{\{w\leq X_{u-,\ell}\}}
        \, \tN_\ell(\dd u, \dd\bz,  \dd w) \\
   &\phantom{\ase}
    + \int_0^\infty \int_{\cU_d}
       \ee^{-\lambda u} \langle\bv, \br\rangle \, \tM(\dd u, \dd\br) ,
  \end{aligned}
 \end{equation}
 where the improper integrals are convergent in \ $L_1$ \ and almost surely.
Especially,
 \ $\EE(w_{\bv,\bX_0}) = \langle\bv, \EE(\bX_0) + \lambda^{-1} \tBbeta\rangle$.
 Particularly, if \ $\langle\bv, \EE(\bX_0) + \lambda^{-1} \tBbeta\rangle
 \ne 0$, \ then \  $\PP(w_{\bv,\bX_0} = 0) <  1$.
\  Further, \ $w_{\bu,\bX_0} \ase 0$ \ if and only if \ $\bX_0 \ase \bzero$
 \ and \ $\tBbeta = \bzero$ \ (equivalently, \ $\bX_0 \ase \bzero$,
 \ $\Bbeta = \bzero$ \ and \ $\nu = 0$).

If the moment condition \eqref{moment_condition_xlogx} does not hold for
 \ $\lambda = s(\tbB)$, \ then \ $\ee^{-s(\tbB)t} \langle \bu, \bX_t\rangle \as 0$
 \ as \ $t \to \infty$, \ i.e., \ $\PP(w_{\bu,\bX_0} = 0) = 1$.

If the moment condition \eqref{moment_condition_xlogx} does not hold for \ $\lambda = s(\tbB)$,
 \ then \ $\ee^{-s(\tbB)t} \langle \bu, \bX_t\rangle$ \ does not converge in \ $L_1$ \ as \ $t\to\infty$,
 \ provided that \ $\PP(\bX_0 = \bzero)<1$ \ or \ $\tBbeta \ne \bzero$. \
If \ $\PP(\bX_0 = \bzero)=1$ \ and \ $\tBbeta = \bzero$, \ then \ $\PP(\bX_t=\bzero)=1$ \ for all \ $t\in\RR_+$.
\end{Thm}

Note that the asymptotic behavior of the second moment
 \ $\EE(|\langle\bv, \bX_t\rangle|^2)$ \ as \ $t \to \infty$ \ explains the role of
 the assumption \ $\Re(\lambda) \in \bigl(\frac{1}{2} s(\tbB), s(\tbB)\bigr]$ \ in
 Theorem \ref{convCBIasL1}, see Proposition \ref{second_moment_asymptotics_CBI}.

\noindent
\textbf{Proof of Theorem \ref{convCBIasL1}.}
By Lemma \ref{mart}, the process
 \ $\bigl(\ee^{-t\tbB} \bX_t
          - \int_0^t \ee^{-u\tbB} \tBbeta \, \dd u\bigr)_{t\in\RR_+}$
 \ is a martingale with respect to the filtration \ $(\cF_t^\bX)_{t\in\RR_+}$.
\ Moreover, for each \ $t \in \RR_+$, \ we have
 \begin{equation}\label{sub}
  \begin{aligned}
   \ee^{-s(\tbB)t} \langle \bu, \bX_t\rangle
  &= \ee^{-s(\tbB)t} \bu^\top \bX_t
   = \bu^\top \ee^{-t\tbB} \bX_t \\
  &= \bu^\top
     \biggl(\ee^{-t\tbB} \bX_t - \int_0^t \ee^{-v\tbB} \tBbeta \, \dd v\biggr)
     + \bu^\top \int_0^t \ee^{-v\tbB} \tBbeta \, \dd v \\
  &= \bu^\top
     \biggl(\ee^{-t\tbB} \bX_t - \int_0^t \ee^{-v\tbB} \tBbeta \, \dd v\biggr)
     + \langle \bu, \tBbeta \rangle \int_0^t \ee^{-s(\tbB)v} \, \dd v ,
  \end{aligned}
 \end{equation}
 where the function
 \ $\RR_+ \ni t
    \mapsto \langle \bu, \tBbeta \rangle \int_0^t \ee^{-s(\tbB)v} \, \dd v
    \in \RR_+$
 \ is increasing, since \ $\bu \in \RR_{++}^d$ \ and \ $\tBbeta \in \RR_+^d$.
\ Consequently, \ $(\ee^{-s(\tbB)t} \langle \bu, \bX_t \rangle)_{t\in\RR_+}$ \ is a
 submartingale with respect to the filtration \ $(\cF_t^\bX)_{t\in\RR_+}$.
Due to Theorem 4.6 in Barczy et al.\ \cite{BarLiPap2}, \ $(\bX_t)_{t\in\RR_+}$ \ and hence
 \ $(\ee^{-s(\tbB)t} \langle \bu, \bX_t \rangle)_{t\in\RR_+}$ \ have c\`adl\`ag sample paths almost surely.
Using again \ $\bu \in \RR_{++}^d$ \ and \eqref{sub}, we get
 \begin{align*}
  \EE(|\ee^{-s(\tbB)t} \langle \bu, \bX_t \rangle|)
  &= \EE(\ee^{-s(\tbB)t} \langle \bu, \bX_t \rangle)
   = \EE(\langle \bu, \bX_0 \rangle)
     + \langle \bu, \tBbeta \rangle \int_0^t \ee^{-s(\tbB)v} \, \dd v \\
  &\leq \|\bu\| \EE(\|\bX_0\|)
        + \langle \bu, \tBbeta \rangle \int_0^\infty \ee^{-s(\tbB)v} \, \dd v
   = \|\bu\| \EE(\|\bX_0\|) + \frac{\langle \bu, \tBbeta \rangle}{s(\tbB)}
 \end{align*}
 for all \ $t\in\RR_+$, \ thus we conclude
 \ $\sup_{t\in\RR_+} \EE(|\ee^{-s(\tbB)t} \langle \bu, \bX_t \rangle|) < \infty$.
\ Hence, by the submartingale convergence theorem, there exists a non-negative
 random variable \ $w_{\bu,\bX_0}$ \ with \ $\EE(w_{\bu,\bX_0}) < \infty$ \ such
 that \eqref{convwu} holds.

If \ $\lambda \in \sigma(\tbB)$ \ such that
 \ $\Re(\lambda) \in \bigl(\frac{1}{2} s(\tbB), s(\tbB)\bigr]$ \ and the
 moment condition \eqref{moment_condition_xlogx} holds, and \ $\bv \in \CC^d$ \ is
 a left eigenvector of \ $\tbB$ \ corresponding to the eigenvalue \ $\lambda$,
 \ then first we show the \ $L_1$-convergence of
 \ $\ee^{-\lambda t} \langle\bv, \bX_t\rangle$ \ as \ $t \to \infty$ \ towards the
 right hand side of \eqref{reprCBIw} together with the \ $L_1$-convergence of the
 improper integrals in \eqref{reprCBIw}.
Note that the condition
 \ $\Re(\lambda) \in \bigl(\frac{1}{2} s(\tbB), s(\tbB)\bigr]$ \ yields
 \ $\Re(\lambda) > 0$, \ so \ $\lambda \ne 0$.
\ For each \ $t \in \RR_+$, \ by Lemma \ref{SDE_transform_sol}, we have the
 representation
 \begin{equation}\label{reprZ}
  \ee^{-\lambda t} \langle\bv, \bX_t\rangle
  = \langle\bv, \bX_0\rangle + Z_t^{(1)} + Z_t^{(2)} + Z_t^{(3)} + Z_t^{(4)}
    + Z_t^{(5)}
 \end{equation}
 with
 \begin{align*}
  Z_t^{(1)}
  &:= \langle\bv, \tBbeta\rangle \int_0^t \ee^{-\lambda u} \, \dd u , \\
  Z_t^{(2)}
  &:= \sum_{\ell=1}^d
       \langle\bv, \be_\ell\rangle
       \int_0^t
        \ee^{-\lambda u} \sqrt{2 c_\ell X_{u,\ell}} \, \dd W_{u,\ell} , \\
  Z_t^{(3)}
  &:= \sum_{\ell=1}^d
       \int_0^t \int_{\cU_d} \int_{\cU_1}
        \ee^{-\lambda u} \langle\bv, \bz\rangle
        \bbone_{\{\|\bz\|<\ee^{\Re(\lambda)u}\}} \bbone_{\{w\leq X_{u-,\ell}\}}
        \, \tN_\ell(\dd u, \dd\bz, \dd w) ,
 \end{align*}
 \begin{align*}
  Z_t^{(4)}
  &:= \sum_{\ell=1}^d
       \int_0^t \int_{\cU_d} \int_{\cU_1}
        \ee^{-\lambda u} \langle\bv, \bz\rangle
        \bbone_{\{\|\bz\|\geq\ee^{\Re(\lambda)u}\}} \bbone_{\{w\leq X_{u-,\ell}\}}
        \, \tN_\ell(\dd u, \dd\bz, \dd w) , \\
  Z_t^{(5)}
  &:= \int_0^t \int_{\cU_d}
       \ee^{-\lambda u} \langle\bv, \br\rangle \, \tM(\dd u, \dd\br) .
 \end{align*}
Hence the \ $L_1$-convergence of \ $\ee^{-\lambda t} \langle\bv, \bX_t\rangle$ \ as
 \ $t \to \infty$ \ towards the right hand side of \eqref{reprCBIw} together with
 the \ $L_1$-convergence of the improper integrals in \eqref{reprCBIw} will follow
 from the convergences \ $D_t^{(j)} \mean 0$ \ as \ $t \to \infty$ \ for every
 \ $j \in \{1, 2, 3, 4, 5\}$ \ with
 \begin{align*}
  D_t^{(1)}
  &:= \frac{\langle\bv, \tBbeta\rangle}{\lambda}
      - \langle\bv, \tBbeta\rangle \int_0^t \ee^{-\lambda u} \, \dd u , \\
  D_t^{(2)}
  &:= \sum_{\ell=1}^d
       \langle\bv, \be_\ell\rangle
       \int_t^\infty
        \ee^{-\lambda u} \sqrt{2 c_\ell X_{u,\ell}} \, \dd W_{u,\ell} ,\\
  D_t^{(3)}
  &:= \sum_{\ell=1}^d
       \int_t^\infty \int_{\cU_d} \int_{\cU_1}
        \ee^{-\lambda u} \langle\bv, \bz\rangle
        \bbone_{\{\|\bz\|<\ee^{\Re(\lambda)u}\}} \bbone_{\{w\leq X_{u-,\ell}\}}
        \, \tN_\ell(\dd u, \dd\bz, \dd w) , \\
  D_t^{(4)}
  &:= \sum_{\ell=1}^d
       \int_t^\infty \int_{\cU_d} \int_{\cU_1}
        \ee^{-\lambda u} \langle\bv, \bz\rangle
        \bbone_{\{\|\bz\|\geq\ee^{\Re(\lambda)u}\}} \bbone_{\{w\leq X_{u-,\ell}\}}
        \, \tN_\ell(\dd u, \dd\bz, \dd w) ,\\
  D_t^{(5)}
  &:= \int_t^\infty \int_{\cU_d}
       \ee^{-\lambda u} \langle\bv, \br\rangle \, \tM(\dd u, \dd\br)
 \end{align*}
 for \ $t \in \RR_+$.
\ We have
 \begin{equation}\label{D_1_conv}
  D_t^{(1)} = \langle\bv, \tBbeta\rangle \int_t^\infty \ee^{-\lambda u} \, \dd u
            \to 0 \qquad \text{as \ $t \to \infty$.}
 \end{equation}
Moreover, for each \ $t \in \RR_+$, \ we have
 \[
   \EE\biggl(\int_t^\infty
              |\ee^{-\lambda u}|^2 2 c_\ell X_{u,\ell} \, \dd u \biggr)
   = 2 c_\ell
        \int_t^\infty \ee^{-2\Re(\lambda)u} \EE(X_{u,\ell}) \, \dd u .
 \]
By formulae \eqref{EXcond} and \eqref{C}, for each \ $v \in \RR_+$ \ and
 \ $\ell \in \{1, \ldots, d\}$, \ we get
 \begin{equation}\label{EX_bound}
  \begin{aligned}
   \EE(X_{v,\ell})
   &= \EE(\be_\ell^\top \bX_v)
    = \EE\biggl(\be_\ell^\top \ee^{v\tbB} \bX_0
              + \be_\ell^\top \int_0^v \ee^{u\tbB} \tBbeta \, \dd u\biggr) \\
   &\leq \|\ee^{v\tbB}\| \EE(\|\bX_0\|)
         + \|\tBbeta\| \int_0^v \|\ee^{u\tbB}\| \, \dd u \\
   &\leq C_3 \ee^{s(\tbB)v} \EE(\|\bX_0\|)
         + C_3 \|\tBbeta\| \int_0^v \ee^{s(\tbB)u} \, \dd u \\
   &= C_3 \ee^{s(\tbB)v} \EE(\|\bX_0\|)
      + C_3 \|\tBbeta\| \frac{\ee^{s(\tbB)v}-1}{s(\tbB)}
    \leq C_4 \ee^{s(\tbB)v}
  \end{aligned}
 \end{equation}
 with \ $C_4 := C_3 \EE(\|\bX_0\|) + \frac{C_3\|\tBbeta\|}{s(\tbB)}$.
\ By \eqref{EX_bound}, for each \ $t \in \RR_+$ \ and
 \ $\ell \in \{1, \ldots, d\}$, \ we obtain
 \begin{align*}
  &\EE\biggl(\int_t^\infty
              |\ee^{-\lambda u}|^2 2 c_\ell X_{u,\ell} \, \dd u \biggr)
   \leq 2 C_4 c_\ell
        \int_t^\infty \ee^{-2\Re(\lambda)u} \ee^{s(\tbB)u} \, \dd u \\
  &= 2 C_4 c_\ell
        \int_t^\infty \ee^{-(2\Re(\lambda)-s(\tbB))u} \, \dd u
   = \frac{2C_4c_\ell}{2\Re(\lambda)-s(\tbB)} \ee^{-(2\Re(\lambda)-s(\tbB))t}
   < \infty ,
 \end{align*}
 thus, by the independence of
 \ $(W_{t,1})_{t\in\RR_+}$, \ldots, $(W_{t,d})_{t\in\RR_+}$ \ and It\^{o}'s isometry for It\^{o}'s integrals
  (see, e.g., Ikeda and Watanabe \cite[Chapter II, Proposition 2.2]{IkeWat}),
  \begin{align*}
   &\EE(|D_t^{(2)}|^2)
    = \EE\Biggl(\sum_{\ell=1}^d
                \int_t^\infty
                 \Re(\langle \bv, \be_\ell\rangle \ee^{-\lambda u})
                 \sqrt{2 c_\ell X_{u,\ell}}
                 \, \dd W_{u,\ell}\Biggr)^2 \\
   &\phantom{\EE(|D_t^{(2)}|^2)=}
      + \EE\Biggl(\sum_{\ell=1}^d
                   \int_t^\infty
                    \Im(\langle \bv, \be_\ell\rangle \ee^{-\lambda u})
                    \sqrt{2 c_\ell X_{u,\ell}}
                    \, \dd W_{u,\ell}\Biggr)^2\\
   &= \sum_{\ell=1}^d
       \int_t^\infty
        (\Re(\langle \bv, \be_\ell\rangle \ee^{-\lambda u}))^2
        2 c_\ell \EE(X_{u,\ell})
        \, \dd u
      + \sum_{\ell=1}^d
         \int_t^\infty
          (\Im(\langle \bv, \be_\ell\rangle \ee^{-\lambda u}))^2
          2 c_\ell \EE(X_{u,\ell})
          \, \dd u \\
   &= \sum_{\ell=1}^d
       \int_t^\infty
        |\langle \bv, \be_\ell\rangle \ee^{-\lambda u}|^2 2 c_\ell \EE(X_{u,\ell})
        \, \dd u
     = 2 \sum_{\ell=1}^d
         |\langle\bv, \be_\ell\rangle|^2 c_\ell
         \int_t^\infty |\ee^{-\lambda u}|^2 \EE(X_{u,\ell}) \, \dd u
  \end{align*}
  \begin{align*}
   &\leq 2 C_4 \|\bv\|^2
         \sum_{\ell=1}^d
          c_\ell \int_t^\infty \ee^{-2\Re(\lambda)u} \ee^{s(\tbB)u} \, \dd u
    = \frac{2C_4\|\bv\|^2}{2\Re(\lambda)-s(\tbB)}
      \Biggl(\sum_{\ell=1}^d c_\ell\Biggr)
      \ee^{-(2\Re(\lambda)-s(\tbB))t} .
  \end{align*}
Consequently, we have
 \begin{equation}\label{D_2_conv_L2}
  D_t^{(2)} \qmean 0 \qquad \text{as \ $t \to \infty$,}
 \end{equation}
 hence we conclude
 \begin{equation}\label{D_2_conv_L1}
  D_t^{(2)} \mean 0 \qquad \text{as \ $t \to \infty$.}
 \end{equation}
By \eqref{EX_bound}, for each \ $t \in \RR_+$, \ we have
 \begin{align*}
  &\sum_{\ell=1}^d
      \EE\biggl(\int_t^\infty \int_{\cU_d} \int_{\cU_1}
                 |\ee^{-\lambda u}|^2 |\langle\bv, \bz\rangle|^2
                 \bbone_{\{\|\bz\|<\ee^{\Re(\lambda)u}\}}
                 \bbone_{\{w\leq X_{u,\ell}\}}
                 \, \dd u \, \mu_\ell(\dd\bz) \, \dd w\biggr) \\
  &= \sum_{\ell=1}^d
      \int_t^\infty \int_{\cU_d}
       \ee^{-2\Re(\lambda)u} |\langle\bv, \bz\rangle|^2
       \bbone_{\{\|\bz\|<\ee^{\Re(\lambda)u}\}} \EE(X_{u,\ell})
       \, \dd u \, \mu_\ell(\dd\bz)
   \leq C_4 \|\bv\|^2 K_t^{(3)}
 \end{align*}
 with
 \begin{equation}\label{D3_bound}
   K_t^{(3)} := \sum_{\ell=1}^d
                \int_t^\infty \int_{\cU_d}
                 \ee^{-(2\Re(\lambda)-s(\tbB))u} \|\bz\|^2
                 \bbone_{\{\|\bz\|<\ee^{\Re(\lambda)u}\}}
                 \, \dd u \, \mu_\ell(\dd\bz)
   \leq K_0^{(3)} .
 \end{equation}
We show that \ $K_0^{(3)} < \infty$.
\ For each \ $\ell \in \{1, \ldots, d\}$, \ using Fubini's theorem, we obtain
 \begin{align*}
  &\int_0^\infty \int_{\cU_d}
    \ee^{-(2\Re(\lambda)-s(\tbB))u} \|\bz\|^2 \bbone_{\{\|\bz\|<1\}}
    \, \dd u \, \mu_\ell(\dd\bz) \\
  &= \int_0^\infty \ee^{-(2\Re(\lambda)-s(\tbB))u} \, \dd u
     \int_{\cU_d} \|\bz\|^2 \bbone_{\{\|\bz\|<1\}} \, \mu_\ell(\dd\bz) \\
  &= \frac{1}{2\Re(\lambda)-s(\tbB)}
     \int_{\cU_d}
      \|\bz\|^2 \bbone_{\{\|\bz\|<1\}} \, \mu_\ell(\dd\bz)
   < \infty
 \end{align*}
 by Definition \ref{Def_admissible}, and
 \begin{align*}
  &\int_0^\infty \int_{\cU_d}
    \ee^{-(2\Re(\lambda)-s(\tbB))u} \|\bz\|^2
    \bbone_{\{1\leq\|\bz\|<\ee^{\Re(\lambda)u}\}}
    \, \dd u \, \mu_\ell(\dd\bz) \\
  &= \int_{\cU_d}
      \biggl(\int_{\frac{1}{\Re(\lambda)}\log(\|\bz\|)}^\infty
              \ee^{-(2\Re(\lambda)-s(\tbB))u} \, \dd u\biggr)
      \|\bz\|^2 \bbone_{\{\|\bz\|\geq1\}}
      \, \mu_\ell(\dd\bz) \\
  &= \int_{\cU_d}
      \frac{1}{2\Re(\lambda)-s(\tbB)}
      \|\bz\|^{-\frac{2\Re(\lambda)-s(\tbB)}{\Re(\lambda)}}
      \|\bz\|^2 \bbone_{\{\|\bz\|\geq1\}}
      \, \mu_\ell(\dd\bz) \\
  &= \frac{1}{2\Re(\lambda)-s(\tbB)}
     \int_{\cU_d} \|\bz\|^{\frac{s(\tbB)}{\Re(\lambda)}}
      \bbone_{\{\|\bz\|\geq1\}} \, \mu_\ell(\dd\bz)
   < \infty
 \end{align*}
 by the moment condition \eqref{moment_condition_xlogx} (in case of
 \ $\Re(\lambda) \in \bigl(\frac{1}{2} s(\tbB), s(\tbB)\bigr)$) \ and by
 Definition \ref{Def_admissible} (in case of \ $\Re(\lambda) = s(\tbB)$ \ or
 equivalently \ $\lambda = s(\tbB)$).
Thus we obtain \ $K_0^{(3)} < \infty$.
\ Consequently, by page 63 in Ikeda and Watanabe \cite{IkeWat}, for each
 \ $t \in \RR_+$, \ we conclude
 \begin{align*}
  \EE(|D_t^{(3)}|^2)
  &= \sum_{\ell=1}^d
      \EE\biggl(\int_t^\infty \int_{\cU_d} \int_{\cU_1}
                 |\ee^{-\lambda u}|^2 |\langle\bv, \bz\rangle|^2
                 \bbone_{\{\|\bz\|<\ee^{\Re(\lambda)u}\}}
                 \bbone_{\{w\leq X_{u,\ell}\}}
                 \, \dd u \, \mu_\ell(\dd\bz) \, \dd w\biggr) \\
  &\leq C_4 \|\bv\|^2 K_t^{(3)} <\infty.
 \end{align*}
We have
 \[
   K_t^{(3)}
   = K_0^{(3)} - \sum_{\ell=1}^d
     \int_0^t \int_{\cU_d}
      \ee^{-(2\Re(\lambda)-s(\tbB))u} \|\bz\|^2
        \bbone_{\{\|\bz\|<\ee^{\Re(\lambda)u}\}}
        \, \dd u \, \mu_\ell(\dd\bz)
 \]
 yielding
 \begin{equation}\label{K_3_conv}
  K_t^{(3)} \to 0 \qquad \text{as \ $t \to \infty$,}
 \end{equation}
 thus \ $\EE(|D_t^{(3)}|^2) \to 0$ \ as \ $t \to \infty$.
\ This implies \ $D_t^{(3)} \qmean 0$ \ as \ $t \to \infty$, \ thus we conclude
 \begin{equation}\label{D_3_conv_L1}
  D_t^{(3)} \mean 0 \qquad \text{as \ $t \to \infty$.}
 \end{equation}
Further, for each \ $t \in \RR_+$, \ we get
 \begin{align*}
  D_t^{(4)}
  &= \sum_{\ell=1}^d
      \int_t^\infty \int_{\cU_d} \int_{\cU_1}
       \ee^{-\lambda u} \langle\bv, \bz\rangle
       \bbone_{\{\|\bz\|\geq\ee^{\Re(\lambda)u}\}}
       \bbone_{\{w\leq X_{u-,\ell}\}}
       \, N_\ell(\dd u, \dd\bz, \dd w) \\
  &\quad
     - \sum_{\ell=1}^d
        \int_t^\infty \int_{\cU_d} \int_{\cU_1}
         \ee^{-\lambda u} \langle\bv, \bz\rangle
         \bbone_{\{\|\bz\|\geq\ee^{\Re(\lambda)u}\}}
         \bbone_{\{w\leq X_{u,\ell}\}}
         \, \dd u \, \mu_\ell(\dd\bz) \, \dd w .
 \end{align*}
Indeed, for each \ $t \in \RR_+$ \ and \ $\ell \in \{1, \ldots, d\}$, \ we have
 \begin{align*}
   &\left\vert\int_t^\infty \int_{\cU_d} \int_{\cU_1}
         \ee^{-\lambda u} \langle\bv, \bz\rangle
         \bbone_{\{\|\bz\|\geq\ee^{\Re(\lambda)u}\}}
         \bbone_{\{w\leq X_{u,\ell}\}}
         \, \dd u \, \mu_\ell(\dd\bz) \, \dd w \right\vert\\
   &\leq  \Vert\bv \Vert \int_t^\infty \int_{\cU_d} \int_{\cU_1}
              \ee^{-\Re(\lambda)u} \|\bz\|
              \bbone_{\{\|\bz\|\geq\ee^{\Re(\lambda)u}\}}
                    \bbone_{\{w\leq X_{u,\ell}\}}
                    \, \dd u \, \mu_\ell(\dd\bz) \, \dd w
                 \in\RR_+
 \end{align*}
 almost surely, since, by Fubini's theorem, \eqref{EX_bound} and \eqref{moment_condition_xlogx},
  we get
 \begin{align*}
  &\EE\biggl(\int_t^\infty \int_{\cU_d} \int_{\cU_1}
              \ee^{-\Re(\lambda)u} \|\bz\|
              \bbone_{\{\|\bz\|\geq\ee^{\Re(\lambda)u}\}}
                    \bbone_{\{w\leq X_{u,\ell}\}}
                    \, \dd u \, \mu_\ell(\dd\bz) \, \dd w\biggr) \\
  &= \int_t^\infty \int_{\cU_d}
          \ee^{-\Re(\lambda)u} \|\bz\| \bbone_{\{\|\bz\|\geq\ee^{\Re(\lambda)u}\}}
          \EE(X_{u,\ell})
          \, \dd u \, \mu_\ell(\dd\bz) \\
  &\leq C_4
        \int_{\cU_d}
         \biggl(\int_t^\infty
                 \ee^{(s(\tbB)-\Re(\lambda))u}
                 \bbone_{\left\{u\leq\frac{1}{\Re(\lambda)}\log(\|\bz\|)\right\}}
                 \, \dd u\biggr)
         \|\bz\| \bbone_{\{\|\bz\|\geq 1\}} \, \mu_\ell(\dd\bz)\\
  &\begin{cases}
    = C_4
      \int_{\cU_d}
       \Bigl(\frac{\log(\|\bz\|)}{s(\tbB)} - t\Bigr)
       \bbone_{\bigl\{\frac{\log(\|\bz\|)}{s(\tbB)}\geq t\bigr\}} \|\bz\|
       \bbone_{\{\|\bz\|\geq 1\}}
       \, \mu_\ell(\dd\bz) \\
    \leq \frac{C_4}{s(\tbB)}
         \int_{\cU_d}
           \|\bz\| \log(\|\bz\|) \bbone_{\{\|\bz\|\geq\ee^{s(\tbB)t}\}}
           \, \mu_\ell(\dd\bz)
    < \infty
     & \text{if \ $\Re(\lambda) = s(\tbB)$,} \\[1mm]
    \text{and} \\
    = C_4
      \int_{\cU_d}
       \frac{\|\bz\|^{\frac{s(\tbB)-\Re(\lambda)}{\Re(\lambda)}}
             - \ee^{(s(\tbB)-\Re(\lambda))t}}
            {s(\tbB)-\Re(\lambda)}
       \|\bz\| \bbone_{\bigl\{\frac{\log(\|\bz\|)}{\Re(\lambda)}\geq t\bigr\}}
       \, \mu_\ell(\dd\bz) \\[3mm]
    \leq \frac{C_4}{s(\tbB)-\Re(\lambda)}
         \int_{\cU_d}
          \|\bz\|^{\frac{s(\tbB)}{\Re(\lambda)}}
          \bbone_{\{\|\bz\|\geq\ee^{\Re(\lambda)t}\}}
          \, \mu_\ell(\dd\bz)
    < \infty
     & \text{if \ $\Re(\lambda) \in \bigl(\frac{1}{2} s(\tbB), s(\tbB)\bigr)$}
   \end{cases}
 \end{align*}
 by the moment condition \eqref{moment_condition_xlogx} and part (vi) of Definition \ref{Def_admissible}.
Consequently, we obtain
 \begin{equation}\label{D4_bound}
  \begin{aligned}
   |D_t^{(4)}|
   &\leq \|\bv\|
         \sum_{\ell=1}^d
          \int_t^\infty \int_{\cU_d} \int_{\cU_1}
           \ee^{-\Re(\lambda)u} \|\bz\| \bbone_{\{\|\bz\|\geq\ee^{\Re(\lambda)u}\}}
           \bbone_{\{w\leq X_{u-,\ell}\}}
           \, N_\ell(\dd u, \dd\bz, \dd w) \\
   &\quad
         + \|\bv\|
           \sum_{\ell=1}^d
            \int_t^\infty \int_{\cU_d} \int_{\cU_1}
             \ee^{-\Re(\lambda)u} \|\bz\|
             \bbone_{\{\|\bz\| \geq\ee^{\Re(\lambda)u}\}}
             \bbone_{\{w\leq X_{u,\ell}\}}
             \, \dd u \, \mu_\ell(\dd\bz) \, \dd w
    =: \|\bv\| K_t^{(4)}
  \end{aligned}
 \end{equation}
 almost surely, where, by page 62 in Ikeda and Watanabe \cite{IkeWat}, for each \ $t \in \RR_+$,
 \ we have
 \begin{align*}
  \EE(K_t^{(4)})
  &= 2 \sum_{\ell=1}^d
     \EE\biggl(\int_t^\infty \int_{\cU_d} \int_{\cU_1}
                \ee^{-\Re(\lambda)u} \|\bz\|
                \bbone_{\{\|\bz\|\geq\ee^{\Re(\lambda)u}\}}
                \bbone_{\{w\leq X_{u,\ell}\}}
                \, \dd u \, \mu_\ell(\dd\bz) \, \dd w\biggr) \\
  &\leq \begin{cases}
         \frac{2C_4}{s(\tbB)}
         \sum_{\ell=1}^d
          \int_{\cU_d}
           \|\bz\| \log(\|\bz\|) \bbone_{\{\|\bz\|\geq\ee^{s(\tbB)t}\}}
           \, \mu_\ell(\dd\bz)
          & \text{if \ $\Re(\lambda) = s(\tbB)$,} \\[2mm]
         \frac{2C_4}{s(\tbB)-\Re(\lambda)}
         \sum_{\ell=1}^d
          \int_{\cU_d}
           \|\bz\|^{\frac{s(\tbB)}{\Re(\lambda)}}
           \bbone_{\{\|\bz\|\geq\ee^{\Re(\lambda)t}\}}
           \, \mu_\ell(\dd\bz)
          & \text{if \ $\Re(\lambda)
                        \in \bigl(\frac{1}{2} s(\tbB), s(\tbB)\bigr)$.}
        \end{cases}
 \end{align*}
Thus the moment condition \eqref{moment_condition_xlogx} yields
 that \ $\EE(K_t^{(4)}) < \infty$ \ for all \ $t \in \RR_+$ \ and, by the dominated convergence theorem,
 \begin{equation}\label{K_4_conv}
  K_t^{(4)} \mean 0 \qquad \text{as \ $t \to \infty$,}
 \end{equation}
 and hence
 \begin{equation}\label{D_4_conv_L1}
  D_t^{(4)} \mean 0 \qquad \text{as \ $t \to \infty$.}
 \end{equation}
In a similar way, for each \ $t \in \RR_+$, \ we get
 \[
   D_t^{(5)}
   = \int_t^\infty \int_{\cU_d}
      \ee^{-\lambda u} \langle\bv, \br\rangle
      \, M(\dd u, \dd\br)
     - \int_t^\infty \int_{\cU_d}
        \ee^{-\lambda u} \langle\bv, \br\rangle
        \, \dd u \, \nu(\dd\br) ,
 \]
 since
 \begin{align*}
   \left\vert
   \int_t^\infty \int_{\cU_d}
        \ee^{-\lambda u} \langle\bv, \br\rangle
        \, \dd u \, \nu(\dd\br)
   \right\vert
  &\leq \|\bv\| \int_t^\infty \int_{\cU_d} \ee^{-\Re(\lambda)u} \|\br\| \, \dd u \, \nu(\dd\br)\\
  & = \frac{\|\bv\|}{\Re(\lambda)} \ee^{-\Re(\lambda)t}
     \int_{\cU_d} \|\br\| \, \nu(\dd\br)
   < \infty
 \end{align*}
 by the moment condition \eqref{moment_condition_m_new} and by Definition
 \ref{Def_admissible}.
 Consequently, we obtain
 \begin{equation}\label{D5_bound}
  \begin{aligned}
   |D_t^{(5)}|
   &\leq \|\bv\|
         \int_t^\infty \int_{\cU_d}
          \ee^{-\Re(\lambda)u} \|\br\| \, M(\dd u, \dd\br) \\
   &\quad
         + \|\bv\|
           \int_t^\infty \int_{\cU_d}
            \ee^{-\Re(\lambda)u} \|\br\| \, \dd u \, \nu(\dd\br)
    =: \|\bv\| K_t^{(5)}
  \end{aligned}
 \end{equation}
 almost surely, where, by page 62 in Ikeda and Watanabe \cite{IkeWat},
 for each \ $t \in \RR_+$, \ we obtain
 \[
   \EE(K_t^{(5)})
   = \frac{2}{\Re(\lambda)} \ee^{-\Re(\lambda)t}
     \int_{\cU_d} \|\br\| \, \nu(\dd\br) < \infty .
 \]
Hence we conclude
 \begin{equation}\label{K_5_conv}
  K_t^{(5)} \mean 0 \qquad \text{as \ $t \to \infty$,}
 \end{equation}
 implying
 \begin{equation}\label{D_5_conv_L1}
  D_t^{(5)} \mean 0 \qquad \text{as \ $t \to \infty$.}
 \end{equation}
The convergences \eqref{D_1_conv}, \eqref{D_2_conv_L1}, \eqref{D_3_conv_L1},
 \eqref{D_4_conv_L1} and \eqref{D_5_conv_L1} yield the \ $L_1$-convergence of
 \ $\ee^{-\lambda t} \langle\bv, \bX_t\rangle$ \ towards the right hand side of
 \eqref{reprCBIw} as \ $t \to \infty$ \ together with the \ $L_1$-convergence of
 the improper integrals in \eqref{reprCBIw}.
In fact, it turned out that \ $D^{(2)}_t$ \ and \ $D^{(3)}_t$ \ converge to
 \ $0$ \ in \ $L_2$ \ as \ $t\to\infty$, \ but \ $D^{(4)}_t$ \ and \ $D^{(5)}_t$
 \ converge to \ $0$ \ only in \ $L_1$ \ as \ $t\to\infty$.

Next we show the almost sure convergence of
 \ $\ee^{-\lambda t} \langle\bv, \bX_t\rangle$ \ as \ $t \to \infty$ \ together
 with the almost sure convergence of the improper integrals in \eqref{reprCBIw}.
For each \ $t \in \RR_+$, \ we use the representation \eqref{reprZ}.
We have
 \begin{equation}\label{Z_0_conv}
  Z_t^{(1)} \to \langle\bv, \tBbeta\rangle \int_0^\infty \ee^{-\lambda u} \, \dd u
            = \frac{\langle\bv, \tBbeta\rangle}{\lambda}
  \qquad \text{as \ $t \to \infty$.}
 \end{equation}
As in case of \ $(D_t^{(2)})_{t\in\RR_+}$, \ for each \ $t \in \RR_+$, \ one can
 derive
 \begin{equation*}
  \begin{aligned}
   \EE(|Z_t^{(2)}|^2)
   &= 2 \sum_{\ell=1}^d
        |\langle\bv, \be_\ell\rangle|^2 c_\ell
        \int_0^t |\ee^{-\lambda u}|^2 \EE(X_{u,\ell}) \, \dd u \\
   &\leq 2 C_4 \|\bv\|^2 \sum_{\ell=1}^d c_\ell
         \int_0^\infty \ee^{-2\Re(\lambda)u} \ee^{s(\tbB)u} \, \dd u
    = \frac{2C_4\|\bv\|^2}{2\Re(\lambda)-s(\tbB)} \sum_{\ell=1}^d c_\ell
    < \infty ,
  \end{aligned}
 \end{equation*}
 hence the real and imaginary parts of \ $(Z_t^{(2)})_{t\in\RR_+}$ \ are
 \ $L_2$-bounded martingales.
As in case of \ $(D_t^{(3)})_{t\in\RR_+}$, \ for each \ $t \in \RR_+$, \ one can
 derive
 \begin{equation*}
  \begin{aligned}
   &\EE(|Z_t^{(3)}|^2)
    = \sum_{\ell=1}^d
       \int_0^t \int_{\cU_d}
        \ee^{-2\Re(\lambda)u} |\langle\bv, \bz\rangle|^2
        \bbone_{\{\|\bz\|<\ee^{\Re(\lambda)u}\}} \EE(X_{u,\ell})
        \, \dd u \, \mu_\ell(\dd\bz) \\
   &\leq C_4 \|\bv\|^2
         \sum_{\ell=1}^d
          \int_0^\infty \int_{\cU_d}
          \ee^{-(2\Re(\lambda)-s(\tbB))u} \|\bz\|^2
          \bbone_{\{\|\bz\|<\ee^{\Re(\lambda)u}\}}
          \, \dd u \, \mu_\ell(\dd\bz)
    = C_4 \|\bv\|^2 K_0^{(3)} < \infty .
  \end{aligned}
 \end{equation*}
Consequently, the real and imaginary parts of \ $(Z_t^{(3)})_{t\in\RR_+}$ \ are
 \ $L_2$-bounded martingales.
As in case of \ $(D_t^{(4)})_{t\in\RR_+}$, \ for each \ $t \in \RR_+$, \ one can
 derive
 \[
   \EE(|Z_t^{(4)}|)
   \leq 2 \|\bv\|
        \sum_{\ell=1}^d
         \int_0^\infty \int_{\cU_d}
          \ee^{-\Re(\lambda)u} \|\bz\| \bbone_{\{\|\bz\|\geq\ee^{\Re(\lambda)u}\}}
          \EE(X_{u,\ell})
          \, \dd u \, \mu_\ell(\dd\bz)
   = \|\bv\| \EE(K_0^{(4)})
   < \infty ,
 \]
 hence the real and imaginary parts of \ $(Z_t^{(4)})_{t\in\RR_+}$ \ are
 \ $L_1$-bounded martingales.
As in case of \ $(D_t^{(5)})_{t\in\RR_+}$, \ for each \ $t \in \RR_+$, \ one can
 derive
 \[
   \EE(|Z_t^{(5)}|)
   \leq 2 \|\bv\|
        \int_0^\infty \int_{\cU_d}
         \ee^{-\Re(\lambda)u} \|\br\| \, \dd u \, \nu(\dd\br)
   = \|\bv\| \EE(K_0^{(5)})
   < \infty
 \]
 hence the real and
 imaginary parts of \ $(Z_t^{(5)})_{t\in\RR_+}$ \ are \ $L_1$-bounded martingales.
Consequently, by the martingale convergence theorem, we conclude that the real and
 imaginary parts of the martingales \ $(Z_t^{(j)})_{t\in\RR_+}$,
 \ $j \in \{2, 3, 4, 5\}$, \ are almost sure convergent \ as \ $t \to \infty$,
 \ hence, by \eqref{Z_0_conv}, we conclude the almost sure convergence of
 \ $\ee^{-\lambda t} \langle\bv, \bX_t\rangle$ \ as \ $t \to \infty$ \ together
 with the almost sure convergence of the improper integrals in \eqref{reprCBIw}.

We have already showed the \ $L_1$-convergence of
 \ $\ee^{-\lambda t} \langle\bv, \bX_t\rangle$ \ as \ $t \to \infty$ \ towards the
 right hand side of \eqref{reprCBIw}, so the almost sure convergence of
 \ $\ee^{-\lambda t} \langle\bv, \bX_t\rangle$ \ as \ $t \to \infty$ \ yields the almost sure convergence in \eqref{convwv} as well.

By the convergence
 \ $\ee^{-\lambda t} \langle\bv, \bX_t\rangle \mean w_{\bv,\bX_0}$ \ as
 \ $t \to \infty$ \ in \eqref{convwv}, we obtain
 \ $\EE(\ee^{-\lambda t} \langle\bv, \bX_t\rangle) \to \EE(w_{\bv,\bX_0})$ \ as
 \ $t \to \infty$.
\ On the other hand, for each \ $t \in \RR_+$, \ using the representation
 \eqref{reprZ} and the martingale property of the processes
 \ $(Z_t^{(j)})_{t\in\RR_+}$, \ $j \in \{2, 3, 4, 5\}$, \ we have
 \ $\EE(\ee^{-\lambda t} \langle\bv, \bX_t\rangle)
     = \langle\bv, \EE(\bX_0)\rangle + Z_t^{(1)}
     = \langle\bv, \EE(\bX_0)\rangle
       + \langle\bv, \tBbeta\rangle \int_0^t \ee^{-\lambda u} \, \dd u$,
 \ hence we obtain
  \ $\EE(\ee^{-\lambda t} \langle\bv, \bX_t\rangle)
     \to \langle\bv, \EE(\bX_0)\rangle
         + \frac{\langle\bv, \tBbeta\rangle}{\lambda}$
 \ as \ $t \to \infty$.
\ Consequently,
 \ $\EE(w_{\bv,\bX_0})
    = \langle\bv, \EE(\bX_0)\rangle
      + \frac{\langle\bv, \tBbeta\rangle}{\lambda}$.
\ From here, we can see that if
 \ $\langle\bv, \EE(\bX_0)\rangle + \frac{\langle\bv, \tBbeta\rangle}{\lambda}\neq 0$,
 \ then \ $\PP(w_{\bv,\bX_0}=0) < 1$.

Next, we prove that \ $w_{\bu,\bX_0} \ase 0$ \ if and only if
 \ $\bX_0 \ase \bzero$ \ and \ $\tBbeta = \bzero$.
\ Since \ $w_{\bu,\bX_0}$ \ is non-negative, we have \ $w_{\bu,\bX_0} \ase 0$
 \ if and only if \ $\EE(w_{\bu,\bX_0}) = 0$.
Since \ $\bu \in \RR_{++}^d$, \ $\PP(\bX_0 \in \RR_+^d) = 1$ \ and
 \ $\tBbeta \in \RR_+^d$, \ we have
 \ $\EE(w_{\bu,\bX_0})
    = \langle\bu, \EE(\bX_0)\rangle + \frac{\langle\bu,\tBbeta\rangle}{s(\tbB)}
    = 0$
 \ if and only if \ $\EE(\bX_0) = \bzero$ \ and \ $\tBbeta = \bzero$,
 \ yielding the assertion in question.

Next, we prove that if the moment condition \eqref{moment_condition_xlogx} does
 not hold for \ $\lambda = s(\tbB)$, \ then \ $\PP(w_{\bu,\bX_0} = 0) = 1$.
\ For each \ $t, T \in \RR_+$, \ by Lemma \ref{decomposition_CBI}, we have
 \ $\bX_{t+T} \distre \bX_t^{(1)} + \bX_t^{(2, T)}$, \ where
 \ $(\bX_s^{(1)})_{s\in\RR_+}$ \ and \ $(\bX_s^{(2,T)})_{s\in\RR_+}$ \ are
 independent multi-type CBI processes with \ $\PP(\bX_0^{(1)} = \bzero) = 1$,
 \ $\bX_0^{(2,T)} \distre \bX_T$, \ and with parameters
 \ $(d, \bc, \Bbeta, \bB, \nu, \bmu)$ \ and \ $(d, \bc, \bzero, \bB, 0, \bmu)$,
 \ respectively.
Taking the inner product with \ $\bu$ \ and multiplying by
 \ $\ee^{-s(\tbB)(t+T)}$, \ we obtain
 \[
   \ee^{-s(\tbB)(t+T)} \langle\bu, \bX_{t+T}\rangle
   \distre \ee^{-s(\tbB)T} (\ee^{-s(\tbB)t} \langle\bu, \bX_t^{(1)}\rangle)
           + \ee^{-s(\tbB)T} (\ee^{-s(\tbB)t} \langle\bu, \bX_t^{(2,T)}\rangle) ,
   \qquad t, T \in\RR_+ .
 \]
Letting \ $t \to \infty$, \ by \eqref{convwu}, we obtain
 \[
   w_{\bu,\bX_0}
   \distre \ee^{-s(\tbB)T}  w_{\bu,\bzero}^{(1)}
           + \ee^{-s(\tbB)T} w_{\bu,\bX_0^{(2,T)}}^{(2,T)}  , \qquad T \in \RR_+ ,
 \]
 where
 \ $w_{\bu,\bzero}^{(1)}
    := \lim_{t\to\infty} \ee^{-s(\tbB)t} \langle\bu, \bX_t^{(1)}\rangle$
 \ and
 \ $w_{\bu,\bX_0^{(2,T)}}^{(2,T)}
    := \lim_{t\to\infty} \ee^{-s(\tbB)t} \langle\bu, \bX_t^{(2,T)}\rangle$
 \ almost surely.
Due to Kyprianou et al.\ \cite[Theorem 1.3/(ii)]{KypPalRen} and the law of total
 probability, we have \ $\PP(w_{\bu,\bX_0^{(2,T)}}^{(2,T)}  = 0) = 1$, \ since
 \ $(\bX_s^{(2,T)})_{s\in\RR_+}$ \ is a CB process and the moment condition
 \eqref{moment_condition_xlogx} does not hold for \ $\lambda = s(\tbB)$.
\ Consequently, for each \ $T \in \RR_+$, \ we have
 \ $w_{\bu,\bX_0} \distre \ee^{-s(\tbB)T}  w_{\bu,\bzero}^{(1)}$.
\ The convergence \ $\ee^{-s(\tbB)T}  w_{\bu,\bzero}^{(1)} \as 0$ \ as
 \ $T \to \infty$ \ yields \ $\ee^{-s(\tbB)T}  w_{\bu,\bzero}^{(1)} \distr 0$ \ as
 \ $T \to \infty$, \ thus we conclude that \ $w_{\bu,\bX_0} \ase 0$, \ hence the
 proof is complete.

Finally, we prove that if the moment condition \eqref{moment_condition_xlogx} does not hold for \ $\lambda = s(\tbB)$,
 \ then \ $\ee^{-s(\tbB)t} \langle \bu, \bX_t\rangle$ \ does not converge in \ $L_1$ \ as \ $t\to\infty$,
 \ provided that \ $\PP(\bX_0 = \bzero)<1$ \ or \ $\tBbeta \ne \bzero$. \
On the contrary, let us suppose that \ $\ee^{-s(\tbB)t} \langle \bu, \bX_t\rangle$ \ converges in \ $L_1$ \ as \ $t\to\infty$.
Recall that if the moment condition \eqref{moment_condition_xlogx} does not hold for \ $\lambda = s(\tbB)$, \ then
 \ $\ee^{-s(\tbB)t} \langle \bu, \bX_t\rangle \as 0$ \ as \ $t\to\infty$, \ which yields that
 \ $\ee^{-s(\tbB)t} \langle \bu, \bX_t\rangle$ \ could converge only to \ $0$ \ in \ $L_1$ \ as \ $t\to\infty$.
\ Especially, \ $\EE(\ee^{-s(\tbB)t} \langle \bu, \bX_t\rangle)$ \  would converge to \ $0$ \ as \ $t\to\infty$.
\ Using \eqref{sub}, we have
 \[
   \EE(\ee^{-s(\tbB)t} \langle \bu, \bX_t\rangle)
     = \langle\bu,\EE(\bX_0)\rangle + \langle \bu, \tBbeta \rangle \int_0^t \ee^{-s(\tbB)v} \, \dd v
     \to \langle\bu,\EE(\bX_0) \rangle + \frac{\langle \bu, \tBbeta \rangle}{s(\tbB)}
     \qquad \text{as \ $t\to\infty$},
 \]
 so \ $\bu^\top\EE(\bX_0) + \frac{\langle \bu, \tBbeta \rangle}{s(\tbB)} = 0$ \ should hold.
Since \ $\bu\in\RR_{++}^d$, \ $\PP(\bX_0\in\RR_+^d)=1$ \ and \ $\tBbeta\in\RR_+^d$, \ this would imply that
 \ $\bX_0\ase \bzero$ \ and \ $\tBbeta=\bzero$ \ (equivalently, \ $\bX_0\ase \bzero$, \ $\Bbeta=\bzero$ \ and \ $\nu=0$)
 \ leading us to a contradiction.
\proofend

In the next remark we explain why we do not have any result in the case when the moment condition \eqref{moment_condition_xlogx}
 does not hold for \ $\lambda \in \big(\frac{1}{2}s(\tbB), s(\tbB)\big)$.

\begin{Rem}\label{Rem_lambda_not_sB}
If the moment condition \eqref{moment_condition_xlogx} does not hold for \ $\lambda \in \big(\frac{1}{2}s(\tbB), s(\tbB)\big)$, \ then
 we do not know whether \ $\ee^{-\lambda t} \langle \bv, \bX_t \rangle$ \ converges almost surely or not as \ $t\to\infty$,
 \ where \ $\bv \in \CC^d$ \ is a left eigenvector  of \ $\tbB$ \ corresponding to the eigenvalue \ $\lambda$.
Provided that it converges almost surely to a complex random variable \ $w_{\bv,\bX_0}$, \ then, similarly as in the proof
 of Theorem \ref{convCBIasL1}, we would have
 \[
   w_{\bv,\bX_0}
   \distre \ee^{-\lambda T}  w_{\bv,\bzero}^{(1)}
           + \ee^{-\lambda T} w_{\bv,\bX_0^{(2,T)}}^{(2,T)}  , \qquad T \in \RR_+ ,
 \]
 where
 \ $w_{\bv,\bzero}^{(1)}
    := \lim_{t\to\infty} \ee^{-\lambda t} \langle\bv, \bX_t^{(1)}\rangle$
 \ and
 \ $w_{\bv,\bX_0^{(2,T)}}^{(2,T)}
    := \lim_{t\to\infty} \ee^{-\lambda t} \langle\bv, \bX_t^{(2,T)}\rangle$
 \ almost surely, and \ $(\bX_s^{(1)})_{s\in\RR_+}$ \ and \ $(\bX_s^{(2,T)})_{s\in\RR_+}$ \ are
 independent multi-type CBI processes with \ $\PP(\bX_0^{(1)} = \bzero) = 1$, \ $\bX_0^{(2,T)} \distre \bX_T$,
 \ and with parameters \ $(d, \bc, \Bbeta, \bB, \nu, \bmu)$ \ and \ $(d, \bc, \bzero, \bB, 0, \bmu)$, \ respectively.
However, contrary to the proof of Theorem \ref{convCBIasL1}, we do not know whether \ $\PP(w_{\bv,\bX_0^{(2,T)}}^{(2,T)} = 0)=1$ \ holds or not.
In the proof of Theorem \ref{convCBIasL1} we used that the corresponding result for \ $\lambda = s(\tbB)$ \ holds for CB processes due to
 Kyprianou et al. \cite[Theorem 1.3/(ii)]{KypPalRen} which is based on a so-called spine decomposition technique.
Unfortunately, we do not know whether this technique could be adapted to the case of \ $\lambda \in \big(\frac{1}{2}s(\tbB), s(\tbB)\big)$ \ or not.
We also do not know if the moment condition \eqref{moment_condition_xlogx} does not hold for \ $\lambda \in \big(\frac{1}{2}s(\tbB), s(\tbB)\big)$, \
 then \ $\ee^{-\lambda t} \langle \bv, \bX_t \rangle$ \ converges in \ $L_1$ \ or not as \ $t\to\infty$.
In the proof of Theorem \ref{convCBIasL1}, the corresponding $L_1$-convergence in the case of \ $\lambda=s(\tbB)$ \ is based on the
 almost sure convergence of \ $\ee^{- s(\tbB)t} \langle \bu, \bX_t \rangle$ \ as \ $t\to\infty$.
The above mentioned questions remain open problems.
\proofend
\end{Rem}

\begin{Thm}\label{convCBIwasL1}
Let \ $(\bX_t)_{t\in\RR_+}$ \ be a supercritical and irreducible multi-type CBI
 process with parameters \ $(d, \bc, \Bbeta, \bB, \nu, \bmu)$ \ such that
 \ $\EE(\|\bX_0\|) < \infty$ \ and the moment condition
 \eqref{moment_condition_m_new} holds.
Then
 \[
   \ee^{-s(\tbB)t} \bX_t \to w_{\bu,\bX_0} \tbu \qquad
   \text{as \ $t \to \infty$ \ almost surely,}
 \]
 where \ $w_{\bu,\bX_0}$ \ is introduced in \eqref{convwv}.
If the moment condition \eqref{moment_condition_xlogx} does not hold for
 \ $\lambda = s(\tbB)$, \ then \ $w_{\bu,\bX_0} \ase 0$.
\ If the moment condition \eqref{moment_condition_xlogx} holds for
 \ $\lambda = s(\tbB)$, \ then \ $\ee^{-s(\tbB)t} \bX_t \to w_{\bu,\bX_0} \tbu$
 \ as \ $t\to\infty$ \ in \ $L_1$ \ as well.
If the moment condition \eqref{moment_condition_xlogx} does not hold for \ $\lambda = s(\tbB)$,
 \ then \ $\ee^{-s(\tbB)t}\bX_t$ \ does not converge in \ $L_1$ \ as \ $t\to\infty$,
 \ provided that \ $\PP(\bX_0 = \bzero)<1$ \ or \ $\tBbeta \ne \bzero$. \
If \ $\PP(\bX_0 = \bzero)=1$ \ and \ $\tBbeta = \bzero$, \ then \ $\PP(\bX_t=\bzero)=1$ \ for all \ $t\in\RR_+$.
\end{Thm}

\noindent
\textbf{Proof.}
First, let us suppose that the moment condition \eqref{moment_condition_xlogx} does
 not hold for \ $\lambda = s(\tbB)$.
\ Then for each \ $\ell \in \{1, \ldots, d\}$, \ by \eqref{convwu}, we have
 \[
   \limsup_{t\to\infty} \ee^{-s(\tbB)t} \be_\ell^\top \bX_t
   \leq \limsup_{t\to\infty}
         \frac{1}{u_\ell} \ee^{-s(\tbB)t} \langle\bu, \bX_t\rangle
   = 0
 \]
 almost surely, yielding that \ $\ee^{-s(\tbB)t} \bX_t \to \bzero$ \ as
 \ $t \to \infty$ \ almost surely.

In what follows, let us suppose that the moment condition
 \eqref{moment_condition_xlogx} holds for \ $\lambda = s(\tbB)$.
\ For each \ $t, T \in \RR_+$, \ put
 \begin{equation}\label{Delta_t_t+T}
  \bDelta_{t,t+T}
  := \ee^{-s(\tbB)(t+T)} \bX_{t+T} - \ee^{-s(\tbB)(t+T)} \ee^{T\tbB} \bX_t
   = \ee^{-s(\tbB)(t+T)} (\bX_{t+T} - \ee^{T\tbB} \bX_t) .
 \end{equation}
We are going to carry out the proof in several steps. As an
initial step, 	we show that for each \ $T \in \RR_+$, \ we have
\begin{equation}\label{Lemma_L1}
\bDelta_{t,t+T} \mean \bzero \qquad \text{as \ $t \to \infty$.}
\end{equation}
Using \eqref{Lemma_L1} we prove the \ $L_1$ \ convergence of
 \ $\ee^{-s(\tbB)t} \bX_t$ \ towards \ $w_{\bu,\bX_0} \tbu$ \ as \ $t\to\infty$ \ (see \eqref{convCBIwL1}).
Then we show the almost sure convergence of \  $\bDelta_{t,t+T}$ \ and that of
 \ $\ee^{-s(\tbB)t} \bX_t$ \ along lattice times (see \eqref{Lemma} and \eqref{lattice}).
Finally, we prove almost sure convergence
 of \ $\ee^{-s(\tbB)t} \bX_t$ \ towards \ $w_{\bu,\bX_0} \tbu$ \ as \ $t\to\infty$.

By Lemma \ref{SDE_transform_sol}, we obtain the representation
 \begin{equation}\label{reprJ}
  \bDelta_{t,t+T}
  = \bJ_{t,t+T}^{(1)} + \bJ_{t,t+T}^{(2)} + \bJ_{t,t+T}^{(3)} + \bJ_{t,t+T}^{(4)}
    + \bJ_{t,t+T}^{(5)}
 \end{equation}
 for all \ $t, T \in \RR_+$ \ with
 \begin{align*}
  \bJ_{t,t+T}^{(1)}
  &:= \ee^{-s(\tbB)(t+T)} \int_t^{t+T} \ee^{(t+T-v)\tbB} \tBbeta \, \dd v , \\
  \bJ_{t,t+T}^{(2)}
  &:= \ee^{-s(\tbB)(t+T)}
      \sum_{\ell=1}^d
       \int_t^{t+T}
        \ee^{(t+T-v)\tbB} \be_\ell \sqrt{2 c_\ell X_{v,\ell}}
        \, \dd W_{v,\ell} , \\
  \bJ_{t,t+T}^{(3)}
  &:= \ee^{-s(\tbB)(t+T)}
      \sum_{\ell=1}^d
       \int_t^{t+T} \int_{\cU_d} \int_{\cU_1}
        \ee^{(t+T-v)\tbB} \bz \bbone_{\{\|\bz\|<\ee^{s(\tbB)v}\}}
        \bbone_{\{w\leq X_{v-,\ell}\}}
        \, \tN_\ell(\dd v, \dd\bz, \dd w) ,\\
  \bJ_{t,t+T}^{(4)}
  &:= \ee^{-s(\tbB)(t+T)}
      \sum_{\ell=1}^d
       \int_t^{t+T} \int_{\cU_d} \int_{\cU_1}
        \ee^{(t+T-v)\tbB} \bz \bbone_{\{\|\bz\|\geq\ee^{s(\tbB)v}\}}
        \bbone_{\{w\leq X_{v-,\ell}\}}
        \, \tN_\ell(\dd v, \dd\bz, \dd w) , \\
  \bJ_{t,t+T}^{(5)}
  &:= \ee^{-s(\tbB)(t+T)}
      \int_t^{t+T} \int_{\cU_d} \ee^{(t+T-v)\tbB} \br \, \tM(\dd v, \dd\br)
 \end{align*}
 for all \ $t, T \in \RR_+$.
\ For each \ $t, T \in \RR_+$, \ we obtain
 \[
   \bJ_{t,t+T}^{(1)}
   = \ee^{-s(\tbB)(t+T)} \int_0^T \ee^{(T-u)\tbB} \tBbeta \, \dd u ,
 \]
 hence, for each \ $T \in \RR_+$, \ we conclude
 \begin{equation}\label{J_1_conv}
  \bJ_{t,t+T}^{(1)} \to \bzero \qquad \text{as \ $t \to \infty$.}
 \end{equation}
By \eqref{C} and \eqref{EX_bound}, for each \ $t, T \in \RR_+$ \ and
 \ $\ell \in \{1, \ldots, d\}$, \ we obtain
 \begin{align*}
  &\EE\biggl(\int_t^{t+T}
             \bigl\|\ee^{(t+T-v)\tbB} \be_\ell\bigr\|^2 2 c_\ell X_{v,\ell}
             \, \dd v \biggr)
   = 2 c_\ell
       \int_t^{t+T} \bigl\|\ee^{(t+T-v)\tbB} \be_\ell\bigr\|^2 \EE(X_{v,\ell})
       \, \dd v \\
  &\leq 2 C_3^2 C_4 c_\ell
        \int_t^{t+T} \ee^{2s(\tbB)(t+T-v)} \ee^{s(\tbB)v} \, \dd v \\
  &\leq 2 C_3^2 C_4 c_\ell
        \ee^{2s(\tbB)(t+T)} \int_t^\infty \ee^{-s(\tbB)v} \, \dd v
   = \frac{2C_3^2 C_4c_\ell}{s(\tbB)} \ee^{s(\tbB)(t+2T)}
   < \infty ,
 \end{align*}
 thus
 \begin{equation}\label{J_2_bound}
  \begin{aligned}
   \EE(\|\bJ_{t,t+T}^{(2)}\|^2)
   &= \ee^{-2s(\tbB)(t+T)}
      \sum_{\ell=1}^d
       \EE\biggl(\int_t^{t+T}
                  \bigl\|\ee^{(t+T-v)\tbB} \be_\ell\bigr\|^2 2 c_\ell X_{v,\ell}
                  \, \dd v \biggr) \\
   &\leq \frac{2C_3^2 C_4}{s(\tbB)} \ee^{-2s(\tbB)(t+T)}
         \sum_{\ell=1}^d c_\ell \ee^{s(\tbB)(t+2T)}
    = \frac{2 C_3^2 C_4}{s(\tbB)} \biggl(\sum_{\ell=1}^d c_\ell\biggr)
      \ee^{-s(\tbB)t} .
 \end{aligned}
 \end{equation}
Consequently, for each \ $T \in \RR_+$, \ we conclude
 \begin{equation}\label{J_2_conv_L2}
  \bJ_{t,t+T}^{(2)} \qmean \bzero \qquad \text{as \ $t \to \infty$,}
 \end{equation}
 and hence
 \begin{equation}\label{J_2_conv_L1}
  \bJ_{t,t+T}^{(2)} \mean \bzero \qquad \text{as \ $t \to \infty$.}
 \end{equation}
By \eqref{C}, \eqref{EX_bound} and \eqref{D3_bound}, for each \ $t, T \in \RR_+$ \ and
 \ $\ell \in \{1, \ldots, d\}$, \ we get
 \begin{align*}
  &\EE\biggl(\int_t^{t+T} \int_{\cU_d} \int_{\cU_1}
              \bigl\|\ee^{(t+T-v)\tbB} \bz\bigr\|^2
              \bbone_{\{\|\bz\|<\ee^{s(\tbB)v}\}}
              \bbone_{\{w\leq X_{v,\ell}\}}
              \, \dd v \, \mu_\ell(\dd\bz) \, \dd w\biggr) \\
  &= \int_t^{t+T} \int_{\cU_d}
      \bigl\|\ee^{(t+T-v)\tbB} \bz\bigr\|^2
      \bbone_{\{\|\bz\|<\ee^{s(\tbB)v}\}} \EE(X_{v,\ell})
      \, \dd v \, \mu_\ell(\dd\bz) \\
  &\leq C_3^2 C_4
        \int_t^{t+T} \int_{\cU_d}
         \ee^{s(\tbB)2(t+T-v)} \|\bz\|^2 \bbone_{\{\|\bz\|<\ee^{s(\tbB)v}\}}
         \ee^{s(\tbB)v}
         \, \dd v \, \mu_\ell(\dd\bz) \\
  &\leq C_3^2 C_4 \ee^{2s(\tbB)(t+T)}
        \int_t^\infty \int_{\cU_d}
         \ee^{-s(\tbB)v} \|\bz\|^2 \bbone_{\{\|\bz\|<\ee^{s(\tbB)v}\}}
         \, \dd v \, \mu_\ell(\dd\bz)
   < \infty ,
 \end{align*}
 hence, by page 62 in Ikeda and Watanabe \cite{IkeWat} and \eqref{D3_bound}, we
 have
 \begin{equation}\label{J_3_bound-2}
  \begin{aligned}
   &\EE(\|\bJ_{t,t+T}^{(3)}\|^2) \\
   &= \ee^{-2s(\tbB)(t+T)}
      \sum_{\ell=1}^d
       \EE\biggl(\int_t^{t+T} \int_{\cU_d} \int_{\cU_1}
                  \bigl\|\ee^{(t+T-v)\tbB} \bz\bigr\|^2
                  \bbone_{\{\|\bz\|<\ee^{s(\tbB)v}\}}
                  \bbone_{\{w\leq X_{v,\ell}\}}
                  \, \dd v \, \mu_\ell(\dd\bz) \, \dd w\biggr) \\
   &\leq C_3^2 C_4
         \sum_{\ell=1}^d
          \int_t^\infty \int_{\cU_d}
           \ee^{-s(\tbB)v} \|\bz\|^2 \bbone_{\{\|\bz\|<\ee^{s(\tbB)v}\}}
           \, \dd v \, \mu_\ell(\dd\bz)
    = C_3^2 C_4 K_t^{(3)} .
  \end{aligned}
 \end{equation}
Thus, by \eqref{K_3_conv}, we obtain \ $\bJ_{t,t+T}^{(3)} \qmean \bzero$ \ as
 \ $t \to \infty$ \ for each \ $T \in \RR_+$.
\ Consequently, for each \ $T \in \RR_+$, \ we conclude
 \begin{equation}\label{J_3_conv_L1}
  \bJ_{t,t+T}^{(3)} \mean \bzero \qquad \text{as \ $t \to \infty$.}
 \end{equation}
Further, similarly as in case of \ $(D^{(4)}_t)_{t\in\RR_+}$, \ for each \ $t, T \in \RR_+$, \ we have
 \begin{align*}
  \bJ_{t,t+T}^{(4)}
  &= \ee^{-s(\tbB)(t+T)}
     \sum_{\ell=1}^d
      \int_t^{t+T} \int_{\cU_d} \int_{\cU_1}
       \ee^{(t+T-v)\tbB} \bz \bbone_{\{\|\bz\|\geq\ee^{s(\tbB)v}\}}
       \bbone_{\{w\leq X_{v-,\ell}\}}
       \, N_\ell(\dd v, \dd\bz, \dd w) \\
  &\quad
     - \ee^{-s(\tbB)(t+T)}
       \sum_{\ell=1}^d
        \int_t^{t+T} \int_{\cU_d} \int_{\cU_1}
         \ee^{(t+T-v)\tbB} \bz \bbone_{\{\|\bz\|\geq\ee^{s(\tbB)v}\}}
         \bbone_{\{w\leq X_{v,\ell}\}}
         \, \dd v \, \mu_\ell(\dd\bz) \, \dd w
 \end{align*}
 almost surely, thus for each \ $i \in \{1, \ldots, d\}$, \ by \eqref{C}, we
 obtain
 \begin{align*}
  |\be_i^\top \bJ_{t,t+T}^{(4)}|
  &\leq \ee^{-s(\tbB)(t+T)}
        \sum_{\ell=1}^d
         \int_t^{t+T} \int_{\cU_d} \int_{\cU_1}
          \be_i^\top \ee^{(t+T-v)\tbB} \bz
          \bbone_{\{\|\bz\|\geq\ee^{s(\tbB)v}\}} \bbone_{\{w\leq X_{v-,\ell}\}}
          \, N_\ell(\dd v, \dd\bz, \dd w) \\
  &\quad
        + \ee^{-s(\tbB)(t+T)}
          \sum_{\ell=1}^d
           \int_t^{t+T} \int_{\cU_d} \int_{\cU_1}
            \be_i^\top \ee^{(t+T-v)\tbB} \bz
            \bbone_{\{\|\bz\|\geq\ee^{s(\tbB)v}\}} \bbone_{\{w\leq X_{v,\ell}\}}
            \, \dd v \, \mu_\ell(\dd\bz) \, \dd w \\
  &\leq C_3 \sum_{\ell=1}^d
             \int_t^\infty \int_{\cU_d} \int_{\cU_1}
              \ee^{-s(\tbB)v} \|\bz\| \bbone_{\{\|\bz\|\geq\ee^{s(\tbB)v}\}}
              \bbone_{\{w\leq X_{v-,\ell}\}}
              \, N_\ell(\dd v, \dd\bz, \dd w) \\
  &\quad
        + C_3 \sum_{\ell=1}^d
               \int_t^\infty \int_{\cU_d} \int_{\cU_1}
                \ee^{-s(\tbB)v} \|\bz\| \bbone_{\{\|\bz\|\geq\ee^{s(\tbB)v}\}}
                \bbone_{\{w\leq X_{v,\ell}\}}
                \, \dd v \, \mu_\ell(\dd\bz) \, \dd w ,
 \end{align*}
 hence, by \eqref{D4_bound}, we get
 \begin{equation}\label{J_4_bound}
  |\be_i^\top \bJ_{t,t+T}^{(4)}| \leq C_3 K_t^{(4)}
 \end{equation}
 almost surely.
Consequently, by \eqref{K_4_conv}, for each \ $T \in \RR_+$, \ we conclude
 \begin{equation}\label{J_4_conv_L1}
  \bJ_{t,t+T}^{(4)} \mean \bzero \qquad \text{ as \ $t \to \infty$.}
 \end{equation}
In a similar way, for each \ $t, T \in \RR_+$, \ we have
 \[
   \bJ_{t,t+T}^{(5)}
   = \ee^{-s(\tbB)(t+T)}
     \int_t^{t+T} \int_{\cU_d} \ee^{(t+T-v)\tbB} \br \, M(\dd v, \dd\br)
     - \ee^{-s(\tbB)(t+T)}
       \int_t^{t+T} \int_{\cU_d} \ee^{(t+T-v)\tbB} \br \, \dd v \, \nu(\dd\br)
 \]
 almost surely, thus for each \ $i \in \{1, \ldots, d\}$, \ by \eqref{C}, we obtain
 \begin{align*}
  &|\be_i^\top \bJ_{t,t+T}^{(5)}| \\
  &\leq \ee^{-s(\tbB)(t+T)}
        \int_t^{t+T} \int_{\cU_d}
         \be_i^\top \ee^{(t+T-v)\tbB} \br \, M(\dd v, \dd\br)
        + \ee^{-s(\tbB)(t+T)}
          \int_t^{t+T} \int_{\cU_d}
           \be_i^\top \ee^{(t+T-v)\tbB} \br \, \dd v \, \nu(\dd\br) \\
  &\leq C_3 \int_t^\infty \int_{\cU_d}
             \ee^{-s(\tbB)v} \|\br\| \, M(\dd v, \dd\br)
        + C_3 \int_t^\infty \int_{\cU_d}
               \ee^{-s(\tbB)v} \|\br\| \, \dd v \, \nu(\dd\br) ,
 \end{align*}
 thus, by \eqref{D5_bound},
 \begin{equation}\label{J_5_bound}
  |\be_i^\top \bJ_{t,t+T}^{(5)}| \leq C_3 K_t^{(5)}
 \end{equation}
 almost surely.
Consequently, by \eqref{K_5_conv}, for each \ $T \in \RR_+$, \ we conclude
 \begin{equation}\label{J_5_conv_L1}
  \bJ_{t,t+T}^{(5)} \mean \bzero \qquad \text{as \ $t \to \infty$.}
 \end{equation}
The convergences \eqref{J_1_conv}, \eqref{J_2_conv_L1}, \eqref{J_3_conv_L1},
 \eqref{J_4_conv_L1} and \eqref{J_5_conv_L1} yield \eqref{Lemma_L1}.
In fact, it turned out that, for each \ $T \in \RR_+$, \ we have
 \ $\bJ_{t,t+T}^{(2)} \qmean \bzero$ \ and
 \ $\bJ_{t,t+T}^{(3)} \qmean \bzero$ \ as \ $t \to \infty$, \ but only
 \ $\bJ_{t,t+T}^{(4)} \mean \bzero$ \ and
 \ $\bJ_{t,t+T}^{(5)} \mean \bzero$ \ as \ $t \to \infty$.

Next we prove
 \begin{equation}\label{convCBIwL1}
  \ee^{-s(\tbB)t} \bX_t \mean w_{\bu,\bX_0} \tbu \qquad
  \text{as \ $t \to \infty$}
 \end{equation}
 by \eqref{Lemma_L1} and \eqref{convwv} with \ $\lambda = s(\tbB)$ \ and
 \ $\bv = \bu$.
\ For each \ $t, T \in \RR_+$ \ and \ $i \in \{1, \ldots, d\}$, \ using
 \eqref{Delta_t_t+T} and the identity \ $\bI_d = \sum_{j=1}^d \be_j \be_j^\top$,
 \ we have
 \begin{align*}
  &\ee^{-s(\tbB)(t+T)} \be_i^\top \bX_{t+T}
   = \be_i^\top \bDelta_{t,t+T}
     + \ee^{-s(\tbB)(t+T)} \be_i^\top \ee^{T\tbB} \bX_t \\
  &= \be_i^\top \bDelta_{t,t+T}
     + \sum_{j=1}^d
        (\be_i^\top \ee^{-s(\tbB)T} \ee^{T\tbB} \be_j)
        \ee^{-s(\tbB)t} \be_j^\top \bX_t \\
  &= \be_i^\top \bDelta_{t,t+T}
     + \sum_{j=1}^d
        (\be_i^\top \tbu \bu^\top \be_j) \ee^{-s(\tbB)t} \be_j^\top \bX_t
     + \sum_{j=1}^d
        \bigl[\be_i^\top (\ee^{-s(\tbB)T} \ee^{T\tbB} - \tbu \bu^\top) \be_j\bigr]
        \ee^{-s(\tbB)t} \be_j^\top \bX_t \\
  &= \be_i^\top \bDelta_{t,t+T}
     + (\be_i^\top \tbu) \ee^{-s(\tbB)t} \bu^\top \bX_t
     + \sum_{j=1}^d
        \bigl[\be_i^\top (\ee^{-s(\tbB)T} \ee^{T\tbB} - \tbu \bu^\top) \be_j\bigr]
        \ee^{-s(\tbB)t} \be_j^\top \bX_t .
 \end{align*}
By \eqref{C}, for each \ $T \in \RR_+$ \ and \ $i, j \in \{1, \ldots, d\}$, \ we
 have
 \[
   \bigl|\be_i^\top (\ee^{-s(\tbB)T} \ee^{T\tbB} - \tbu \bu^\top) \be_j\bigr|
   \leq \|\ee^{-s(\tbB)T} \ee^{T\tbB} - \tbu \bu^\top\|
   \leq C_1 \ee^{-C_2T}
   \leq C_5 \ee^{-C_2T} \be_i^\top \tbu \bu^\top \be_j ,
 \]
 where
 \[
   C_5 := \max_{i,j\in\{1,\ldots,d\}} \frac{C_1}{\be_i^\top \tbu \bu^\top \be_j}
        = \frac{C_1}{\min\limits_{i,j\in\{1,\ldots,d\}}\be_i^\top\tbu\bu^\top\be_j}
        = \frac{C_1}
               {\min\limits_{i\in\{1,\ldots,d\}}\be_i^\top\tbu
                \min\limits_{j\in\{1,\ldots,d\}}\bu^\top\be_j}
        \in \RR_{++} ,
 \]
 since \ $\bu, \tbu \in \RR^d_{++}$.
\ Hence for each \ $t, T \in \RR_+$ \ and \ $i \in \{1, \ldots, d\}$, \ we have
 \begin{equation}\label{difference}
  \begin{aligned}
   &|\ee^{-s(\tbB)(t+T)} \be_i^\top \bX_{t+T} - w_{\bu,\bX_0} \be_i^\top \tbu| \\
   &\leq |\be_i^\top \bDelta_{t,t+T}|
         + |(\be_i^\top \tbu)
            (\ee^{-s(\tbB)t} \bu^\top \bX_t - w_{\bu,\bX_0})| \\
   &\quad
         + \sum_{j=1}^d
            C_5 \ee^{-C_2T}
            (\be_i^\top \tbu \bu^\top \be_j) \ee^{-s(\tbB)t} \be_j^\top \bX_t \\
   &\leq |\be_i^\top \bDelta_{t,t+T}|
         + \|\tbu\|
           |\ee^{-s(\tbB)t} \langle\bu, \bX_t\rangle - w_{\bu,\bX_0}|
         + C_5 \|\tbu\| \ee^{-C_2T-s(\tbB)t} \langle\bu, \bX_t\rangle .
  \end{aligned}
 \end{equation}
For each \ $t, T \in \RR_+$ \ and \ $i \in \{1, \ldots, d\}$, \ by
 \eqref{difference}, we obtain
 \begin{align*}
  &\EE(|\ee^{-s(\tbB)(t+T)} \be_i^\top \bX_{t+T}
        - w_{\bu,\bX_0} \be_i^\top \tbu|) \\
  &\leq \EE(|\be_i^\top \bDelta_{t,t+T}|)
        + \|\tbu\|
          \EE(|\ee^{-s(\tbB)t} \langle\bu, \bX_t\rangle - w_{\bu,\bX_0}|)
        + C_5 \|\tbu\| \ee^{-C_2T-s(\tbB)t} \EE(\langle\bu, \bX_t\rangle) .
 \end{align*}
By \eqref{Lemma_L1} and \eqref{convwv} with \ $\lambda = s(\tbB)$ \ and
 \ $\bv = \bu$, \ for each \ $T \in \RR_+$ \ and \ $i \in \{1, \ldots, d\}$, \ we
 obtain
 \begin{align*}
  &\limsup_{t\to\infty}
    \EE(|\ee^{-s(\tbB)t} \be_i^\top \bX_t - w_{\bu,\bX_0} \be_i^\top \tbu|)
   = \limsup_{t\to\infty}
      \EE(|\ee^{-s(\tbB)(t+T)} \be_i^\top \bX_{t+T}
           - w_{\bu,\bX_0} \be_i^\top \tbu|) \\
  &\leq C_5 \|\tbu\| \ee^{-C_2T}
        \limsup_{t\to\infty} \EE(\ee^{-s(\tbB)t} \langle\bu, \bX_t\rangle)
   = C_5 \|\tbu\| \ee^{-C_2T} \EE(w_{\bu,\bX_0}) ,
 \end{align*}
 hence, by \ $T \to \infty$, \ we conclude \eqref{convCBIwL1}.

Next we show that for each \ $m \in \NN$ \ and
 \ $\delta \in \RR_{++}$, \ we have
 \begin{equation}\label{Lemma}
  \bDelta_{n\delta,(n+m)\delta} \as \bzero \qquad \text{as \ $n \to \infty$.}
 \end{equation}
For each \ $m \in \NN$, \ $\delta \in \RR_{++}$ \ and
 \ $\vare \in \RR_{++}$, \ by \eqref{J_2_bound}, we obtain
 \begin{equation}\label{J_2_sumbound}
  \begin{aligned}
   \sum_{n=1}^\infty \PP(\|\bJ_{n\delta,(n+m)\delta}^{(2)}\| > \vare)
   &\leq \frac{1}{\vare^2}
         \sum_{n=1}^\infty \EE(\|\bJ_{n\delta,(n+m)\delta}^{(2)}\|^2) \\
   &\leq \frac{2 C_3^2 C_4}{\vare^2s(\tbB)} \biggl(\sum_{\ell=1}^d c_\ell\biggr)
         \sum_{n=1}^\infty \ee^{-s(\tbB)n\delta}
     < \infty ,
  \end{aligned}
 \end{equation}
 hence, by the Borel--Cantelli lemma, for each \ $m \in \NN$ \ and
 \ $\delta \in \RR_{++}$, \ we conclude
 \begin{equation}\label{J_2_conv}
  \bJ_{n\delta,(n+m)\delta}^{(2)} \as \bzero \qquad \text{as \ $n \to \infty$.}
 \end{equation}
For each \ $m \in \NN$, \ $\delta \in \RR_{++}$ \ and
 \ $\vare \in \RR_{++}$, \ by \eqref{J_3_bound-2}, we obtain
 \begin{equation}\label{K_3_bound-1}
  \sum_{n=1}^\infty \PP(\|\bJ_{n\delta,(n+m)\delta}^{(3)}\| > \vare)
   \leq \frac{1}{\vare^2}
        \sum_{n=1}^\infty \EE(\|\bJ_{n\delta,(n+m)\delta}^{(3)}\|^2)
   \leq \frac{C_3^2C_4}{\vare^2}
        \sum_{n=1}^\infty K_{n\delta}^{(3)} .
 \end{equation}
We will show
 \begin{equation}\label{K_3_bound}
  \sum_{n=1}^\infty K_{n\delta}^{(3)} < \infty .
 \end{equation}
By Fubini's theorem, we have
 \begin{align*}
  \sum_{n=1}^\infty K_{n\delta}^{(3)}
  &= \sum_{n=1}^\infty
      \sum_{\ell=1}^d
       \int_{n\delta}^\infty
        \int_{\cU_d}
        \ee^{-s(\tbB)u} \|\bz\|^2
         \bbone_{\{\|\bz\|<\ee^{s(\tbB)u}\}}
         \, \dd u \, \mu_\ell(\dd\bz) \\
  &= \sum_{\ell=1}^d \sum_{n=1}^\infty \sum_{k=n}^\infty
      \int_{k\delta}^{(k+1)\delta} \int_{\cU_d}
       \ee^{-s(\tbB)u} \|\bz\|^2
       \bbone_{\{\|\bz\|<\ee^{s(\tbB)u}\}}
       \, \dd u \, \mu_\ell(\dd\bz)
 \end{align*}
 \begin{align*}
  &= \sum_{\ell=1}^d \sum_{k=1}^\infty
      k \int_{k\delta}^{(k+1)\delta} \int_{\cU_d}
         \ee^{-s(\tbB)u} \|\bz\|^2
         \bbone_{\{\|\bz\|<\ee^{s(\tbB)u}\}}
         \, \dd u \, \mu_\ell(\dd\bz) \\
  &\leq \frac{1}{\delta}
        \sum_{\ell=1}^d
         \int_0^\infty \int_{\cU_d}
          u \ee^{-s(\tbB)u} \|\bz\|^2
          \bbone_{\{\|\bz\|<\ee^{s(\tbB)u}\}}
          \, \dd u \, \mu_\ell(\dd\bz) .
 \end{align*}
Here, for each \ $\ell \in \{1, \ldots, d\}$, \ using Fubini's theorem, we have
 \begin{align*}
  \int_0^\infty \int_{\cU_d}
    u \ee^{-s(\tbB)u} \|\bz\|^2 \bbone_{\{\|\bz\|<1\}}
    \, \dd u \, \mu_\ell(\dd\bz)
  &= \int_0^\infty u \ee^{-s(\tbB)u} \, \dd u \int_{\cU_d}
      \|\bz\|^2 \bbone_{\{\|\bz\|<1\}} \, \mu_\ell(\dd\bz) \\
  &= \frac{1}{s(\tbB)^2}
     \int_{\cU_d} \|\bz\|^2 \bbone_{\{\|\bz\|<1\}} \, \mu_\ell(\dd\bz)
   < \infty
 \end{align*}
 by Definition \ref{Def_admissible}, and
 \begin{align*}
  &\int_0^\infty \int_{\cU_d}
    u \ee^{-s(\tbB)u} \|\bz\|^2
    \bbone_{\{1\leq\|\bz\|<\ee^{s(\tbB)u}\}}
    \, \dd u \, \mu_\ell(\dd\bz) \\
  &= \int_{\cU_d}
      \biggl(\int_{\frac{1}{s(\tbB)}\log(\|\bz\|)}^\infty
             u \ee^{-s(\tbB)u} \, \dd u\biggr)
      \|\bz\|^2 \bbone_{\{\|\bz\|\geq1\}}
      \, \mu_\ell(\dd\bz) \\
  &= \int_{\cU_d}
      \frac{\|\bz\|^{-1}}{s(\tbB)^2}
      (1 + \log(\|\bz\|))
      \|\bz\|^2 \bbone_{\{\|\bz\|\geq1\}}
      \, \mu_\ell(\dd\bz) \\
  &\leq \frac{1}{s(\tbB)^2}
        \int_{\cU_d} \|\bz\|
        \bbone_{\{\|\bz\|\geq1\}} \, \mu_\ell(\dd\bz)
        + \frac{1}{s(\tbB)^2}
          \int_{\cU_d}
           \|\bz\| \log(\|\bz\|)
           \bbone_{\{\|\bz\|\geq1\}}
           \, \mu_\ell(\dd\bz)
   < \infty
 \end{align*}
 by the moment condition \eqref{moment_condition_xlogx} with \ $\lambda = s(\tbB)$.
\ Thus, for each \ $\delta \in \RR_{++}$, \ we obtain \eqref{K_3_bound}.
Hence, for each \ $m \in \NN$ \ and \ $\delta \in \RR_{++}$, \ by
 \eqref{K_3_bound-1} and by the Borel--Cantelli lemma, we conclude
 \begin{equation}\label{J_3_conv}
  \bJ_{n\delta,(n+m)\delta}^{(3)} \as \bzero \qquad \text{as \ $n \to \infty$.}
 \end{equation}
By \eqref{J_4_bound}, for each \ $m, n \in \NN$, \ $\delta \in \RR_{++}$
 \ and \ $i \in \{1, \ldots, d\}$, \ we have
 \ $|\be_i^\top \bJ_{n\delta,(n+m)\delta}^{(4)}| \leq C_3 K_{n\delta}^{(4)}$.
\ For each \ $\delta \in \RR_{++}$, \ the function
 \ $\NN \ni n \mapsto K_{n\delta}^{(4)}$ \ is decreasing almost surely, hence, by
 \eqref{K_4_conv}, we obtain
 \begin{equation}\label{J_3_conv-1}
  K_{n\delta}^{(4)} \as 0 \qquad \text{as \ $n \to \infty$.}
 \end{equation}
Consequently, for each \ $m \in \NN$ \ and \ $\delta \in \RR_{++}$, \ we conclude
 \begin{equation}\label{J_4_conv}
  \bJ_{n\delta,(n+m)\delta}^{(4)} \as \bzero \qquad \text{as \ $n \to \infty$.}
 \end{equation}
By \eqref{J_5_bound}, for each \ $m, n \in \NN$, \ $\delta \in \RR_{++}$
 \ and \ $i \in \{1, \ldots, d\}$, \ we have
 \ $|\be_i^\top \bJ_{n\delta,(n+m)\delta}^{(5)}| \leq C_3 K_{n\delta}^{(5)}$.
\ For each \ $\delta \in \RR_{++}$, \ the function
 \ $\NN \ni n \mapsto K_{n\delta}^{(5)}$ \ is decreasing almost surely, hence, by
 \eqref{K_5_conv}, we obtain
 \begin{equation}\label{J_5_conv-1}
  K_{n\delta}^{(5)} \as 0 \qquad \text{as \ $n \to \infty$.}
 \end{equation}
Consequently, for each \ $m \in \NN$ \ and \ $\delta \in \RR_{++}$, \ we conclude
 \begin{equation}\label{J_5_conv}
  \bJ_{n\delta,(n+m)\delta}^{(5)} \as \bzero \qquad \text{as \ $n \to \infty$.}
 \end{equation}
The representation \eqref{reprJ} and the convergences \eqref{J_1_conv},
 \eqref{J_2_conv}, \eqref{J_3_conv}, \eqref{J_4_conv} and \eqref{J_5_conv} yield
 \eqref{Lemma}.

Next, using the almost sure convergences \eqref{Lemma} and \eqref{convwu}, we will
 show the almost sure convergence of \ $(\ee^{-s(\tbB)t} \bX_t)_{t\in\RR_+}$
 \ along lattice times, i.e., we will prove that for each \ $\delta \in \RR_{++}$,
 \ we have
 \begin{equation}\label{lattice}
  \ee^{-s(\tbB)n\delta} \bX_{n\delta} \as w_{\bu,\bX_0} \tbu \qquad
  \text{as \ $n \to \infty$.}
 \end{equation}
By \eqref{difference}, \eqref{Lemma} and \eqref{convwu}, for each \ $m \in \NN$,
 \ $\delta \in \RR_{++}$ \ and \ $i \in \{1, \ldots, d\}$, \ we obtain
 \begin{align*}
  &\limsup_{n\to\infty}
    |\ee^{-s(\tbB)n\delta} \be_i^\top \bX_{n\delta}
     - w_{\bu,\bX_0} \be_i^\top \tbu|
   = \limsup_{n\to\infty}
      |\ee^{-s(\tbB)(n+m)\delta} \be_i^\top \bX_{(n+m)\delta}
       - w_{\bu,\bX_0} \be_i^\top \tbu| \\
  &\leq C_5 \|\tbu\| \ee^{-C_2m\delta}
        \limsup_{n\to\infty}
         \ee^{-s(\tbB)n\delta} \langle\bu, \bX_{n\delta}\rangle
   = C_5 \|\tbu\| \ee^{-C_2m\delta} w_{\bu,\bX_0}
 \end{align*}
 almost surely, hence, by \ $m \to \infty$, \ we conclude \eqref{lattice}.

The aim of the following discussion is to derive
 \begin{equation}\label{convCBIwas}
  \ee^{-s(\tbB)t} \bX_t \as w_{\bu,\bX_0} \tbu \qquad
  \text{as \ $t \to \infty$}
 \end{equation}
 by the help of the almost sure convergence \eqref{lattice} of
 \ $(\ee^{-s(\tbB)n\delta} \bX_{n\delta})_{n\in\NN}$ \ for all
 \ $\delta \in \RR_{++}$ \ together with the almost sure convergence \eqref{convwu}
 of \ $(\ee^{-s(\tbB)t} \langle\bu, \bX_t\rangle)_{t\in\RR_+}$.
\ For each \ $t \in \RR_+$, \ $i \in \{1, \ldots, d\}$, \ $n \in \NN$ \ and
 \ $\delta \in \RR_{++}$, \ we have
 \begin{align*}
  |\ee^{-s(\tbB)t} \be_i^\top \bX_t - w_{\bu,\bX_0} \be_i^\top \tbu|
  &\leq |\ee^{-s(\tbB)t} \be_i^\top \bX_t
         - \ee^{-s(\tbB)t} \be_i^\top \ee^{((n+1)\delta-t)\tbB} \bX_t| \\
  &\quad
        + |\ee^{-s(\tbB)t} \be_i^\top \ee^{((n+1)\delta-t)\tbB} \bX_t
           - w_{\bu,\bX_0} \be_i^\top \tbu| ,
 \end{align*}
 hence for each \ $i \in \{1, \ldots, d\}$, \ we get
 \begin{align*}
  &\limsup_{t\to\infty}
    |\ee^{-s(\tbB)t} \be_i^\top \bX_t - w_{\bu,\bX_0} \be_i^\top \tbu| \\
  &\leq \limsup_{\delta\downarrow0} \limsup_{n\to\infty}
         \sup_{t\in[n\delta,(n+1)\delta)}
          |\ee^{-s(\tbB)t} \be_i^\top \bX_t - w_{\bu,\bX_0} \be_i^\top \tbu|
   \leq I_i^{(1)} + I_i^{(2)}
 \end{align*}
 with
 \begin{align*}
  I_i^{(1)}
  &:= \limsup_{\delta\downarrow0} \limsup_{n\to\infty}
       \sup_{t\in[n\delta,(n+1)\delta)}
        |\ee^{-s(\tbB)t} \be_i^\top \bX_t
         - \ee^{-s(\tbB)t} \be_i^\top \ee^{((n+1)\delta-t)\tbB} \bX_t| \\
  &\phantom{:}
    = \limsup_{\delta\downarrow0} \limsup_{n\to\infty}
       \sup_{t\in[n\delta,(n+1)\delta)}
        |\ee^{-s(\tbB)t} \be_i^\top (\bI_d - \ee^{((n+1)\delta-t)\tbB}) \bX_t| , \\
  I_i^{(2)}
  &:= \limsup_{\delta\downarrow0} \limsup_{n\to\infty}
       \sup_{t\in[n\delta,(n+1)\delta)}
        |\ee^{-s(\tbB)t} \be_i^\top \ee^{((n+1)\delta-t)\tbB} \bX_t
         - w_{\bu,\bX_0} \be_i^\top \tbu| .
 \end{align*}
For each \ $i \in \{1, \ldots, d\}$, \ we have
 \ $I_i^{(1)} \leq I^{(1,1)} I^{(1,2)}$ \ with
 \begin{align*}
  I^{(1,1)}
  &:= \limsup_{\delta\downarrow0} \limsup_{n\to\infty}
       \sup_{t\in[n\delta,(n+1)\delta)}
        \|\bI_d - \ee^{((n+1)\delta-t)\tbB}\| \\
  &\phantom{:}
    = \limsup_{\delta\downarrow0} \limsup_{n\to\infty}
       \sup_{v\in(0,\delta]} \|\bI_d - \ee^{v\tbB}\|
    = \limsup_{\delta\downarrow0} \sup_{v\in(0,\delta]}
       \|\bI_d - \ee^{v\tbB}\|
    = 0 , \\
  I^{(1,2)}
  &:= \limsup_{\delta\downarrow0} \limsup_{n\to\infty}
       \sup_{t\in[n\delta,(n+1)\delta)}
        \ee^{-s(\tbB)t} \|\bX_t\| .
 \end{align*}
For each \ $\bx \in \RR^d_+$, \ we have
 \begin{equation}\label{x}
  \|\bx\|^2
  = \sum_{j=1}^d (\be_j^\top \bx)^2
  = \sum_{j=1}^d \frac{[(\be_j^\top \bu)(\be_j^\top \bx)]^2}{(\be_j^\top \bu)^2}
  \leq \sum_{j=1}^d \frac{\langle\bu, \bx\rangle^2}{(\be_j^\top \bu)^2}
  = C_6 \langle\bu, \bx\rangle^2
 \end{equation}
 with \ $C_6 := \sum_{j=1}^d \frac{1}{(\be_j^\top \bu)^2} \in \RR_{++}$, \ since
 \ $\bu \in \RR^d_{++}$.
\ Thus, by the almost sure convergence \eqref{convwu}, we obtain
 \[
   I^{(1,2)}
   \leq \sqrt{C_6}
        \limsup_{\delta\downarrow0} \limsup_{n\to\infty}
        \sup_{t\in[n\delta,(n+1)\delta)}
         \ee^{-s(\tbB)t} \langle\bu, \bX_t\rangle
   = \sqrt{C_6} w_{\bu,\bX_0}
   < \infty
 \]
 almost surely, and hence, we conclude \ $I_i^{(1)} = 0$ \ almost surely for each
 \ $i \in \{1, \ldots, d\}$.

In order to show \ $I_i^{(2)} = 0$ \ almost surely for each
 \ $i \in \{1, \ldots, d\}$, \ by Lemma \ref{SDE_transform_sol}, we consider the
 representation
 \[
   \ee^{-s(\tbB)t} \be_i^\top \ee^{((n+1)\delta-t)\tbB} \bX_t
   - w_{\bu,\bX_0} \be_i^\top \tbu
   = \sum_{j=0}^5 J_{t,i}^{(n,\delta,j)}
 \]
 with
 \begin{align*}
  J_{t,i}^{(n,\delta,0)}
  &:= \ee^{-s(\tbB)t} \be_i^\top \ee^{\delta\tbB} \bX_{n\delta}
      - w_{\bu,\bX_0} \be_i^\top \tbu , \\
  J_{t,i}^{(n,\delta,1)}
  &:= \ee^{-s(\tbB)t}
      \int_{n\delta}^t \be_i^\top \ee^{((n+1)\delta-v)\tbB} \tBbeta \, \dd v , \\
  J_{t,i}^{(n,\delta,2)}
  &:= \ee^{-s(\tbB)t}
      \sum_{\ell=1}^d
       \int_{n\delta}^t
        \be_i^\top \ee^{((n+1)\delta-v)\tbB} \be_\ell \sqrt{2 c_\ell X_{v,\ell}}
        \, \dd W_{v,\ell} , \\
  J_{t,i}^{(n,\delta,3)}
  &:= \ee^{-s(\tbB)t}
      \sum_{\ell=1}^d
       \int_{n\delta}^t \int_{\cU_d} \int_{\cU_1}
        \be_i^\top \ee^{((n+1)\delta-v)\tbB} \bz
        \bbone_{\{\|\bz\|<\ee^{s(\tbB)v}\}}
        \bbone_{\{w\leq X_{v-,\ell}\}}
        \, \tN_\ell(\dd v, \dd\bz, \dd w) , \\
  J_{t,i}^{(n,\delta,4)}
  &:= \ee^{-s(\tbB)t}
      \sum_{\ell=1}^d
       \int_{n\delta}^t \int_{\cU_d} \int_{\cU_1}
        \be_i^\top \ee^{((n+1)\delta-v)\tbB}
        \bz \bbone_{\{\|\bz\|\geq\ee^{s(\tbB)v}\}}
        \bbone_{\{w\leq X_{v-,\ell}\}}
        \, \tN_\ell(\dd v, \dd\bz, \dd w) , \\
  J_{t,i}^{(n,\delta,5)}
  &:= \ee^{-s(\tbB)t}
      \int_{n\delta}^t \int_{\cU_d}
       \be_i^\top \ee^{((n+1)\delta-v)\tbB} \br \, \tM(\dd v, \dd\br)
 \end{align*}
 for all \ $i \in \{1, \ldots, d\}$, \ $n \in \NN$, \ $\delta \in \RR_{++}$ \ and
 \ $t \in [n\delta, (n+1)\delta)$.
\ For each \ $i \in \{1, \ldots, d\}$, \ $n \in \NN$, \ $\delta \in \RR_{++}$ \ and
 \ $t \in [n\delta, (n+1)\delta)$, \ we have
 \[
   J_{t,i}^{(n,\delta,0)}
   \leq J_{t,i}^{(n,\delta,0,1)} + J_{t,i}^{(n,\delta,0,2)}
        + J_i^{(n,\delta,0,3)}
 \]
 with
 \begin{align*}
  J_{t,i}^{(n,\delta,0,1)}
  &:= \ee^{-s(\tbB)t} \be_i^\top \ee^{\delta\tbB} \bX_{n\delta}
      - \ee^{-s(\tbB)t} \be_i^\top \bX_{n\delta}
    = \ee^{-s(\tbB)t} \be_i^\top (\ee^{\delta\tbB} - \bI_d) \bX_{n\delta} , \\
  J_{t,i}^{(n,\delta,0,2)}
  &:= \ee^{-s(\tbB)t} \be_i^\top \bX_{n\delta}
      - \ee^{-s(\tbB)n\delta} \be_i^\top \bX_{n\delta}
    = (\ee^{-s(\tbB)(t-n\delta)} - 1) \ee^{-s(\tbB)n\delta}
      \be_i^\top \bX_{n\delta} , \\
  J_i^{(n,\delta,0,3)}
  &:= \ee^{-s(\tbB)n\delta} \be_i^\top \bX_{n\delta}
      - w_{\bu,\bX_0} \be_i^\top \tbu .
 \end{align*}
For each \ $i \in \{1, \ldots, d\}$, \ we have
 \begin{align*}
  &\limsup_{\delta\downarrow0} \limsup_{n\to\infty}
    \sup_{t\in[n\delta,(n+1)\delta)}
     |J_{t,i}^{(n,\delta,0,1)}| \\
  &\leq \biggl(\limsup_{\delta\downarrow0} \|\ee^{\delta\tbB} - \bI_d\|\biggr)
        \biggl(\limsup_{\delta\downarrow0} \limsup_{n\to\infty}
                \sup_{t\in[n\delta,(n+1)\delta)}
                 \ee^{-s(\tbB)t} \|\bX_{n\delta}\|\biggr) ,
 \end{align*}
 where \ $\limsup_{\delta\downarrow0} \|\ee^{\delta\tbB} - \bI_d\| = 0$, \ and, by
 \eqref{x} and \eqref{convwu},
 \begin{align*}
  \limsup_{\delta\downarrow0} \limsup_{n\to\infty}
   \sup_{t\in[n\delta,(n+1)\delta)}
    \ee^{-s(\tbB)t} \|\bX_{n\delta}\|
  &\leq \sqrt{C_6} \limsup_{\delta\downarrow0} \limsup_{n\to\infty}
        \ee^{-s(\tbB)n\delta} \langle\bu, \bX_{n\delta}\rangle \\
  &= \sqrt{C_6} w_{\bu,\bX_0}
   < \infty
 \end{align*}
 almost surely, hence we get
 \ $\limsup_{\delta\downarrow0} \limsup_{n\to\infty}
     \sup_{t\in[n\delta,(n+1)\delta)}
      |J_{t,i}^{(n,\delta,0,1)}|
    = 0$
 \ almost surely.
Moreover, for each \ $i \in \{1, \ldots, d\}$, \ we have
 \begin{align*}
  &\limsup_{\delta\downarrow0} \limsup_{n\to\infty}
    \sup_{t\in[n\delta,(n+1)\delta)}
     |J_{t,i}^{(n,\delta,0,2)}| \\
  &\leq \biggl(\limsup_{\delta\downarrow0} \limsup_{n\to\infty}
                \sup_{t\in[n\delta,(n+1)\delta)}
                 |\ee^{-s(\tbB)(t-n\delta)} - 1|\biggr)
        \biggl(\limsup_{\delta\downarrow0} \limsup_{n\to\infty}
                  \ee^{-s(\tbB)n\delta} \|\bX_{n\delta}\|\biggr) ,
 \end{align*}
 where
 \[
   \limsup_{\delta\downarrow0} \limsup_{n\to\infty}
    \sup_{t\in[n\delta,(n+1)\delta)}
     |\ee^{-s(\tbB)(t-n\delta)} - 1|
   = \limsup_{\delta\downarrow0} |\ee^{-s(\tbB)\delta} - 1|
   = 0 ,
 \]
 and, by \eqref{x} and \eqref{convwu},
 \[
   \limsup_{\delta\downarrow0} \limsup_{n\to\infty}
     \ee^{-s(\tbB)n\delta} \|\bX_{n\delta}\|
   \leq \sqrt{C_6} \limsup_{\delta\downarrow0} \limsup_{n\to\infty}
        \ee^{-s(\tbB)n\delta} \langle\bu, \bX_{n\delta}\rangle
   = \sqrt{C_6} w_{\bu,\bX_0}
   < \infty
 \]
 almost surely, hence we obtain
 \ $\limsup_{\delta\downarrow0} \limsup_{n\to\infty}
     \sup_{t\in[n\delta,(n+1)\delta)}
      |J_{t,i}^{(n,\delta,0,2)}|
    = 0$
 \ almost surely.
Further, for each \ $i \in \{1, \ldots, d\}$, \ by \eqref{lattice}, we get
 \[
   \limsup_{\delta\downarrow0} \limsup_{n\to\infty}
    \sup_{t\in[n\delta,(n+1)\delta)}
     |J_i^{(n,\delta,0,3)}|
   = \limsup_{\delta\downarrow0} \limsup_{n\to\infty} |J_i^{(n,\delta,0,3)}|
   = 0 \qquad \text{almost surely,}
 \]
 and we conclude
 \begin{equation}\label{J0}
  \limsup_{\delta\downarrow0} \limsup_{n\to\infty}
   \sup_{t\in[n\delta,(n+1)\delta)}
    |J_{t,i}^{(n,\delta,0)}|
  = 0 \qquad \text{almost surely.}
 \end{equation}
Moreover, for all \ $i \in \{1, \ldots, d\}$, \ $n \in \NN$,
 \ $\delta \in \RR_{++}$ \ and \ $t \in [n\delta, (n+1)\delta)$, \ by \eqref{C},
 \begin{align*}
  |J_{t,i}^{(n,\delta,1)}|
  &\leq \ee^{-s(\tbB)t}
         \int_{n\delta}^t |\be_i^\top \ee^{((n+1)\delta-v)\tbB} \tBbeta| \, \dd v
   \leq \|\tBbeta\| \ee^{-s(\tbB)t}
        \int_{n\delta}^t \|\ee^{((n+1)\delta-v)\tbB}\| \, \dd v       \\
  &\leq C_3 \|\tBbeta\| \ee^{-s(\tbB)t}
        \int_{n\delta}^t \ee^{s(\tbB)((n+1)\delta-v)} \, \dd v
   \leq C_3 \|\tBbeta\| \ee^{s(\tbB)((n+1)\delta-t)}
        \int_{n\delta}^\infty \ee^{-s(\tbB)v} \, \dd v \\
  &= \frac{C_3\|\tBbeta\|}{s(\tbB)} \ee^{s(\tbB)((n+1)\delta-t)}
        \ee^{-s(\tbB)n\delta}
   = \frac{C_3\|\tBbeta\|}{s(\tbB)} \ee^{s(\tbB)(\delta-t)} ,
 \end{align*}
 thus
 \begin{equation}\label{J1}
  \begin{aligned}
   \limsup_{\delta\downarrow0} \limsup_{n\to\infty}
    \sup_{t\in[n\delta,(n+1)\delta)}
     |J_{t,i}^{(n,\delta,1)}|
   &\leq \frac{C_3\|\tBbeta\|}{s(\tbB)}
         \limsup_{\delta\downarrow0} \limsup_{n\to\infty}
          \sup_{t\in[n\delta,(n+1)\delta)}
           \ee^{s(\tbB)(\delta-t)} \\
   &= \frac{C_3\|\tBbeta\|}{s(\tbB)}
      \limsup_{\delta\downarrow0} \limsup_{n\to\infty} \ee^{-s(\tbB)(n-1)\delta}
    = 0 .
  \end{aligned}
 \end{equation}
Further, for each \ $i \in \{1, \ldots, d\}$, \ $n \in \NN$,
 \ $\delta \in \RR_{++}$ \ and \ $t \in [n\delta, (n+1)\delta)$, \ we have
 \ $|J_{t,i}^{(n,\delta,2)}| \leq |M_{t,i}^{(n,\delta,2)}|$, \ where
 \[
   M_{t,i}^{(n,\delta,2)}
   := \ee^{-s(\tbB)n\delta}
      \sum_{\ell=1}^d
       \int_{n\delta}^t
        \be_i^\top \ee^{((n+1)\delta-v)\tbB} \be_\ell \sqrt{2 c_\ell X_{v,\ell}}
        \, \dd W_{v,\ell} ,
   \qquad t \in \RR_+ ,
 \]
 is a square integrable martingale (which can be checked as in case of
 \ $(Z_t^{(2)})_{t\in\RR_+}$ \ in the proof of Theorem \ref{convCBIasL1}).
By the maximal inequality for submartingales and by \eqref{J_2_bound}, for each
 \ $i \in \{1, \ldots, d\}$, \ $n \in \NN$, \ $\delta \in \RR_{++}$ \ and
 \ $\vare \in \RR_{++}$, \ we obtain
 \begin{align*}
  &\PP\biggl(\sup_{t\in[n\delta,(n+1)\delta)} |J_{t,i}^{(n,\delta,2)}|
             > \vare\biggr)
   \leq \PP\biggl(\sup_{t\in[n\delta,(n+1)\delta)} |M_{t,i}^{(n,\delta,2)}|
                  > \vare\biggr)
   \leq \frac{1}{\vare^2} \EE\bigl[(M_{(n+1)\delta,i}^{(n,\delta,2)})^2\bigr] \\
  &= \frac{1}{\vare^2} \ee^{-2s(\tbB)n\delta}
     \sum_{\ell=1}^d
      \EE\biggl(\int_{n\delta}^{(n+1)\delta}
                 |\be_i^\top \ee^{((n+1)\delta-v)\tbB} \be_\ell|^2
                 2 c_\ell X_{v,\ell}
                 \, \dd v\biggr) \\
  &= \frac{2}{\vare^2} \ee^{-2s(\tbB)n\delta}
     \sum_{\ell=1}^d
      c_\ell \int_{n\delta}^{(n+1)\delta}
              |\be_i^\top \ee^{((n+1)\delta-v)\tbB} \be_\ell|^2 \EE(X_{v,\ell})
              \, \dd v \\
  &\leq \frac{2}{\vare^2} \ee^{-2s(\tbB)n\delta}
        \sum_{\ell=1}^d
         c_\ell \int_{n\delta}^{(n+1)\delta}
                 \|\ee^{((n+1)\delta-v)\tbB} \be_\ell\|^2 \EE(X_{v,\ell}) \, \dd v
   = \frac{1}{\vare^2} \ee^{2s(\tbB)\delta}
     \EE(\|\bJ_{n\delta,(n+1)\delta}^{(2)}\|^2) ,
 \end{align*}
 hence, for each \ $i \in \{1, \ldots, d\}$,
 \ $\delta \in \RR_{++}$ \ and \ $\vare \in \RR_{++}$, \ we obtain
 \[
   \sum_{n=1}^\infty
    \PP\biggl(\sup_{t\in[n\delta,(n+1)\delta)} |J_{t,i}^{(n,\delta,2)}|
              > \vare\biggr)
   \leq \frac{1}{\vare^2} \ee^{2s(\tbB)\delta}
        \sum_{n=1}^\infty \EE(\|\bJ_{n\delta,(n+1)\delta}^{(2)}\|^2)
   < \infty ,
 \]
 by \eqref{J_2_sumbound}.
By the Borel--Cantelli lemma, for each
 \ $i \in \{1, \ldots, d\}$ \ and \ $\delta \in \RR_{++}$, \ we conclude
 \[
   \sup_{t\in[n\delta,(n+1)\delta)} |J_{t,i}^{(n,\delta,2)}| \as 0 \qquad
   \text{as \ $n \to \infty$,}
 \]
 and hence
 \begin{equation}\label{J2}
  \limsup_{\delta\downarrow0} \limsup_{n\to\infty} \sup_{t\in[n\delta,(n+1)\delta)}
   |J_{t,i}^{(n,\delta,2)}|
  = 0 \qquad \text{almost surely.}
 \end{equation}
In a similar way, for each \ $i \in \{1, \ldots, d\}$, \ $n \in \NN$,
 \ $\delta \in \RR_{++}$ \ and \ $t \in [n\delta, (n+1)\delta)$, \ we have
 \ $|J_{t,i}^{(n,\delta,3)}| \leq |M_{t,i}^{(n,\delta,3)}|$, \ where
 \[
   M_{t,i}^{(n,\delta,3)}
   := \ee^{-s(\tbB)n\delta}
      \sum_{\ell=1}^d
       \int_{n\delta}^t \int_{\cU_d} \int_{\cU_1}
        \be_i^\top \ee^{((n+1)\delta-v)\tbB} \bz
        \bbone_{\{\|\bz\|<\ee^{s(\tbB)v}\}}
        \bbone_{\{w\leq X_{v-,\ell}\}}
        \, \tN_\ell(\dd v, \dd\bz, \dd w)
 \]
 for \ $t \in \RR_+$ \ defines a square integrable martingale (which can be checked
 as in case of \ $(Z_t^{(3)})_{t\in\RR_+}$ \ in the proof of Theorem
 \ref{convCBIasL1}).
By the maximal inequality for submartingales and by \eqref{J_3_bound-2},
 \eqref{K_3_bound-1} and \eqref{K_3_bound}, for each
 \ $i \in \{1, \ldots, d\}$, \ $\delta \in \RR_{++}$ \ and
 \ $\vare \in \RR_{++}$, \ we obtain
 \begin{align*}
  \sum_{n=1}^\infty
   \PP\biggl(\sup_{t\in[n\delta,(n+1)\delta)} |J_{t,i}^{(n,\delta,3)}|
             > \vare\biggr)
  &\leq \sum_{n=1}^\infty
         \PP\biggl(\sup_{t\in[n\delta,(n+1)\delta)} |M_{t,i}^{(n,\delta,3)}|
                   > \vare\biggr)
   \leq \frac{1}{\vare^2}
        \sum_{n=1}^\infty \EE\bigl[(M_{(n+1)\delta,i}^{(n,\delta,3)})^2\bigr] \\
  &\leq \frac{1}{\vare^2} \ee^{2s(\tbB)\delta}
        \sum_{n=1}^\infty \EE(\|\bJ_{n\delta,(n+1)\delta}^{(3)}\|^2)
   \leq \frac{C_3^2 C_4}{\vare^2} \ee^{2s(\tbB)\delta}
        \sum_{n=1}^\infty K_{n\delta}^{(3)}
   < \infty ,
 \end{align*}
 thus, by the Borel--Cantelli lemma, for each
 \ $i \in \{1, \ldots, d\}$ \ and \ $\delta \in \RR_{++}$, \ we conclude
 \[
   \sup_{t\in[n\delta,(n+1)\delta)} |J_{t,i}^{(n,\delta,3)}| \as 0 \qquad
   \text{as \ $n \to \infty$,}
 \]
 and hence
 \begin{equation}\label{J3}
  \limsup_{\delta\downarrow0} \limsup_{n\to\infty} \sup_{t\in[n\delta,(n+1)\delta)}
   |J_{t,i}^{(n,\delta,3)}|
  = 0 \qquad \text{almost surely.}
 \end{equation}
Next, for each \ $i \in \{1, \ldots, d\}$, \ $n \in \NN$, \ $\delta \in \RR_{++}$
 \ and \ $t \in [n\delta, (n+1)\delta)$, \ as in case of
 \ $\bJ_{t,t+T}^{(4)}$, \ we obtain
 \begin{align*}
  |J_{t,i}^{(n,\delta,4)}|
  &\leq C_3 \ee^{-s(\tbB)t}
        \sum_{\ell=1}^d
             \int_{n\delta}^\infty \int_{\cU_d} \int_{\cU_1}
              \ee^{s(\tbB)((n+1)\delta-v)} \|\bz\|
              \bbone_{\{\|\bz\|\geq\ee^{s(\tbB)v}\}}
              \bbone_{\{w\leq X_{v-,\ell}\}}
              \, N_\ell(\dd v, \dd\bz, \dd w) \\
  &\quad
        + C_3 \ee^{-s(\tbB)t}
          \sum_{\ell=1}^d
               \int_{n\delta}^\infty \int_{\cU_d} \int_{\cU_1}
                \ee^{s(\tbB)((n+1)\delta-v)} \|\bz\|
                \bbone_{\{\|\bz\|\geq\ee^{s(\tbB)v}\}}
                \bbone_{\{w\leq X_{v,\ell}\}}
                \, \dd v \, \mu_\ell(\dd\bz) \, \dd w
 \end{align*}
 almost surely, hence, by \eqref{D4_bound},
 \begin{align*}
  \sup_{t\in[n\delta,(n+1)\delta)} |J_{t,i}^{(n,\delta,4)}|
   \leq C_3 \ee^{s(\tbB)\delta}
        \sum_{\ell=1}^d
             \int_{n\delta}^\infty \int_{\cU_d} \int_{\cU_1}
              \ee^{-s(\tbB)v} \|\bz\| \bbone_{\{\|\bz\|\geq\ee^{s(\tbB)v}\}}
              \bbone_{\{w\leq X_{v-,\ell}\}}
              \, N_\ell(\dd v, \dd\bz, \dd w)
 \end{align*}
 \begin{align*}
  + C_3 \ee^{s(\tbB)\delta}
          \sum_{\ell=1}^d
               \int_{n\delta}^\infty \int_{\cU_d} \int_{\cU_1}
                \ee^{-s(\tbB)v} \|\bz\| \bbone_{\{\|\bz\|\geq\ee^{s(\tbB)v}\}}
                \bbone_{\{w\leq X_{v,\ell}\}}
                \, \dd v \, \mu_\ell(\dd\bz) \, \dd w
   = C_3 \ee^{s(\tbB)\delta} K_{n\delta}^{(4)}
 \end{align*}
 almost surely.
By \eqref{J_3_conv-1}, for each \ $i \in \{1, \ldots, d\}$, \ we conclude
 \begin{equation}\label{J4}
  \limsup_{\delta\downarrow0} \limsup_{n\to\infty} \sup_{t\in[n\delta,(n+1)\delta)}
   |J_{t,i}^{(n,\delta,4)}|
  \leq \limsup_{\delta\downarrow0} \limsup_{n\to\infty}
        C_3 \ee^{s(\tbB)\delta} K_{n\delta}^{(4)}
  = 0 \qquad \text{almost surely.}
 \end{equation}
Finally, for each \ $i \in \{1, \ldots, d\}$, \ $n \in \NN$,
 \ $\delta \in \RR_{++}$ \ and \ $t \in [n\delta, (n+1)\delta)$, \ as in case of
 \ $\bJ_{t,t+T}^{(5)}$, \ we obtain
 \begin{align*}
  |J_{t,i}^{(n,\delta,5)}|
  &\leq C_3 \ee^{-s(\tbB)t}
        \int_{n\delta}^\infty \int_{\cU_d}
         \ee^{s(\tbB)((n+1)\delta-v)} \|\br\| \, M(\dd v, \dd\br) \\
  &\quad
        + C_3 \ee^{-s(\tbB)t}
          \int_{n\delta}^\infty \int_{\cU_d}
           \ee^{s(\tbB)((n+1)\delta-v)} \|\br\| \, \dd v \, \nu(\dd\br)
 \end{align*}
 almost surely, hence, by \eqref{D5_bound},
 \begin{align*}
  \sup_{t\in[n\delta,(n+1)\delta)} |J_{t,i}^{(n,\delta,5)}|
  &\leq C_3 \ee^{s(\tbB)\delta}
        \int_{n\delta}^\infty \int_{\cU_d}
         \ee^{-s(\tbB)v} \|\br\| \, M(\dd v, \dd\br) \\
  &\quad
        + C_3 \ee^{s(\tbB)\delta}
          \int_{n\delta}^\infty \int_{\cU_d}
           \ee^{-s(\tbB)v} \|\br\| \, \dd v \, \nu(\dd\bz)
   = C_3 \ee^{s(\tbB)\delta} K_{n\delta}^{(5)}
 \end{align*}
 almost surely.
By \eqref{J_5_conv-1}, for each \ $i \in \{1, \ldots, d\}$, \ we conclude
 \begin{equation}\label{J5}
  \limsup_{\delta\downarrow0} \limsup_{n\to\infty} \sup_{t\in[n\delta,(n+1)\delta)}
   |J_{t,i}^{(n,\delta,5)}|
  \leq \limsup_{\delta\downarrow0} \limsup_{n\to\infty}
        C_3 \ee^{s(\tbB)\delta} K_{n\delta}^{(5)}
  = 0 \qquad \text{almost surely.}
 \end{equation}
The convergences \eqref{J0}, \eqref{J1}, \eqref{J2}, \eqref{J3}, \eqref{J4} and
 \eqref{J5} yield \eqref{convCBIwas}.

Finally, we check that if the moment condition \eqref{moment_condition_xlogx} does not hold for \ $\lambda = s(\tbB)$,
 \ then \ $\ee^{-s(\tbB)t}\bX_t$ \ does not converge in \ $L_1$ \ as \ $t\to\infty$,
 \ provided that \ $\PP(\bX_0 = \bzero)<1$ \ or \ $\tBbeta \ne \bzero$. \
On the contrary, let us suppose that \ $\ee^{-s(\tbB)t}\bX_t$ \ converges in \ $L_1$ \ as \ $t\to\infty$.
\ Then, especially, \ $\ee^{-s(\tbB)t}\langle \bu, \bX_t \rangle $ \ converges in \ $L_1$ \ as \ $t\to\infty$,
 \ which leads us to a contradiction by Theorem \ref{convCBIasL1}.
\proofend

Next we present \ $L_2$-convergence results for supercritical and irreducible
 multi-type CBI processes.

\begin{Thm}\label{convCBIL2}
Let \ $(\bX_t)_{t\in\RR_+}$ \ be a supercritical and irreducible multi-type CBI
 process with parameters \ $(d, \bc, \Bbeta, \bB, \nu, \bmu)$ \ such that
 \ $\EE(\|\bX_0\|^2) < \infty$ \ and the moment conditions
 \begin{equation}\label{moment_condition_CBI2}
  \sum_{\ell=1}^d
   \int_{\cU_d} \|\bz\|^2 \bbone_{\{\|\bz\|\geq1\}} \, \mu_\ell(\dd\bz) < \infty ,
  \qquad
  \int_{\cU_d} \|\br\|^2 \bbone_{\{\|\br\|\geq1\}} \, \nu(\dd \br) < \infty
 \end{equation}
 hold.
Then for each \ $\lambda \in \sigma(\tbB)$ \ with
 \ $\Re(\lambda) \in \bigl(\frac{1}{2} s(\tbB), s(\tbB)\bigr]$ \ and for each
 left eigenvector \ $\bv \in \CC^d$ \ of \ $\tbB$ \ corresponding to the eigenvalue
 \ $\lambda$, \ we have
 \begin{equation}\label{convwvL2}
  \ee^{-\lambda t} \langle \bv, \bX_t \rangle \qmean w_{\bv,\bX_0} \qquad
  \text{as \ $t \to \infty$}
 \end{equation}
 (especially, \ $\EE(|w_{\bv,\bX_0}|^2) < \infty$), where \ $w_{\bv,\bX_0}$ \ is
 introduced in \eqref{convwv}, and the improper integrals in \eqref{reprCBIw} are
 convergent in \ $L_2$.

Moreover,
 \begin{equation}\label{convCBIwL2}
  \ee^{-s(\tbB)t} \bX_t \qmean w_{\bu,\bX_0} \tbu \qquad
  \text{as \ $t \to \infty$.}
 \end{equation}
\end{Thm}

\noindent
\textbf{Proof.}
First, note that the moment conditions of Theorem \ref{convCBIasL1} hold, so,
 especially, we have the representation \eqref{reprCBIw}, where the improper
 integrals are convergent almost surely and in \ $L_1$ \ as well.
In order to show \eqref{convwvL2}, we consider the representation
 \[
   w_{\bv,\bX_0} - \ee^{-\lambda t} \langle\bv, \bX_t\rangle
   = D_t^{(1)} + D_t^{(2)} + D_t^{(3)} + D_t^{(4)} + D_t^{(5)} , \qquad
   t \in \RR_+ ,
 \]
 where \ $D_t^{(j)}$, \ $j\in\{1,2,3,4,5\}$, \ $t\in\RR_+$, \ are defined in the proof of Theorem \ref{convCBIasL1}.
Recall that by \eqref{D_1_conv} and \eqref{D_2_conv_L2}, we have
 \ $D_t^{(1)} \to 0$ \ and \ $D_t^{(2)} \qmean 0$ \ as \ $t \to \infty$.
\ By \eqref{EX_bound}, for each \ $t \in \RR_+$ \ and
 \ $\ell \in \{1, \ldots, d\}$, \ we have
 \begin{align*}
  &\EE\biggl(\int_t^\infty \int_{\cU_d} \int_{\cU_1}
                |\ee^{-\lambda u}|^2 |\langle\bv, \bz\rangle|^2
                \bbone_{\{w\leq X_{u,\ell}\}}
                \, \dd u \, \mu_\ell(\dd\bz) \, \dd w\biggr) \\
  &= \int_t^\infty \int_{\cU_d}
      \ee^{-2\Re(\lambda)u} |\langle\bv, \bz\rangle|^2 \EE(X_{u,\ell})
      \, \dd u \, \mu_\ell(\dd\bz)\\
  &\leq C_4 \|\bv\|^2
        \int_t^\infty \int_{\cU_d}
         \ee^{-(2\Re(\lambda)-s(\tbB))u} \|\bz\|^2 \, \dd u \, \mu_\ell(\dd\bz) \\
  &= C_4 \|\bv\|^2
     \int_t^\infty \ee^{-(2\Re(\lambda)-s(\tbB))u} \, \dd u
     \int_{\cU_d} \|\bz\|^2 \, \mu_\ell(\dd\bz) \\
  &= \frac{C_4\|\bv\|^2}{2\Re(\lambda)-s(\tbB)}
     \ee^{-(2\Re(\lambda)-s(\tbB))t}
     \int_{\cU_d} \|\bz\|^2 \, \mu_\ell(\dd\bz)
   < \infty ,
 \end{align*}
 since
 \ $\int_{\cU_d} \|\bz\|^2 \, \mu_\ell(\dd\bz)
    = \int_{\cU_d} \|\bz\|^2 \bbone_{\{\|\bz\|<1\}} \, \mu_\ell(\dd\bz)
      + \int_{\cU_d} \|\bz\|^2 \bbone_{\{\|\bz\|\geq1\}} \, \mu_\ell(\dd\bz)
    < \infty$
 \ by Definition \ref{Def_admissible} and by the moment condition
 \eqref{moment_condition_CBI2}.
Hence, by page 63 in Ikeda and Watanabe \cite{IkeWat}, for each \ $t \in \RR_+$,
 \ we get
 \begin{align*}
  \EE(|D_t^{(3)} + D_t^{(4)}|^2)
  &= \sum_{\ell=1}^d
      \EE\biggl(\int_t^\infty \int_{\cU_d} \int_{\cU_1}
                |\ee^{-\lambda u}|^2 |\langle\bv, \bz\rangle|^2
                \bbone_{\{w\leq X_{u,\ell}\}}
                \, \dd u \, \mu_\ell(\dd\bz) \, \dd w\biggr) \\
  &\leq \frac{C_4\|\bv\|^2}{2\Re(\lambda)-s(\tbB)}
        \ee^{-(2\Re(\lambda)-s(\tbB))t}
         \sum_{\ell=1}^d \int_{\cU_d} \|\bz\|^2 \, \mu_\ell(\dd\bz) ,
 \end{align*}
 yielding
 \begin{equation}\label{D_3-4_conv_L2}
  D_t^{(3)} + D_t^{(4)} \qmean 0 \qquad \text{as \ $t \to \infty$.}
 \end{equation}
Moreover, again by page 63 in Ikeda and Watanabe \cite{IkeWat}, for each
 \ $t \in \RR_+$, \ we have
 \begin{align*}
  \EE(|D_t^{(5)}|^2)
  = \int_t^\infty \int_{\cU_d}
     |\ee^{-\lambda u}|^2 |\langle\bv, \br\rangle|^2
     \, \dd u \, \nu(\dd\br)
  \leq \frac{\|\bv\|^2}{2\Re(\lambda)} \ee^{-2\Re(\lambda)t}
     \int_{\cU_d} \|\br\|^2 \, \nu(\dd\br)
   < \infty
 \end{align*}
 since \ $\Re(\lambda) > 0$ \ and
 \ $\int_{\cU_d} \|\br\|^2 \, \nu(\dd\br)
    = \int_{\cU_d} \|\br\|^2 \bbone_{\{\|\br\|<1\}} \, \nu(\dd\br)
      + \int_{\cU_d} \|\br\|^2 \bbone_{\{\|\br\|\geq1\}} \, \nu(\dd\br)
    < \infty$
 \ by Definition \ref{Def_admissible} and by the moment condition
 \eqref{moment_condition_CBI2}, implying
 \begin{equation}\label{D_5_conv_L2}
  D_t^{(5)} \qmean 0 \qquad \text{as \ $t \to \infty$.}
 \end{equation}
The convergences \eqref{D_1_conv}, \eqref{D_2_conv_L2}, \eqref{D_3-4_conv_L2}
 and \eqref{D_5_conv_L2} yield \eqref{convwvL2}.

In order to show \eqref{convCBIwL2}, using the moment condition
 \eqref{moment_condition_CBI2}, first we prove that for each \ $T \in \RR_+$,
 \ we have
 \begin{equation}\label{Lemma_L2}
  \bDelta_{t,t+T} \qmean \bzero \qquad \text{as \ $t \to \infty$,}
 \end{equation}
 where \ $\bDelta_{t,t+T}$, \ $t, T \in \RR_+$, \ are introduced in
 \eqref{Delta_t_t+T}.
We use the representation \eqref{reprJ} for \ $\bDelta_{t,t+T}$,
 \ $t, T \in \RR_+$.
\ Recall that by \eqref{J_1_conv} and \eqref{J_2_conv_L2}, we have
 \ $\bJ_{t,t+T}^{(1)} \to \bzero$ \ and \ $\bJ_{t,t+T}^{(2)} \qmean \bzero$ \ as
 \ $t \to \infty$ \ for all \ $T \in \RR_+$.
\ By \eqref{C} and \eqref{EX_bound}, for each \ $t, T \in \RR_+$ \ and
 \ $\ell \in \{1, \ldots, d\}$, \ we get
 \begin{align*}
  &\EE\biggl(\int_t^{t+T} \int_{\cU_d} \int_{\cU_1}
              \bigl\|\ee^{(t+T-v)\tbB} \bz\bigr\|^2
              \bbone_{\{w\leq X_{v,\ell}\}}
              \, \dd v \, \mu_\ell(\dd\bz) \, \dd w\biggr) \\
  &= \int_t^{t+T} \int_{\cU_d}
      \bigl\|\ee^{(t+T-v)\tbB} \bz\bigr\|^2 \EE(X_{v,\ell})
      \, \dd v \, \mu_\ell(\dd\bz)\\
  &\leq C_3^2 C_4
        \int_t^{t+T} \int_{\cU_d}
         \ee^{s(\tbB)2(t+T-v)} \|\bz\|^2 \ee^{s(\tbB)v}
         \, \dd v \, \mu_\ell(\dd\bz) \\
  &\leq C_3^2 C_4 \ee^{2s(\tbB)(t+T)}
        \int_t^\infty \int_{\cU_d}
         \ee^{-s(\tbB)v} \|\bz\|^2 \, \dd v \, \mu_\ell(\dd\bz) \\
  &= C_3^2 C_4 \ee^{2s(\tbB)(t+T)}
     \int_t^\infty \ee^{-s(\tbB)v} \, \dd v
     \int_{\cU_d} \|\bz\|^2 \, \mu_\ell(\dd\bz) \\
  &= \frac{C_3^2 C_4}{s(\tbB)} \ee^{s(\tbB)(t+2T)}
     \int_{\cU_d} \|\bz\|^2 \, \mu_\ell(\dd\bz)
   < \infty ,
 \end{align*}
 hence, by page 62 in Ikeda and Watanabe \cite{IkeWat}, we
 have
 \begin{align*}
  \EE(\|\bJ_{t,t+T}^{(3)} + \bJ_{t,t+T}^{(4)}\|^2)
   = \ee^{-2s(\tbB)(t+T)}
     \sum_{\ell=1}^d
      \EE\biggl(\int_t^{t+T} \int_{\cU_d} \int_{\cU_1}
                 \bigl\|\ee^{(t+T-v)\tbB} \bz\bigr\|^2
                 \bbone_{\{w\leq X_{v,\ell}\}}
                 \, \dd v \, \mu_\ell(\dd\bz) \, \dd w\biggr)
  \end{align*}
  \begin{align*}
  &\leq \frac{C_3^2 C_4}{s(\tbB)} \ee^{-s(\tbB)t}
        \sum_{\ell=1}^d
         \int_{\cU_d} \|\bz\|^2 \, \mu_\ell(\dd\bz) .
 \end{align*}
Consequently, for each \ $T \in \RR_+$, \ we conclude
 \begin{equation}\label{J_3-4_conv_L2}
  \bJ_{t,t+T}^{(3)} + \bJ_{t,t+T}^{(4)} \qmean \bzero \qquad
  \text{as \ $t \to \infty$.}
 \end{equation}
Moreover, by page 62 in Ikeda and Watanabe \cite{IkeWat} and \eqref{C}, for each
 \ $t, T \in \RR_+$, \ we have
 \begin{align*}
  &\EE(\|\bJ_{t,t+T}^{(5)}\|^2)
   = \ee^{-2s(\tbB)(t+T)}
     \int_t^{t+T} \int_{\cU_d}
       \bigl\|\ee^{(t+T-v)\tbB} \br\bigr\|^2
       \, \dd v \, \nu(\dd\br) \\
  &\leq C_3^2 \ee^{-2s(\tbB)(t+T)}
        \int_t^{t+T} \int_{\cU_d}
         \ee^{s(\tbB)2(t+T-v)} \|\br\|^2
         \, \dd v \, \nu(\dd\br)
   \leq \frac{C_3^2}{2s(\tbB)} \ee^{-2s(\tbB)t}
     \int_{\cU_d} \|\br\|^2 \, \nu(\dd\br) .
 \end{align*}
Consequently, for each \ $T \in \RR_+$, \ we get
 \begin{equation}\label{J_5_conv_L2}
  \bJ_{t,t+T}^{(5)} \qmean \bzero \qquad \text{as \ $t \to \infty$.}
 \end{equation}
The convergences \eqref{J_1_conv}, \eqref{J_2_conv_L2}, \eqref{J_3-4_conv_L2},
 and \eqref{J_5_conv_L2} yield \eqref{Lemma_L2}.

Finally we prove \eqref{convCBIwL2}.
For each \ $t, T \in \RR_+$ \ and \ $i \in \{1, \ldots, d\}$, \ by
 \eqref{difference}, we obtain
 \begin{align*}
  &\EE(|\ee^{-s(\tbB)(t+T)} \be_i^\top \bX_{t+T}
        - w_{\bu,\bX_0} \be_i^\top \tbu|^2) \\
  &\leq 3 \EE(|\be_i^\top \bDelta_{t,t+T}|^2)
        + 3 \|\tbu\|^2
          \EE(|\ee^{-s(\tbB)t} \langle\bu, \bX_t\rangle - w_{\bu,\bX_0}|^2)
        + 3 C_5^2 \|\tbu\|^2 \ee^{-2C_2T-2s(\tbB)t}
          \EE(\langle\bu, \bX_t\rangle^2) .
 \end{align*}
By \eqref{Lemma_L2} and \eqref{convwvL2} with \ $\lambda = s(\tbB)$ \ and
 \ $\bv = \bu$, \ for each \ $T \in \RR_+$ \ and \ $i \in \{1, \ldots, d\}$, \ we
 obtain
 \begin{align*}
  &\limsup_{t\to\infty}
    \EE(|\ee^{-s(\tbB)t} \be_i^\top \bX_t - w_{\bu,\bX_0} \be_i^\top \tbu|^2)
   = \limsup_{t\to\infty}
      \EE(|\ee^{-s(\tbB)(t+T)} \be_i^\top \bX_{t+T}
           - w_{\bu,\bX_0} \be_i^\top \tbu|^2) \\
  &\leq 3 C_5^2 \|\tbu\|^2 \ee^{-2C_2T}
        \limsup_{t\to\infty}
         \EE((\ee^{-s(\tbB)t} \langle\bu, \bX_t\rangle)^2)
   = 3 C_5^2 \|\tbu\|^2 \ee^{-2C_2T} \EE(w_{\bu,\bX_0}^2) ,
 \end{align*}
 hence, by \ $T \to \infty$, \ we conclude \eqref{convCBIwL2}.
\proofend

\appendix

\vspace*{5mm}

\noindent{\bf\Large Appendix}

\section{On a decomposition of CBI processes}
\label{deco_CBI}

\begin{Lem}\label{decomposition_CBI}
If \ $(\bX_t)_{t\in\RR_+}$ \ is a multi-type CBI process with parameters
 \ $(d, \bc, \Bbeta, \bB, \nu, \bmu)$, \ then for each \ $t, T \in \RR_+$, \ we
 have \ $\bX_{t+T} \distre \bX_t^{(1)} + \bX_t^{(2,T)}$, \ where
 \ $(\bX_s^{(1)})_{s\in\RR_+}$ \ and \ $(\bX_s^{(2,T)})_{s\in\RR_+}$ \ are
 independent multi-type CBI processes with \ $\PP(\bX_0^{(1)} = \bzero) = 1$,
 \ $\bX_0^{(2,T)} \distre \bX_T$, \ and with parameters
 \ $(d, \bc, \Bbeta, \bB, \nu, \bmu)$ \ and \ $(d, \bc, \bzero, \bB, 0, \bmu)$,
 \ respectively.
\end{Lem}

\noindent
\textbf{Proof.}
It is known that \ $\bv(r, \bv(s, \blambda)) = \bv(r + s, \blambda)$ \ for all
 \ $r, s \in \RR_+$ \ and \ $\blambda \in \RR_+^d$, \ see, e.g., Li
 \cite[page 58]{Li}.
By the independence of \ $(\bX_s^{(1)})_{s\in\RR_+}$ \ and
 \ $(\bX_s^{(2,T)})_{s\in\RR_+}$, \ by \eqref{Laplace_transform} and by the law of
 total probability, for each \ $t, T \in \RR_+$ \ and \ $\blambda \in \RR_+^d$,
 \ we have
 \begin{align*}
  &\EE(\ee^{-\langle\blambda,\bX_t^{(1)}+\bX_t^{(2,T)}\rangle})
   = \EE(\ee^{-\langle\blambda,\bX_t^{(1)}\rangle})
     \EE(\ee^{-\langle\blambda,\bX_t^{(2,T)}\rangle}) \\
  &= \exp\biggl\{- \langle\bzero, \bv(t, \blambda)\rangle
                 - \int_0^t \psi(\bv(s, \blambda)) \, \dd s\biggr\}
     \EE(\ee^{-\langle\bX_0^{(2,T)}, \bv(t,\blambda)\rangle}) \\
  &= \exp\biggl\{- \int_0^t \psi(\bv(s, \blambda)) \, \dd s\biggr\}
     \EE(\ee^{-\langle\bX_T, \bv(t,\blambda)\rangle}) \\
  &= \exp\biggl\{- \int_0^t \psi(\bv(s, \blambda)) \, \dd s\biggr\}
     \EE\biggl(\exp\biggl\{- \langle\bX_0, \bv(T,\bv(t, \blambda))\rangle
                           - \int_0^T
                              \psi(\bv(s, \bv(t, \blambda)))
                              \, \dd s\biggr\}\biggr) \\
  &= \EE\biggl(\exp\biggl\{- \langle\bX_0, \bv(t+T, \blambda)\rangle
                           - \int_0^t \psi(\bv(s, \blambda)) \, \dd s
                           - \int_0^T
                              \psi(\bv(s + t, \blambda)) \, \dd s\biggr\}\biggr)\\
  &= \EE\biggl(\exp\biggl\{- \langle\bX_0, \bv(t+T, \blambda)\rangle
                           - \int_0^t \psi(\bv(s, \blambda)) \, \dd s
                           - \int_t^{t+T}
                              \psi(\bv(u, \blambda)) \, \dd u\biggr\}\biggr)\\
 &  = \EE\biggl(\exp\biggl\{- \langle\bX_0, \bv(t+T, \blambda)\rangle
                           - \int_0^{t+T}
                              \psi(\bv(s, \blambda)) \, \dd s\biggr\}\biggr)
   = \EE(\ee^{-\langle\blambda,\bX_{t+T}\rangle}) ,
 \end{align*}
 hence we obtain the assertion.
\proofend

\section{On the second moment of projections of multi-type CBI processes}
\label{second_moment_CBI}

An explicit formula for the second moment of the projection of a multi-type CBI
 process on the left eigenvectors of its branching mean matrix is presented
 together with its asymptotic behavior in the supercritical and irreducible case.

\begin{Pro}\label{second_moment_asymptotics_CBI}
If \ $(\bX_t)_{t\in\RR_+}$ \ is a multi-type CBI process with parameters
 \ $(d, \bc, \Bbeta, \bB, \nu, \bmu)$ \ such that \ $\EE(\|\bX_0\|^2) < \infty$
 \ and the moment condition \eqref{moment_condition_CBI2} holds, then for each
 left eigenvector \ $\bv \in \CC^d$ \ of \ $\tbB$ \ corresponding to an arbitrary
 eigenvalue \ $\lambda \in \sigma(\tbB)$, \ we have
 \[
   \EE(|\langle\bv, \bX_t\rangle|^2)
   = E_{\bv,\lambda}(t) + \sum_{\ell=1}^d C_{\bv,\ell} I_{\lambda,\ell}(t)
     + I_\lambda(t) \int_{\cU_d} |\langle\bv, \br\rangle|^2 \, \nu(\dd\br) ,
   \qquad t \in \RR_+ ,
 \]
 with
 \begin{align*}
  E_{\bv,\lambda}(t)
  &:= \EE\biggl(\biggl|\ee^{\lambda t} \langle\bv, \bX_0\rangle
                       + \langle\bv, \tBbeta\rangle
                         \int_0^t \ee^{\lambda(t-u)} \, \dd u\biggr|^2\biggr) , \\
  I_{\lambda,\ell}(t)
  &:= \int_0^t \ee^{2\Re(\lambda)(t-u)} \EE(X_{u,\ell}) \, \dd u , \qquad
   \ell \in \{1, \ldots, d\} ,
 \end{align*}
 \begin{align*}
  I_\lambda(t)
  &:= \int_0^t \ee^{2\Re(\lambda)(t-u)} \, \dd u , \\
  C_{\bv,\ell}
  &:= 2 |\langle\bv, \be_\ell\rangle|^2 c_\ell
      + \int_{\cU_d} |\langle\bv, \bz\rangle|^2 \, \mu_\ell(\dd\bz) , \qquad
   \ell \in \{1, \ldots, d\} .
 \end{align*}
If, in addition, \ $(\bX_t)_{t\in\RR_+}$ \ is supercritical and irreducible, then
 we have
 \[
   \lim_{t\to\infty} h(t) \EE(|\langle\bv, \bX_t\rangle|^2) = M_2 ,
 \]
 where
 \begin{gather*}
  h(t) := \begin{cases}
           \ee^{-s(\tbB)t}
            & \text{if \ $\Re(\lambda)
                          \in \bigl(-\infty, \frac{1}{2} s(\tbB)\bigr)$,} \\
           t^{-1} \ee^{-s(\tbB)t}
            & \text{if \ $\Re(\lambda) = \frac{1}{2} s(\tbB)$,} \\
           \ee^{-2\Re(\lambda)t}
            & \text{if \ $\Re(\lambda)
                          \in \bigl(\frac{1}{2} s(\tbB), s(\tbB)\bigr]$,}
          \end{cases}
 \end{gather*}
 and
 \begin{gather*}
  M_2 := \begin{cases}
          \frac{1}{s(\tbB)-2\Re(\lambda)}
          \bigl(\langle\bu, \EE(\bX_0)\rangle
                + \frac{\langle\bu, \tBbeta\rangle}{s(\tbB)}\bigr)
          \sum_{\ell=1}^d C_{\bv,\ell} \langle\be_\ell, \tbu\rangle
           & \text{if \ $\Re(\lambda)
                         \in \bigl(-\infty, \frac{1}{2} s(\tbB)\bigr)$,} \\[2mm]
          \bigl(\langle\bu, \EE(\bX_0)\rangle
                + \frac{\langle\bu, \tBbeta\rangle}{s(\tbB)}\bigr)
          \sum_{\ell=1}^d C_{\bv,\ell} \langle\be_\ell, \tbu\rangle
           & \text{if \ $\Re(\lambda) = \frac{1}{2} s(\tbB)$,} \\[2mm]
          \EE\bigl(\bigl|\langle\bv, \bX_0\rangle
                         + \frac{\langle\bv, \tBbeta\rangle}{\lambda}\bigr|^2\bigr)
          + \frac{1}{2\Re(\lambda)}
            \int_{\cU_d} |\langle\bv, \br\rangle|^2 \, \nu(\dd\br) \\
          + \sum\limits_{\ell=1}^d
             C_{\bv,\ell}
             \be_\ell^\top (2 \Re(\lambda) \bI_d - \tbB)^{-1}
             \bigl(\EE(\bX_0) + \frac{\tBbeta}{2\Re(\lambda)}\bigr)
           & \text{if \ $\Re(\lambda)
                         \in \bigl(\frac{1}{2} s(\tbB), s(\tbB)\bigr]$.}
         \end{cases}
 \end{gather*}
\end{Pro}

Note that Proposition \ref{second_moment_asymptotics_CBI} can be considered as a
 counterpart of Theorem 1 in Section 7 in Chapter V in Athreya and Ney
 \cite{AthNey}.

\noindent
\textbf{Proof of Proposition \ref{second_moment_asymptotics_CBI}.}
We use the representation \eqref{reprZ} for
 \ $\ee^{-\lambda t} \langle\bv, \bX_t\rangle$, \ $t \in \RR_+$, \ where
 \ $Z_t^{(1)}$ \ is deterministic for all \ $t \in \RR_+$.
\ By the independence of \ $\bX_0$, \ $(W_{t,1})_{t\geq0}$, \ldots,
 \ $(W_{t,d})_{t\geq0}$, \ $N_1$, \ldots, \ $N_d$ \ and \ $M$, \ for each
 \ $t \in \RR_+$, \ the random variables \ $\langle\bv, \bX_0\rangle$ \ and
 \ $Z_t^{(j)}$, \ $j \in \{2, 3, 4, 5\}$ \ are conditionally independent with
 respect to \ $(\bX_u)_{u\in[0,t]}$, \ hence for each \ $t \in \RR_+$,
 \[
   \EE\bigl(|\ee^{-\lambda t} \langle\bv, \bX_t\rangle|^2
            \,\big|\, (\bX_u)_{u\in[0,t]}\bigr)
   = \EE\bigl(\bigl|\langle\bv, \bX_0\rangle + Z_t^{(1)}\bigr|^2
              \,\big|\, (\bX_u)_{u\in[0,t]}\bigr)
     + \sum_{j=2}^5
        \EE\bigl(\bigl|Z_t^{(j)}\bigr|^2 \,\big|\, (\bX_u)_{u\in[0,t]}\bigr)
 \]
 almost surely.
We have, for each \ $t \in \RR_+$,
 \begin{align*}
  \EE\bigl(\bigl|\langle\bv, \bX_0\rangle + Z_t^{(1)}\bigr|^2
           \,\big|\, (\bX_u)_{u\in[0,t]}\bigr)
  &= \biggl|\langle\bv, \bX_0\rangle
            + \langle\bv, \tBbeta\rangle
              \int_0^t \ee^{-\lambda u} \, \dd u\biggr|^2 , \\
  \EE\bigl(\bigl|Z_t^{(2)}\bigr|^2 \,\big|\, (\bX_u)_{u\in[0,t]})
  &= 2 \sum_{\ell=1}^d
        |\langle\bv, \be_\ell\rangle|^2 c_\ell
        \int_0^t \ee^{-2\Re(\lambda)u} X_{u,\ell} \, \dd u , \\
  \EE\bigl(\bigl|Z_t^{(3)}\bigr|^2 + \bigl|Z_t^{(4)}\bigr|^2
           \,\big|\, (\bX_u)_{u\in[0,t]}\bigr)
  &= \sum_{\ell=1}^d
      \int_0^t \ee^{-2\Re(\lambda)u} X_{u,\ell} \, \dd u
      \int_{\cU_d} |\langle\bv, \bz\rangle|^2 \, \mu_\ell(\dd\bz) , \\
  \EE\bigl(\bigl|Z_t^{(5)}\bigr|^2 \,\big|\, (\bX_u)_{u\in[0,t]}\bigr)
  &= \int_0^t \ee^{-2\Re(\lambda)u} \, \dd u
     \int_{\cU_d} |\langle\bv, \br\rangle|^2 \, \nu(\dd\br)
 \end{align*}
 almost surely.
Taking the expectation and multiplying by
 \ $|\ee^{\lambda t}|^2 = \ee^{2\Re(\lambda)t}$, \ $t \in \RR_+$, \ we obtain the
 formula for \ $\EE(|\langle\bv, \bX_t\rangle|^2)$, \ $t \in \RR_+$.

Assume now that, in addition, \ $(\bX_t)_{t\in\RR_+}$ \ is supercritical and
 irreducible.
For each \ $t \in \RR_+$, \ we have
 \[
   E_{\bv,\lambda}(t)
   \leq 2 \ee^{2\Re(\lambda)t} \EE(|\langle\bv, \bX_0\rangle|^2)
        + 2 |\langle\bv, \tBbeta\rangle|^2
          \biggl(\int_0^t \ee^{\Re(\lambda)w} \, \dd w\biggr)^2 .
 \]
If \ $\Re(\lambda) \in (-\infty, 0)$, \ then we have
 \[
   \ee^{-s(\tbB)t} E_{\bv,\lambda}(t)
   \leq 2 \ee^{-(s(\tbB)-2\Re(\lambda))t} \EE(|\langle\bv, \bX_0\rangle|^2)
        + \frac{2|\langle\bv, \tBbeta\rangle|^2}{\Re(\lambda)^2} \ee^{-s(\tbB)t}
          (\ee^{\Re(\lambda)t} - 1)^2
   \to 0
 \]
 as \ $t \to \infty$.
\ If \ $\Re(\lambda) = 0$, \ then we have
 \[
   \ee^{-s(\tbB)t} E_{\bv,\lambda}(t)
   \leq 2 \ee^{-s(\tbB)t} \EE(|\langle\bv, \bX_0\rangle|^2)
        + 2 |\langle\bv, \tBbeta\rangle|^2 t^2 \ee^{-s(\tbB)t}
   \to 0 \qquad \text{as \ $t \to \infty$.}
 \]
If \ $\Re(\lambda) \in \bigl(0, \frac{1}{2} s(\tbB)\bigr)$, \ then we have
 \[
   \ee^{-s(\tbB)t} E_{\bv,\lambda}(t)
   \leq 2 \ee^{-(s(\tbB)-2\Re(\lambda))t} \EE(|\langle\bv, \bX_0\rangle|^2)
        + \frac{2|\langle\bv, \tBbeta\rangle|^2}{\Re(\lambda)^2}
          \ee^{-(s(\tbB)-2\Re(\lambda))t}
          (1 - \ee^{-\Re(\lambda)t})^2
   \to 0
 \]
 as \ $t \to \infty$.
\ If \ $\Re(\lambda) = \frac{1}{2} s(\tbB)$ \ then
 \[
   t^{-1} \ee^{-s(\tbB)t} E_{\bv,\lambda}(t)
   \leq 2 t^{-1} \EE(|\langle\bv, \bX_0\rangle|^2)
        + \frac{2|\langle\bv, \tBbeta\rangle|^2}{\Re(\lambda)^2} t^{-1}
          (1 - \ee^{-\Re(\lambda)t})^2
   \to 0
 \]
 as \ $t \to \infty$.
\ If \ $\Re(\lambda) \in \bigl(\frac{1}{2} s(\tbB), s(\tbB)\bigr]$, \ then, by the
 dominated convergence theorem, we have
 \begin{align*}
  \ee^{-2\Re(\lambda)t} E_{\bv,\lambda}(t)
  &= \EE\biggl(\biggl|\langle\bv, \bX_0\rangle
                      + \langle\bv, \tBbeta\rangle
                        \int_0^t \ee^{-\lambda u} \, \dd u\biggr|^2\biggr) \\
  &\to \EE\biggl(\biggl|\langle\bv, \bX_0\rangle
                        + \langle\bv, \tBbeta\rangle
                          \int_0^\infty \ee^{-\lambda u} \, \dd u\biggr|^2\biggr)
     = \EE\biggl(\biggl|\langle\bv, \bX_0\rangle
                        + \frac{\langle\bv, \tBbeta\rangle}
                               {\lambda}\biggr|^2\biggr)
 \end{align*}
 as \ $t \to \infty$.

Moreover, for each \ $u \in \RR_+$ \ and \ $\ell \in \{1, \ldots, d\}$, \ by
 formula \eqref{EXcond}, we have
 \[
   \EE(X_{u,\ell})
   = \be_\ell^\top \ee^{u\tbB} \EE(\bX_0)
     + \int_0^u \be_\ell^\top \ee^{w\tbB} \tBbeta \, \dd w ,
 \]
 thus, for each \ $t \in \RR_+$ \ and \ $\ell \in \{1, \ldots, d\}$, \ we get
 \[
   I_{\lambda,\ell}(t)
   = \be_\ell^\top \bA_{\lambda,1}(t) \EE(\bX_0)
     + \be_\ell^\top \bA_{\lambda,2}(t) \tBbeta
 \]
 with
 \[
   \bA_{\lambda,1}(t)
   := \int_0^t \ee^{2\Re(\lambda)(t-u)} \ee^{u\tbB} \, \dd u , \qquad
   \bA_{\lambda,2}(t)
   := \int_0^t
       \ee^{2\Re(\lambda)(t-u)} \biggl(\int_0^u \ee^{w\tbB} \, \dd w\biggr) \dd u .
 \]
If \ $\Re(\lambda) \in \bigl(-\infty, \frac{1}{2} s(\tbB)\bigr)$, \ then we have
 \[
   \ee^{-s(\tbB)t} \bA_{\lambda,1}(t)
   = \ee^{-(s(\tbB)-2\Re(\lambda))t}
     \int_0^t \ee^{u(\tbB-2\Re(\lambda)\bI_d)} \, \dd u
   \to \frac{1}{s(\tbB)-2\Re(\lambda)} \tbu \bu^\top \qquad
   \text{as \ $t \to \infty$,}
 \]
 since \ $\tbB - 2 \Re(\lambda) \bI_d \in \RR_{(+)}^{d\times d}$ \ is irreducible,
 \ $s(\tbB - 2 \Re(\lambda) \bI_d) = s(\tbB) - 2 \Re(\lambda) > 0$, \ and the left
 (right) Perron vectors of the matrices \ $\ee^{\tbB}$ \ and
 \ $\ee^{\tbB - 2 \Re(\lambda) \bI_d}$ \ coincide, \ see, e.g., the proof of formula (B.2) in Barczy et al.\ \cite{BarPap}.
If \ $\Re(\lambda) = 0$, \ then, by Fubini's theorem, we obtain
 \[
   \ee^{-s(\tbB)t} \bA_{\lambda,2}(t)
   = \ee^{-s(\tbB)t} \int_0^t (t - w) \ee^{w\tbB} \, \dd w
   \to \frac{1}{s(\tbB)^2} \tbu \bu^\top \qquad
   \text{as \ $t \to \infty$,}
 \]
 since, by \eqref{C}, we have
 \begin{align*}
  &\biggl\|\ee^{-s(\tbB)t} \int_0^t (t - w) \ee^{w\tbB} \, \dd w
           - \frac{1}{s(\tbB)^2} \tbu \bu^\top\biggr\| \\
  &\leq \biggl\|\ee^{-s(\tbB)t}
               \int_0^t (t - w)
                \ee^{s(\tbB)w}
                (\ee^{-s(\tbB)w} \ee^{w\tbB} - \tbu \bu^\top) \, \dd w\biggr\| \\
  &\quad
       + \biggl\|\ee^{-s(\tbB)t}
                 \int_0^t (t - w) \ee^{s(\tbB)w} \tbu \bu^\top \, \dd w
                 - \frac{1}{s(\tbB)^2} \tbu \bu^\top\biggr\| \\
  &\leq C_1 \ee^{-s(\tbB)t}
       \int_0^t (t - w) \ee^{s(\tbB)w} \ee^{-C_2w} \, \dd w
       + \biggl\|\frac{\ee^{-s(\tbB)t}}{s(\tbB)^2}
                 (\ee^{s(\tbB)t} - s(\tbB) t - 1) \tbu \bu^\top
                 - \frac{1}{s(\tbB)^2} \tbu \bu^\top\biggr\| \\
  &\leq C_1 \ee^{-s(\tbB)t}
       \int_0^t (t - w) \ee^{s(\tbB)w} \ee^{-\tC_2 w} \, \dd w
       + \biggl\|\frac{\ee^{-s(\tbB)t}}{s(\tbB)^2}
                 (s(\tbB) t + 1) \tbu \bu^\top\biggr\| \\
  &= \frac{C_1\ee^{-s(\tbB)t}}{(s(\tbB)-\tC_2)^2}
        \bigl(\ee^{(s(\tbB)-\tC_2)t} - (s(\tbB)-\tC_2) t - 1\bigr)
       + \frac{\ee^{-s(\tbB)t}}{s(\tbB)^2} (s(\tbB) t + 1) \|\tbu \bu^\top\|
  \to 0
 \end{align*}
 as \ $t \to \infty$, \ where \ $\tC_2 \in (0, C_2 \land s(\tbB))$.
\ Hence, if \ $\Re(\lambda) = 0$, \ then
 \[
   \ee^{-s(\tbB)t} I_{\lambda,\ell}(t)
   \to \frac{1}{s(\tbB)} \be_\ell^\top \tbu \bu^\top
       \biggl(\EE(\bX_0) + \frac{\tBbeta}{s(\tbB)}\biggr)
   = \frac{\langle\be_\ell, \tbu\rangle}{s(\tbB)}
     \biggl(\langle\bu, \EE(\bX_0)\rangle
            + \frac{\langle\bu, \tBbeta\rangle}{s(\tbB)}\biggr)
 \]
 as \ $t \to \infty$.
\ If \ $\Re(\lambda) \in (-\infty, \frac{1}{2} s(\tbB)) \setminus \{0\}$, \ then,
 by Fubini's theorem, we obtain
 \begin{align*}
  \ee^{-s(\tbB)t} \bA_{\lambda,2}(t)
  &= \ee^{-(s(\tbB)-2\Re(\lambda))t}
     \int_0^t
      \ee^{-2\Re(\lambda)u} \biggl(\int_0^u \ee^{w\tbB} \, \dd w\biggr) \dd u \\
  &= \ee^{-(s(\tbB)-2\Re(\lambda))t}
     \int_0^t
      \biggl(\int_w^t \ee^{-2\Re(\lambda)u} \, \dd u\biggr)
      \ee^{w\tbB} \, \dd w \\
  &= \frac{1}{2\Re(\lambda)} \ee^{-(s(\tbB)-2\Re(\lambda))t}
     \int_0^t
      \biggl(\ee^{-2\Re(\lambda)w} - \ee^{-2\Re(\lambda)t}\biggr) \ee^{w\tbB}
      \, \dd w
 \end{align*}
 \begin{align*}
  &= \frac{1}{2\Re(\lambda)}
     \biggl(\ee^{-(s(\tbB)-2\Re(\lambda))t}
            \int_0^t \ee^{(\tbB-2\Re(\lambda)\bI_d)w} \, \dd w
            - \ee^{-s(\tbB)t} \int_0^t \ee^{w\tbB} \, \dd w\biggr) \\
  &\to \frac{1}{2\Re(\lambda)}
       \biggl(\frac{\tbu \bu^\top}{s(\tbB)-2\Re(\lambda)}
              - \frac{\tbu \bu^\top}{s(\tbB)}\biggr)
   = \frac{\tbu \bu^\top}{(s(\tbB)-2\Re(\lambda))s(\tbB)}
 \end{align*}
 as \ $t \to \infty$, \ since \ $s(\tbB)-2\Re(\lambda) > 0$ \ and \ $s(\tbB) > 0$
 \ imply
 \[
   \ee^{-(s(\tbB)-2\Re(\lambda))t}
   \int_0^t \ee^{(\tbB-2\Re(\lambda)\bI_d)w} \, \dd w
   \to \frac{\tbu \bu^\top}{s(\tbB)-2\Re(\lambda)} , \qquad
   \ee^{-s(\tbB)t} \int_0^t \ee^{w\tbB} \, \dd w
   \to \frac{\tbu \bu^\top}{s(\tbB)}
 \]
 as \ $t \to \infty$, \ see, e.g., the proof of Proposition B.1 in
 Barczy et al.\ \cite{BarPap}.
Hence, if \ $\Re(\lambda) \in (-\infty, \frac{1}{2} s(\tbB)) \setminus \{0\}$,
 \ then
 \begin{align*}
  \ee^{-s(\tbB)t} I_{\lambda,\ell}(t)
  &\to \frac{1}{s(\tbB)-2\Re(\lambda)} \be_\ell^\top \tbu \bu^\top
       \biggl(\EE(\bX_0) + \frac{\tBbeta}{s(\tbB)}\biggr) \\
  &= \frac{\langle\be_\ell, \tbu\rangle}{s(\tbB)-2\Re(\lambda)}
     \biggl(\langle\bu, \EE(\bX_0)\rangle
            + \frac{\langle\bu, \tBbeta\rangle}{s(\tbB)}\biggr) \qquad
  \text{as \ $t \to \infty$.}
 \end{align*}
If \ $\Re(\lambda) = \frac{1}{2} s(\tbB)$, \ then we have
 \[
   t^{-1} \ee^{-s(\tbB)t} \bA_{\lambda,1}(t)
   = \frac{1}{t} \int_0^t \ee^{-s(\tbB)u} \ee^{u\tbB} \, \dd u
   \to \tbu \bu^\top \qquad
   \text{as \ $t \to \infty$,}
 \]
 see, e.g., part (v) of Lemma A.2 in Barczy et al.\ \cite{BarPap}.
Moreover, if \ $\Re(\lambda) = \frac{1}{2} s(\tbB)$, \ then, by Fubini's theorem,
 we obtain
 \begin{align*}
  t^{-1} \ee^{-s(\tbB)t} \bA_{\lambda,2}(t)
  &= \frac{1}{t}
     \int_0^t \ee^{-s(\tbB)u} \biggl(\int_0^u \ee^{w\tbB} \, \dd w\biggr) \dd u
   = \frac{1}{t}
     \int_0^t
      \biggl(\int_w^t \ee^{-s(\tbB)u} \, \dd u\biggr)
      \ee^{w\tbB} \, \dd w \\
  &= \frac{1}{s(\tbB)t}
     \int_0^t
      \biggl(\ee^{-s(\tbB)w} - \ee^{-s(\tbB)t}\biggr) \ee^{w\tbB}
      \, \dd w \\
  &= \frac{1}{s(\tbB)}
     \biggl(\frac{1}{t} \int_0^t \ee^{-s(\tbB)w} \ee^{w\tbB} \, \dd w
            - \frac{1}{t} \ee^{-s(\tbB)t} \int_0^t \ee^{w\tbB} \, \dd w\biggr)
   \to \frac{1}{s(\tbB)} \tbu \bu^\top
 \end{align*}
 as \ $t \to \infty$, \ since
 	\begin{equation}\label{Proposition B.1}
 	\ee^{-s(\tbB)t} \int_0^t \ee^{w\tbB} \, \dd w
 	\to \frac{\tbu\bu^\top}{s(\tbB)} \qquad \text{as \ $t \to \infty$,}
 	\end{equation}
 	see, e.g., Barczy et al.\ \cite[formula (B.2)]{BarPap}.
\ Hence, if \ $\Re(\lambda) = \frac{1}{2} s(\tbB)$, \ then
 \[
   t^{-1} \ee^{-s(\tbB)t} I_{\lambda,\ell}(t)
   \to \be_\ell^\top \tbu \bu^\top
       \biggl(\EE(\bX_0) + \frac{\tBbeta}{s(\tbB)}\biggr)
   = \langle\be_\ell, \tbu\rangle
     \biggl(\langle\bu, \EE(\bX_0)\rangle
            + \frac{\langle\bu, \tBbeta\rangle}{s(\tbB)}\biggr)
 \]
 as \ $t \to \infty$.
\ If \ $\Re(\lambda) \in \bigl(\frac{1}{2} s(\tbB), s(\tbB)\bigr]$, \ then we have
 \[
   \ee^{-2\Re(\lambda)t} \bA_{\lambda,1}(t)
   = \int_0^t \ee^{u(\tbB - 2\Re(\lambda)\bI_d)} \, \dd u
   \to \int_0^\infty \ee^{u(\tbB - 2\Re(\lambda)\bI_d)} \, \dd u
   = (2\Re(\lambda)\bI_d-\tbB)^{-1}
 \]
 as \ $t \to \infty$, \ since
 \ $\tbB -  2 \Re(\lambda) \bI_d \in \RR_{(+)}^{d\times d}$ \ is irreducible and
 \ $s(\tbB - 2 \Re(\lambda) \bI_d) = s(\tbB) - 2 \Re(\lambda) < 0$, \ see, e.g.,
 the proof of Proposition B.1 in Barczy et al.\ \cite{BarPap}.
Moreover, if \ $\Re(\lambda) \in \bigl(\frac{1}{2} s(\tbB), s(\tbB)\bigr]$, \ then,
 by Fubini's theorem and \eqref{Proposition B.1}, we obtain
 \begin{align*}
  \ee^{-2\Re(\lambda)t} \bA_{\lambda,2}(t)
  &= \int_0^t
      \ee^{-2\Re(\lambda)u} \biggl(\int_0^u \ee^{w\tbB} \, \dd w\biggr) \dd u
  = \int_0^t
      \biggl(\int_w^t \ee^{-2\Re(\lambda)u} \, \dd u\biggr) \ee^{w\tbB} \, \dd w \\
  &= \frac{1}{2\Re(\lambda)}
     \int_0^t
      (\ee^{-2\Re(\lambda)w} - \ee^{-2\Re(\lambda)t}) \ee^{w\tbB} \, \dd w
  \end{align*}
 \begin{align*}
  &= \frac{1}{2\Re(\lambda)} \int_0^t \ee^{w( \tbB - 2\Re(\lambda)\bI_d)} \, \dd w
     - \ee^{-(2\Re(\lambda)-s(\tbB))t} \ee^{-s(\tbB)t}
       \int_0^t \ee^{w\tbB} \, \dd w \\
  &\to \frac{1}{2\Re(\lambda)} (2\Re(\lambda)\bI_d-\tbB)^{-1} \qquad
   \text{as \ $t \to \infty$.}
 \end{align*}
Hence, if \ $\Re(\lambda) \in \bigl(\frac{1}{2} s(\tbB), s(\tbB)\bigr]$, \ then
 \[
   \ee^{-2\Re(\lambda)t} I_{\lambda,\ell}(t)
   \to \be_\ell^\top (2\Re(\lambda)\bI_d-\tbB)^{-1}
       \biggl(\EE(\bX_0) + \frac{\tBbeta}{2\Re(\lambda)}\biggr) \qquad
   \text{as \ $t \to \infty$.}
 \]

Further, we have
 \[
   I_\lambda(t)
   = \int_0^t \ee^{2\Re(\lambda)w} \, \dd w
   = \begin{cases}
      t & \text{if \ $\Re(\lambda) = 0$,} \\
      \frac{1}{2\Re(\lambda)} (\ee^{2\Re(\lambda)t} - 1)
       & \text{if \ $\Re(\lambda) \in (-\infty, s(\tbB)] \setminus \{0\}$.}
     \end{cases}
 \]
If \ $\Re(\lambda) = 0$, \ then
 \[
   \ee^{-s(\tbB)t} I_\lambda(t) = t \ee^{-s(\tbB)t} \to 0 \qquad
   \text{as \ $t \to \infty$.}
 \]
If \ $\Re(\lambda) \in \bigl(-\infty, \frac{1}{2} s(\tbB)\bigr) \setminus \{0\}$,
 \ then
 \[
   \ee^{-s(\tbB)t} I_\lambda(t)
   = \frac{1}{2\Re(\lambda)} (\ee^{-(s(\tbB)-2\Re(\lambda))t} - \ee^{-s(\tbB)t})
   \to 0 \qquad \text{as \ $t \to \infty$.}
 \]
If \ $\Re(\lambda) = \frac{1}{2} s(\tbB)$, \ then \ $\Re(\lambda) \ne 0$ \ and
 \[
   t^{-1} \ee^{-s(\tbB)t} I_\lambda(t)
   = \frac{1}{s(\tbB)} t^{-1} (1 - \ee^{-s(\tbB)t})
   \to 0 \qquad \text{as \ $t \to \infty$.}
 \]
If \ $\Re(\lambda) \in \bigl(\frac{1}{2} s(\tbB), s(\tbB)\bigr]$, \ then
 \ $\Re(\lambda) \in \RR_{++}$ \ and
 \[
   \ee^{-2\Re(\lambda)t} I_\lambda(t)
   = \frac{1}{2\Re(\lambda)} (1 - \ee^{-2\Re(\lambda)t})
   \to \frac{1}{2\Re(\lambda)} \qquad \text{as \ $t \to \infty$.}
 \]
The proof is complete.
\proofend

\section*{Acknowledgements}
We would like to thank the referees for their comments that helped us to improve the paper.

\end{document}